%% file: biot-cam-paper-kknr.tex
\documentclass{elsarticle}
\input{customized_packages_and_macros}

\journal{**}

\definecolor{ao(english)}{rgb}{0.0, 0.5, 0.0}

\definecolor{gainsboro}{rgb}{0.86, 0.86, 0.86}
\definecolor{lightgray}{rgb}{0.83, 0.83, 0.83} 
\definecolor{silver}{rgb}{0.75, 0.75, 0.75}
\definecolor{ashgrey}{rgb}{0.7, 0.75, 0.71}
\definecolor{battleshipgrey}{rgb}{0.52, 0.52, 0.51}
\definecolor{ao(english)}{rgb}{0.0, 0.5, 0.0}

\def \distr#1 {
\begin{tikzpicture}
  	\node[anchor=south west,inner sep=0] (img) at (0,0) {{#1}};
	\node[below=of img, node distance=0cm,  yshift=1.05cm, xshift=0cm, font=\color{gray}] {{\scriptsize number of cells}};
  	\node[left=of img, node distance=0cm, rotate=90, anchor=center, yshift=-0.8cm,font=\color{gray}] {{\scriptsize local error and indicator}};
\end{tikzpicture}
}

\newcommand{\sveta}[1]{{\color{black}{#1}}}
\newcommand{\authors}[1]{{\color{black}{#1}}}

\usepackage{scalerel}
\usepackage{stackengine,wasysym}

\newcommand\reallywidetilde[1]{\ThisStyle{%
\setbox0=\hbox{$\SavedStyle#1$}%
\stackengine{-.1\LMpt}{$\SavedStyle#1$}{%
\stretchto{\scaleto{\SavedStyle\mkern.2mu\AC}{.5150\wd0}}{.3\ht0}%
}{O}{c}{F}{T}{S}%
}}

\begin{document}

\begin{frontmatter}

\title{Guaranteed and computable error bounds \\
for approximations constructed by an iterative decoupling of the Biot problem}

\author{Kundan Kumar$^*$, 
Svetlana Kyas$^{\dagger}$, 
Jan Martin Nordbotten$^{*}$, and 
Sergey Repin $^{\S}$}

\address{$^*$Department of Mathematics, University of Bergen, Norway; \\
$^{\dagger}$Geothermal Energy and Geofluids, Institute of Geophysics, ETH Z\"urich, Switzerland\\
$^{\S}$ University of Jyvaskyla, Finland; 
Steklov Inst. Math. of the RAS at St. Petersburg; 
Peter the Great Polytechnic University, St. Petersburg, Russia.
}

\begin{abstract}
The paper is concerned with guaranteed 
a posteriori error estimates for a class of evolutionary 
problems related to poroelastic media governed by the quasi-static linear Biot equations. 
The system is decoupled employing the fixed-stress split scheme, which leads to a semi-discrete system
solved iteratively. The error bounds are derived by combining a posteriori estimates for
contractive mappings with those of the functional type for elliptic partial differential equations.
The estimates are applicable for any approximation 
in the admissible functional space and are independent of the discretization method. They are fully computable, 
do not contain mesh dependent constants, and provide reliable global estimates of the error measured in the 
energy norm. Moreover, they suggest efficient error indicators for the distribution 
of local errors, which can be used in adaptive procedures. 
%
\end{abstract}

\begin{keyword}
Biot problem,
fixed-stress split iterative scheme,
a posteriori error estimates,
contraction mappings.
\end{keyword}

\end{frontmatter}

\section{Introduction}
\label{sec:intro}
%
Problems defined in a poroelastic medium contribute to a wide range of application areas including simulation 
of oil reservoirs, 
prediction of the environmental changes, soil subsidence and liquefaction in the earthquake engineering, 
well stability, sand production, waste deposition, hydraulic fracturing, 
${\rm CO}_2$ sequestration, and understanding of the biological tissues in biomechanics. 
%
In recent years, mathematical modeling of poroelastic problems has become a highly important topic
because it helps engineers to understand and predict complicated phenomena arising in such media as well 
as assists in preventing possible future financial calamities. 
However, numerical schemes designed for any of the existing models provide approximations that contain 
errors of different nature, which must be controlled. Therefore, a reliable quantitative analysis of poroelasticity 
problems requires efficient and computable error estimates that could be applied for various approximations 
and computation methods.
%

Mathematical modeling of poroelasticity is usually based on the Biot model that consists of the elasticity 
equation coupled with the equation describing the slow motion of the fluid. Computational errors in one part of the 
model may seriously affect the accuracy of the other one. 
%
Therefore, getting reliable and efficient a posteriori error estimates for coupled problems is, in general, a much
more complicated task than for a single equation.

The Biot model is a system describing the flow and the displacement in a porous medium by momentum and 
mass conservation equations. Initially, it was derived at a macroscopic scale (with inertia effects negligible) 
in the works Terzaghi \cite{Terzaghi1925} and Biot \cite{Biot1941}. The settlement of different types of 
soils was predicted in \cite{Terzaghi1925}, which later was extended to the generalized theory of consolidation 
\cite{Biot1941, Biot1955}.  A comprehensive discussion of the theory of poromechanics can be found in 
\cite{Coussy1995}. Thus, for modeling of the {\em solid displacement} $\vectoru$ and the 
{\em fluid pressure} $p$, we consider the system that governs the coupling of an {\em elastic isotropic 
porous medium} saturated with {\em slightly compressible viscous single-phase fluid}
\begin{equation}
\begin{array}{rl}
-{\rm div} \big(\,\lambda \, (\dvrg \vectoru) \, \tensorI 
+ 2 \,\mu \,\strain (\vectoru) 
- \alpha \, p \, \tensorI \,\big) 
\!\! &= \boldsymbol{f} \quad \mbox{in} \quad Q := \Omega \times (0, T), \\
\partial_t \big(\beta\, p + \alpha \dvrg \vectoru \big) 
- \dvrg \tensorK \nabla p
\!\! &= g \, \quad \mbox{in} \quad Q,
\end{array}
\label{eq:poroelastic-system}
\end{equation}
where $Q$ denotes a space-time cylinder (with bounded domain $\Omega \subset \mathbb{R}^{d}$, ${d} = {\{2, 3\}}$ 
having Lipschitz continuous boundary $\partial \Omega$ and given time-interval $(0, T)$, $0 < T < +\infty$), 
$\boldsymbol{f} \in H^1(0, T; [\L{2}(\Omega)]^d)$ and $g \in L^2(0, T; L^2(\Omega))$ are the body 
force and the volumetric fluid source, respectively\footnote{For convenience of the reader, we collected the definitions related to the functional spaces in Appendix}. 
The first equation in \eqref{eq:poroelastic-system} follows from the balance of linear momentum for the 
{\em total Cauchy stress tensor} 
$$\stress_{\rm por} := \stress(\vectoru) - \alpha \, p \, \tensorI$$ that accounts not only $\vectoru$ but also 
the pressure $p$ scaled by the dimensionless Biot-Willis coefficient $\alpha > 0$. Linear elastic tensor is governed by 
Hooke's law $$\stress (\vectoru) 
:= 2 \,\mu \,\strain (\vectoru) + \lambda \, {\rm tr} \strain (\vectoru) \, \tensorI
= 2 \,\mu \,\strain (\vectoru) + \lambda \, (\dvrg \vectoru) \, \tensorI,$$
where $\strain (\vectoru) := \tfrac{1}{2} \big(\nabla\,\vectoru + (\nabla\,\vectoru)^{\rm T} \big)$ is the 
{\em tensor of small strains}. Here,  $\lambda, \mu > 0$ are the {\em Lam\'e constants}
proportional to Young's modulus $E$ and Poisson's ratio $\nu$ via relations $\mu = \tfrac{E}{2\, (1+\nu)}$
and $\lambda = \tfrac{E\, \nu}{(1 + \nu) \, (1 - 2\nu)}$.
The second equation is the fluid mass conservation (continuity) equation in $Q$. Here,  $\beta$ stands for the 
{\em storage coefficient} and 
\sveta{$\tensorK$ is the {\em permeability tensor} 
assumed to be symmetric, uniformly bounded, anisotropic, 
and heterogeneous in space and constant in 
time, i.e., 
\begin{equation}
\lambda_{\tensorK} |\vectortau|^2 \leq \tensorK(x)  \vectortau \cdot \vectortau \leq \mu_{\tensorK} |\vectortau|^2, \quad 
\lambda_{\tensorK}, \mu_{\tensorK} >0, \quad \text{ for all }
\vectortau \in \mathds{R}^{d}.
\label{eq:eigenvalue-K}
\end{equation}
}

Let $\Sigma = \partial \Omega \times (0, T)$ be a lateral surface of $Q$, whereas 
$\Sigma_0 := \partial \Omega \times \{0\}$ and $\Sigma_T := \partial \Omega \times \{T\}$ define the bottom 
and the top parts of the mantel, such that $\partial Q = \Sigma \cup \Sigma_0 \cup \Sigma_T$.
The initial conditions are assumed to be as follows:
\begin{equation}
p(x, 0) = p^\circ \in H^1(\Omega)\quad \mbox{and} \quad \vectoru(x, 0) = \vectoru^\circ  \in [\H{1}(\Omega)]^d 
\quad \;\mbox{on} \quad \Sigma_0.
\label{eq:ic}
\end{equation}
%
%
We introduce the following partitions of the boundary: 
$\partial \Omega 
= \Sigma^{p}_D \cup \Sigma^{p}_N 
= \Sigma^{\vectoru}_D \cup  \Sigma^{\vectoru}_N,$ 
where $\Sigma^{p}_D$ and $\Sigma^{\vectoru}_D$ must have positive measures, i.e., 
$|\Sigma^{p}_D|, |\Sigma^{\vectoru}_D| >0$, with corresponding boundary conditions (BCs):
\begin{equation}
\begin{alignedat}{2}
p & = p_D \quad \mbox{on} \quad \Sigma^{p}_D, 
\\
- \tensorK \nabla p \cdot \vectorn 
& = z_N \quad \mbox{on} \quad \Sigma^{p}_N, 
\\
\vectoru & = \vectoru_D \quad \mbox{on} \quad \Sigma^{\vectoru}_D, 
\\
\stress_{\rm por} \cdot \vectorn & = {\vectort}_N \quad \mbox{on} \quad \Sigma^{\vectoru}_N.
%
\end{alignedat}
\label{eq:bcs}
\end{equation}
%
For the fluid content $\beta \, p + \alpha \, \dvrg \vectoru$, we prescribe the following initial condition
\begin{equation*}
\eta(x, 0) 
:= \beta \, p(x, 0) + \alpha \, \dvrg \vectoru(x, 0) 
= \beta \, p^\circ + \alpha \, \dvrg \vectoru^\circ,
\label{eq:initial-condition}
\end{equation*}
where $p^\circ$ and $\vectoru^\circ$ are defined in \eqref{eq:ic}. To simplify the exposition, 
we consider only homogeneous BCs, i.e., $p_D, z_N = 0$ and $\vectoru_D, {\vectort}_N = \boldsymbol{0}$
for the time being, even though all results are valid for more general assumptions.
%

The work \cite{Showalter2000} provides the results on existence, uniqueness, and regularity theory for 
\eqref{eq:poroelastic-system}--\eqref{eq:bcs} in the Hilbert space setting, whereas \cite{ShowalterStefanelli2004} 
extends the recent results to a wider class of diffusion problems in the poroelastic media with more general 
material deformation models. Corresponding a priori error estimates can be found in \cite{MuradLoula1992}.
Considered system can be understood as the singular limit of the fully 
dynamic Biot-Allard problem (see the details in \cite{MikelicWheeler2012}), where the acceleration of the solid 
in the mechanics part of \eqref{eq:poroelastic-system} is neglected.

As the Biot model is a coupled system of partial differential equations (PDEs), we have iterative as well as 
monolithic approaches used for solving the problem (see, e.g., \cite{SettariMourits1998}). For the first 
approach, the problem can be reformulated with a contractive operator, which naturally yields 
iteration methods for its solution (see \cite{MikelicWheeler2012}). 
On each step in the time, the flow problem is considered first. It is followed by solving the mechanics using already 
recovered pressure. The procedure is iterated until the desired convergence is reached. Different \authors{alteration} 
of iterative cycles in flow and mechanics, i.e., single- \cite{AlmaniKumarWheeler2017} and multi-rate schemes 
\cite{Kumaretal2015, Almanietal2016}, can be considered. The second approach is fully coupled and considers 
the system with two unknowns simultaneously. 

%
%
The iterative coupling offers several advantages over the monolithic method in the code design. 
In particular, in terms of availability of highly developed discretization methods 
(primal \cite{VermeerVerruijt1981, Reed1984, ZienkiewiczShiomi1984}, 
mixed \cite{MuradLoula1992, MuradLoula1994, MuradThomeeLoula1996}, 
Galerkin least-squares \cite{KorsaweStarke2005}, 
finite volume (FV) \cite{Nordbotten2016},
discontinuous Galerkin (dG) methods \cite{ChenLuoFeng2013}, 
high-order methods \cite{Boffietal2016}
isogeometric analysis \cite{VignolletMayBorst2015}, as well as 
combinations of above-mentioned ones) 
and algebraic solvers (e.g., general Schur complement based preconditioners \cite{ChanPhoonLee2001, Phoonetall2002, WhiteBorja2011, HagaOsnesLangtangen2012, AxelssonBlahetaByczanski2012, Hagaetall2012, CastellettoWhiteTchelepi2015, 
CastellettoWhiteFerronato2016},
and the recently developed robust ones with respect to (w.r.t.) the model parameters 
\cite{Rhebergenetall2014, Rhebergenetall2015, HongKraus2016, LeeMardalWinther2017, BarlandLeeMardalWinther2017}). 
In the fully-coupled methods, constructing  efficient preconditioning techniques for the arising algebraic systems 
\authors{remains a matter of active ongoing scientific research} (see, e.g., 
\cite{Kimetal2009, Whiteetall2016, CastellettoWhiteFerronato2016, GasparRodrigo2017, HongKraus2017}).

The question of a posteriori error control for the poroelastic models has been already addressed using different 
techniques. Application of the residual-based error estimates to the coupled elliptic-parabolic problems can be traced
back to works \cite{Meunier2007, ErnMeunier2009}. Recently, similar error indicators were used in 
\cite{VohralikWheeler2013, Yousef2013, DiPietroVohralikYousef2014, DiPietroFlauraudVohralikYousef2014, 
Ahmedetal2017} 
for immiscible incompressible two- or multi-phase flows in porous media to address the questions of adaptive 
stopping criteria and mesh refinement. 
In \cite{LarssonRunesson2015}, authors suggest an a posteriori error estimator based on the appropriate 
dual problem in space-time for a coupled consolidation problem involving large deformations.
In \cite{Riedlbeck2017, Bertrandetal2017, AhmedRaduNordbotten2019, AhmedNordbottenRadu2020}, 
adaptive space-time algorithms, relying on 
equilibrated fluxes technique, were applied to the Biot's consolidation model (formulated as a 
system with four unknowns). 
\authors{In this work, however, we turn to the functional error estimates (majorants and minorants), which 
are fully computable and provide guaranteed bounds of errors arising in the numerical 
approximations. The derivation of such estimates is based on functional arguments and 
variational formulation of the problem in question. Therefore, the method does not use specific 
properties of approximations (e.g., Galerkin orthogonality) and special properties of the
exact solution (e.g., high regularity). The estimates do not contain mesh dependent constants 
and are valid for any approximation in the natural energy class. Moreover, the majorant also 
yields an efficient error indicator,  that provides mesh adaptation.}


Our main goal is to deduce efficient a posteriori error estimates for the approximation of the system 
\eqref{eq:poroelastic-system} and verify them in application to the numerical problems. 
In \cite{NordbottenRahmanRepinValdman2010}, a posteriori error estimates of the functional type has been 
derived for the stationary Barenblatt--Biot model of porous media. This paper deals with more complicated 
Biot problem presented by an elliptic--parabolic system of partial differential equations. 
Our approach is based on the contraction property of the iterative
method \cite{Banach1922}, which is rather general and not just restricted to fixed-stress 
scheme and functional type estimates of each equation in the Biot system (see, e.g., \cite{RepinDeGruyter2008}). 
To the best knowledge of the authors, it is the first study targeting such a coupling between the elastic 
behavior of the medium and the fluid flow in a context of functional error estimates. 
The main result is presented in Theorem \ref{th:total-estimate}. 
%
Moreover, they not only serve as reliable estimates of the global error measured 
in the energy norm but also as an efficient indicator of the local error distribution over the 
computational domain (confirmed by majorant application to each of the decoupled problems in
\cite{Malietall2014} and references therein). 
The latter property makes functional majorants advantageous in application to 
automated adaptive mesh generation algorithms. 

%

The paper has the following structure: Section \ref{sec:var-formulation-discretization} is dedicated to the 
generalized formulation of the Biot system and its semi-discrete 
counterpart derived after applying the explicit Euler scheme in time. In Sections \ref{sec:fixed-stress-split-iteration-scheme},
we introduce an incremental approach, namely, fixed-stress split scheme, for discretizing considered coupled 
system. In particular, it justifies the optimal choice of parameters in
the iterative scheme and proves that it is a contraction with an explicitly computable convergence rate. 
Sections \ref{sec:error-estimates} and \ref{sec:iterative-estimates} are dedicated to the derivation of auxiliary 
Lemmas used in the proof of Theorems \ref{th:pressure-displacement-estimate} and \ref{th:total-estimate} 
with the general estimates for the approximations generated by the fully decoupled iterative approach. 
Finally, Section \ref{sec:numerical-example} contains a collection of examples, which illustrates the application 
of derived error estimates to the Biot problem. 
%
\section{Variational formulation and discretization}
\label{sec:var-formulation-discretization}

We study approximations of the system \eqref{eq:poroelastic-system}, where 
$\widetilde{\boldsymbol{V}} \equiv H^1(0, T; [\H{1}(\Omega)]^d)$
denotes the space for $u$ (field of displacements) and $\widetilde{W} \equiv H^1(0,T; H^1(\Omega))$ 
is the space for the variable $p$ (pressure).
The generalized setting of \eqref{eq:poroelastic-system} reads: find a pair 
$(\vectoru, p) \in \widetilde{\boldsymbol{V}}_0 \times \widetilde{W}_0$ such that
%
\begin{alignat}{2}
2 \mu \, (\strain(\vectoru), \strain(\vectorv))_Q
+ \lambda \, (\dvrg \vectoru, \dvrg \vectorv)_Q
+ \alpha \, (\nabla p, \vectorv)_Q 
& = (\boldsymbol{f}, \vectorv)_Q, 
\quad \quad \quad \;\;\, 
\forall \vectorv \in \widetilde{\boldsymbol{V}}_{0}, \label{eq:mechanics}\\
%
(\tensorK \nabla p, \nabla w)_Q 
 + ({\partial_t} (\beta \, p + \alpha \, \dvrg \vectoru),  w)_Q 
& = (g, w)_Q, 
\quad \quad \quad  \quad  
 \forall \, w \in \widetilde{W}_0,
\label{eq:flow}
\end{alignat}
where
\begin{alignat*}{2}
\widetilde{\boldsymbol{V}}_{0} &:= 
\big\{ v \in H^1(0, T; [\H{1}(\Omega)]^d) \; \mid \; \vectorv(t)|_{\Sigma^{u}_D} = \boldsymbol{0} \quad \mbox{a.e.} \; t \in (0, T) \big\}, 
\\
\widetilde{W}_{0} & := 
\big\{ w \in H^1(0, T; \H{1}(\Omega)) \quad \mid \; w(t)|_{\Sigma^{p}_D} = 0 \quad \mbox{a.e.} \; t \in (0, T) \big\}.
\end{alignat*}%
The Biot system of type \eqref{eq:mechanics}--\eqref{eq:flow} was analyzed by several authors
to establish existence, uniqueness, and regularity. The first theoretical results on the existence and 
uniqueness of a (weak) solution are presented in \cite{Zenisek1984} for the case of $\beta = 0$. 
Further work in this direction can be found in \cite{Showalter2000, ShowalterSu2001}.
The well-posedness of the quasi-static Biot system is ensured under the above-mentioned assumptions. 
In fact, \cite{MikelicWheeler2013, MikelicWangWheeler2014} established contractive results in suitable 
norms for iterative coupling of \eqref{eq:mechanics}--\eqref{eq:flow}. For an overview of the stability 
of existing iterative algorithms, we refer the reader to \cite{KimTchelepiJuanes2011DrainedSplit,
KimTchelepiJuanes2011FixedStress}.

The system \eqref{eq:mechanics}--\eqref{eq:flow} can be viewed as the two-field formulation of the poroelasticity 
problem. In numerical analysis, there are alternative approaches to treat such a system as 
{\em three- and four-field formulations}. In the case of the three-field model, an additional variable is introduced
to represent the flux in the flow equations, whereas the four-field model considers stress as yet another unknown.
Three-field formulation is rather flexible (since it allows different combinations of discretizations) and is accepted 
by the community as providing more physical approximations of the unknowns than in the two-field case. 
Recently, four-field formulation received increasing attention from the research community, where both equations 
were treated by mixed methods. \authors{The advantages of the latter representation are local conservation of mass and 
momentum balance and more accurate representation of fluxes and stresses.}
The choice of the formulation (from the above-mentioned list) is usually motivated by the considered application 
as well as the restriction on the computational resources.
\authors{For instance, the mixed formulation of \eqref{eq:flow} does not only provide the flux that satisfies the local mass 
conservation property but also generates an effective approximation of this function, which is advantageous to 
functional type error control. It practically minimizes 
the majorant related to the pressure function (see \eqref{eq:majorant-pressure-energy-norm}).
The same method/principle works for the  reconstruction of the stress field. }
The system \eqref{eq:mechanics}--\eqref{eq:flow} is considered in the time-interval $[0, T]$ divided by $N$ 
sub-intervals, such that it forms the corresponding set 
%
$\mathcal{T}_N = \cup_{n = 1}^{N} \overline{{I}^{n}}$, ${I}^{n} = (t^{n-1}, t^{n}).$
%
Let $\vectoru^n(x) \in \boldsymbol{V_0}$ and $p^n(x) \in W_0$, where
\begin{equation}
\boldsymbol{V_0} := \big\{ \vectorv \in \boldsymbol{V} \equiv [\H{1}(\Omega)]^d \, | \, \vectorv \big|_{\Sigma^{\vectoru}_D} = 0 \big\}
\qquad \mbox{and} \qquad
W_0 := \big\{ w \in W \equiv H^1(\Omega) \, | \, w \big|_{\Sigma^{p}_D} = 0 \big\},
\label{eq:spaces}
\end{equation}
respectively, the spatial parts of the solution at $t = t^{n}$. Then, the semi-discrete approximation of
\eqref{eq:mechanics}--\eqref{eq:flow} satisfies the system
\begin{alignat*}{2}
(2 \, \mu \, \strain(\vectoru^n), \strain(\vectorv))_\Omega
+ (\lambda \, \dvrg \vectoru^n, \dvrg \vectorv)_\Omega
+ \alpha \, (\nabla p^n, \vectorv)_\Omega & 
= (\boldsymbol{f}^n, \vectorv)_\Omega, \quad \forall \, \vectorv \in \boldsymbol{V_0},
\\
%
(\tensorK \nabla p^n, \nabla w)_\Omega 
+ \tfrac{1}{\tau^n} \, (\beta (p^{n} - p^{n-1})_\Omega 
+ \alpha \, \dvrg(\vectoru^{n} - \vectoru^{n-1}), w)_\Omega
& = (g^n, w)_\Omega, \quad \forall \, w \in W_0,
\end{alignat*}
where $\tau^n = t^n - t^{n-1}$. This system generates the following problem to be solved on each step 
of the time-incremental method: find the pair $(\vectoru, p)^n \in \boldsymbol{V_0}  \times W_0$
%
\begin{alignat}{2}
(2 \, \mu \, \strain(\vectoru^n), \strain(\vectorv))_\Omega
+ (\lambda \, \dvrg \vectoru^n, \dvrg \vectorv)_\Omega
+ \alpha \, (\nabla p^n, \vectorv)_\Omega & 
= (\boldsymbol{f}^n, \vectorv)_\Omega, \quad \forall \, \vectorv \in \boldsymbol{V_0},
\label{eq:mechanics-semi-discrete-withn}
\\
%
(\tensorK_{\tau^n} \nabla p^n, \nabla w)_\Omega
+ \beta \, (p^{n}, w)_\Omega 
+ \alpha \, (\dvrg \, \vectoru^{n} , w)_\Omega & 
= (\widetilde{g}^n, w)_\Omega, \quad \forall \, w_0 \in W_0,
\label{eq:flow-semi-discrete-withn}
\end{alignat}
where $\tensorK_{\tau^n} := \tau^n \, \tensorK$, the right-hand side of \eqref{eq:flow-semi-discrete-withn} 
is defined as 
\begin{equation}
\widetilde{g}^n 
= \tau^n \, g^n + \beta \, p^{n-1} + \alpha\,  \dvrg \vectoru^{n-1},
\label{eq:g-n}
\end{equation} 
and the pair $(\vectoru, p)^{n-1} \in \boldsymbol{V_0}  \times W_0$
is given by the previous time step. The initial values are chosen as $(\vectoru, p)^{0} = (p^\circ, \vectoru^\circ)$.
Since from now on, we deal only with the semi-discrete counterpart of Biot problem, we omit 
the subscript $\Omega$ in the scalar product. Moreover, we always consider 
\eqref{eq:mechanics-semi-discrete-withn}--\eqref{eq:flow-semi-discrete-withn} on $n$th time step, which 
allows us to omit the superscript $n$ for the rest of the paper and consider the system
 \begin{alignat}{2}
(2 \,\mu \, \strain(\vectoru), \strain(\vectorv))
+ (\lambda \, \dvrg \vectoru, \dvrg \vectorv)
+ \alpha \, (\nabla p, \vectorv) & 
= (\boldsymbol{f}, \vectorv), \quad \forall \, \vectorv \in \boldsymbol{V_0},
\label{eq:mechanics-semi-discrete}
\\
(\tensorK_{\tau} \nabla p, \nabla w)
+ \beta \, (p, w) 
+ \alpha \, (\dvrg \, \vectoru , w) & 
= (\widetilde{g}, w), \quad \forall \, w \in W_0.
\label{eq:flow-semi-discrete}
\end{alignat}
%
This work aims to derive a fully guaranteed a posteriori estimates of the error between the obtained 
approximations $(\tilde{\vectoru}, \tilde{p}) \in \boldsymbol{V_{0h}} \times W_{0h}$, where 
$\boldsymbol{V_{0h}}$ and $W_{0h}$ are discretization spaces of conforming approximations of 
functional spaces $\boldsymbol{V_{0}}$ and $W_{0}$, respectively, and the pair of the exact 
solutions $(\vectoru, p)$ of the Biot system, which is accumulated from the errors on all $N$ time steps, i.e., 
$${e}_{\vectoru} := \vectoru - \tilde{\vectoru} \quad \mbox{and} \quad 
{e}_p := p - \tilde{p}.
$$
%
On each time step, these errors are measured in terms of the combined norm
\begin{alignat}{2}
\! \! \!
\big|\!\big[ ({e}_{\vectoru}, {e}_p) \big]\!\big|
     & \, := |\!|\!|\, {e}_{\vectoru} \,|\!|\!|^2_{\vectoru} + |\!|\!| \, {e}_p \, |\!|\!|^2_{p}.
     \label{eq:error}
\end{alignat}
In turn, each term of the error norm is defined as follows:
%
\begin{alignat}{2}
|\!|\!|\, {e}_{\vectoru} \,|\!|\!|^2_{\vectoru}	& 
:= \| \, \strain({e}_{\vectoru}) \, \|^2_{2\mu}
+ \| \, \dvrg ({e}_{\vectoru})\, \|^2_{\lambda} \quad \mbox{and} \quad 
|\!|\!|\, {e}_p \,|\!|\!|^2_{p} & 
:= \| \, \nabla {e}_p \, \|^2_{\tensorK_{\tau}}
+ \| \, {e}_p \, \|^2_{\beta}, 
\label{eq:error-pressure-displacement}
\end{alignat}
%
%
\sveta{where 
$\| \, w \, \|^2_{\lambda} := \int_\Omega \lambda \, w^2 \, {\rm d}x $, 
$\| \, \strain(\vectorw)  \, \|^2_{2\mu} := \int_\Omega 2\mu \, \strain(\vectorw) : \strain(\vectorw) \,{\rm d}x $, and 
$\|\, \vectorw \, \|^2_{\tensorK_{\tau}} := \int_\Omega \tensorK_\tau \, \vectorw \cdot  \vectorw  \,{\rm d}x $ are 
$L^2$-norms respectively weighted with $2\mu$, $\lambda$, and tensor $\tensorK_\tau$ for any scalar- and 
vector-valued functions} $w$ and $\vectorw$. The global bound of the errors ${e}_{\vectoru}$ and ${e}_p$
contains incremental contributions from each time-interval, i.e., 
\begin{equation}
\sum_{n = 1, \ldots, N} 
\big|\!\big[ ({e}^{n)}_{\vectoru}, {e}^{(n)}_p) \big]\!\big| =: 
\big|\!\big[ ({e}_{\vectoru}, {e}_p) \big]\!\big| 
\leq 
\overline{\rm {M}} \big(\tilde{\vectoru}, \tilde{p}\big)
:= \sum_{n = 1, \ldots, N} \overline{\rm {M}}^{(n)} \big(\tilde{\vectoru}, \tilde{p}\big).
\label{eq:accumulated-error-estimate}
\end{equation}
%

\authors{
For the iterative approach, on each time-step $I^n$, 
the Biot system is decoupled into two sub-problems related to the linear elasticity and a single-phase flow problem.
Then, an iterative procedure is applied to obtain the pair $(\vectoru^{i}, p^{i}) = (\vectoru, p)^{i}$.
Next, each equation is discretized and solved, such that instead of $(\vectoru, p)^{i}$ the pair 
$(\vectoru, p)^{i}_h$, containing the approximation error of the numerical method, is used. 
In Section \ref{sec:error-estimates}, we derive computable a posteriori estimates for this pair of 
approximate solution. 

Functional $\overline{\rm {M}} := \overline{\rm {M}}^{(n)}$ (Theorem \ref{th:total-estimate}) 
is derived by combining the estimates derived for the contractive mapping \cite{Banach1922} 
and the a posteriori error majorants for the elliptic problems (initially introduced \cite{Repin1997, Repin2000}). 
The validity of such estimates is based on the contraction property of specifically constructed linear combination 
of displacement and pressure
$\tfrac{\alpha}{\gamma} \, \dvrg \vectoru^{i} - \tfrac{L}{\gamma} \, p^{i}$, $L, \gamma >0$, 
(the so-called \emph{volumetric mean stress}). The selection of the parameters $L$ and $\gamma$ is justified 
and explained in Section \ref{sec:fixed-stress-split-iteration-scheme}.

\begin{remark}
For the alternative monolithic approach, one solves \eqref{eq:mechanics-semi-discrete}--\eqref{eq:flow-semi-discrete} 
simultaneously for pressure and displacement, reconstructing the pair of approximations 
$(\tilde{\vectoru}, \tilde{p}) = (\vectoru_h, p_h) = (\vectoru, p)_h$. For this case, 
the corresponding computable bound of the error between $(\tilde{\vectoru}, \tilde{p})$ and the exact solution can be derived. 
Such a functional error bound is a combination of a posteriori error estimates for each of the unknowns in 
\eqref{eq:mechanics-semi-discrete} and \eqref{eq:flow-semi-discrete} (see, e.g., 
\cite{NeittaanmakiRepin2004, RepinDeGruyter2008} and references therein). 
For the detailed derivation, we refer the reader to \cite{Kundanetal2018}.
\end{remark}
}
\sveta{
\begin{remark}
We note that 
due to the Korn and Friedrichs' 
inequalities, both $\| \, {e}_{\vectoru} \, \|^2$ and $|\!|\!|\, {e}_{\vectoru} \,|\!|\!|^2_{\vectoru}$ are estimated 
by $\| \,\strain({e}_{\vectoru}) \, \|^2$. Moreover, the physical bound on the Lam\'e parameters is given as 
$d\, \lambda + 2\, \mu > 0$, in the most general case, allowing for the first parameter $\lambda$ to be slightly 
negative for so-called auxetic materials.  In this case, we use the fact that
$$|\!|\!|\, {e}_{\vectoru} \,|\!|\!|^2_{\vectoru}
:= \| \, \strain({e}_{\vectoru}) \, \|^2_{2\mu}
+ \| \, \dvrg ({e}_{\vectoru})\, \|^2_{\lambda} 
\; \mnote{\eqref{eq:inequality-divu-strain}} \leq
(2\, {\mu} + d\, {\lambda} ) \, \| \,\strain({e}_{\vectoru}) \, \|^2$$
holds, and work with the positively-weighted norm $\| \,\strain({e}_{\vectoru}) \, \|^2$. However, as auxetic materials are rare,
the added complexity associated with allowing for such cases has been avoided in this paper.  Consequently, the proofs below are based on the assumption of non-negative Lam\'e parameters. 
\end{remark}
}

%

\section{Fixed-stress splitting scheme}
\label{sec:fixed-stress-split-iteration-scheme}
%

Formal application of the iteration method to \eqref{eq:mechanics-semi-discrete}--\eqref{eq:flow-semi-discrete} 
yields the problem to be solved on the $i$th iteration step:
	\begin{alignat}{2}
	\! \! \! \! \! \! \! \! \! \! \! \! \! \! \! \! \!
	(\tensorK_{\tau} \nabla p^i, \nabla w) 
	+ \beta \, (p^i, w) 
	+ \alpha \, (\dvrg \vectoru^{i-1}, w) 
	& = (\widetilde{g}, w),  \quad \forall \, w \in W_0, 
	\label{eq:flow-iterative} \\
	%
	(2 \, \mu \, \strain (\vectoru^{i}), \strain (\vectorv)) 
	+ (\lambda \, \dvrg \vectoru^{i}, \dvrg \vectorv) 
	+ \alpha \, (\nabla p^{i}, \vectorv) & 
	= (\boldsymbol{f}, \vectorv), \quad \forall \, \vectorv \in \boldsymbol{V_0},
	\label{eq:mechanics-iterative}
	\end{alignat}
%
where the flow equation \eqref{eq:flow-iterative} is solved for $p^{i}$, using $\vectoru^{i-1}$, and 
the elasticity equation \eqref{eq:mechanics-iterative} is used to reconstruct $\vectoru^{i}$ using  $p^{i}$ recovered on 
the previous step.
To obtain the initial data for the iteration procedure, we first set the pressure equal to the 
hydrostatic pressure, i.e., it follows from $\nabla p^{0}  = \rho_f g$ with $\rho_f>0$ denoting the fluid phase density. Whereas 
$\vectoru^{0}$ is reconstructed by \eqref{eq:mechanics-iterative} using $p^{0}$. 
\sveta{The iteration proposed in \eqref{eq:flow-iterative}--\eqref{eq:mechanics-iterative} 
is known as the {\em fixed strain}, and is only conditionally stable.}

To stabilize the iteration scheme \eqref{eq:flow-iterative}--\eqref{eq:mechanics-iterative}, we consider 
the `fixed-stress splitting approach', which properties were initially studied
in \cite{KimTchelepiJuanes2011FixedStress} and 
\cite{MikelicWheeler2013}, respectively. 
This scheme operates with a special quantity: {\em volumetric mean total stress} 
\begin{equation}
\eta^{i} = \tfrac{\alpha}{\gamma} \dvrg \vectoru^{i} - \tfrac{L}{\gamma} \, p^{i} \in W,
\label{eq:relation-sigma-div-u-p}
\end{equation}
where $\gamma$ and $L$ are certain positive tuning parameters. 
These parameters are usually kept constant on each half-time step. 
The optimal choice of $\gamma$ and $L$ allows 
us to prove that this iteration scheme is a contraction in the \sveta{$L^2$-norm 
$\| \delta {\eta}^{i} \|^2$, where $\delta {\eta}^{i} := {\eta}^{i} - {\eta}^{i-1}$}. 
Moreover, it accelerates the iteration procedure by reducing the number of iterations.

\sveta{
By adding $L \, (p^i - p^{i-1})$ into the right-hand side of \eqref{eq:flow-iterative}, 
we rewrite the system \eqref{eq:flow-iterative}--\eqref{eq:mechanics-iterative} using the definition 
\eqref{eq:relation-sigma-div-u-p} as
\begin{alignat}{2}
(\tensorK_{\tau} \nabla {p^{i}}, \nabla w) + 
	(\beta + L) ({p^{i}}, w)
	& = (\widetilde{g} - \gamma \, {\eta}^{i-1} , w) , \qquad \qquad \;\, \forall \, w \in W_0, 
	\label{eq:flow-fully-discrete} \\
\big(2 \, \mu \, \strain ({\vectoru^{i}}), \strain (\vectorv) \big) 
	+ (\lambda \, \dvrg {\vectoru^{i}}, \dvrg \vectorv)
  & = (\boldsymbol{f}^{i} - \alpha  \nabla p^{i}, \vectorv), \qquad \qquad \,\,\forall \, \vectorv \in \boldsymbol{V_0},
\label{eq:mechanics-fully-discrete} 
\end{alignat}
with complimented mixed BCs $p^{i} = 0$ on $\Sigma^{p}_D$ and
$\tensorK_{\tau} \nabla p^{i} \cdot \vectorn = 0$ on $\Sigma^{p}_N$ as well as
$\vectoru^{i} = \boldsymbol{0}$ on $\Sigma^{\vectoru}_D$ and 
$\stress^{i}_{\rm por} \cdot \vectorn = \boldsymbol{0}$ on $\Sigma^{\vectoru}_N$.
%
}

\sveta{
Let  
\begin{equation}
\delta \nabla p^{i} := \nabla p^{i} - \nabla p^{i-1}, \quad  
\strain (\delta \vectoru^{i}) := \strain (\vectoru^{i}) - \strain (\vectoru^{i-1}), \quad
\delta {\eta}^{i} := {\eta}^{i} - {\eta}^{i-1}. 
\label{eq:deltas}
\end{equation}
}
Theorem \ref{th:theorem-contraction} establishes a contraction-type inequality for the norm
$\| \delta {\eta}^{i} \|^2$.
\vskip 15pt
\begin{theorem}[\cite{MikelicWheeler2013,MikelicWangWheeler2014}]
\label{th:theorem-contraction}
\vskip-10pt
If $\gamma = \tfrac{\alpha}{\sqrt{\lambda}}$ and $L \geq \tfrac{\alpha^2}{2 \,\lambda}$, then 
\eqref{eq:flow-fully-discrete}--\eqref{eq:mechanics-fully-discrete} scheme
is a contraction that satisfies the estimate
\begin{equation}
\| \strain (\delta \vectoru^{i}) \|^2_{2\mu}
+ q \, \| \nabla \delta p^{i} \|^2_{\tensorK_{\tau}} 
+ \| \delta {\eta}^{i} \|^2
\leq q^2 \| \delta {\eta}^{i-1} \|^2, \quad q = \tfrac{L}{\beta + L}, 
\label{eq:contration-general-2}
\end{equation}
\sveta{where  $\delta \nabla p^{i}$, $\strain (\delta \vectoru^{i})$, $\delta {\eta}^{i}$ are defined in \eqref{eq:deltas}.}
\end{theorem}

\begin{remark}
\authors{
The estimates in Theorem \ref{th:theorem-contraction}, satisfying contraction estimate \eqref{eq:contration-general-2}, 
also holds for Galerkin approximations $\{\delta \eta_h\}^i \in W_h$, where $W_h$ is a discretization 
space of $W.$  Moreover, in Appendix, we show that similar contraction theorem holds for sequence the
$\{ \delta (\eta - \eta_h)^i \} \in W_h$, where $\eta^i \in W$ is generated by the fixed-stress split iterative scheme defined in 
\eqref{eq:flow-fully-discrete}--\eqref{eq:mechanics-fully-discrete} and $\eta^i_h \in W_h$ is discretization of the latter 
sequence. 
}
\sveta{Generally, it is important to note that all the theorems and lemmas below are formulated for a pair
$(\vectoru, p)^i \in \boldsymbol{V_{0}} \times W_{0}$ that forms a contraction w.r.t. to 
$(\vectoru, p)^{i-1} \in \boldsymbol{V_{0}} \times W_{0}$ and its discrete approximation 
$(\vectoru, p)^i_h \in \boldsymbol{V_{0h}} \times W_{0h}$ that forms a contraction relative to 
$(\vectoru, p)^{i-1}_h \in \boldsymbol{V_{0h}} \times W_{0h}$.}
\end{remark}

\begin{remark}
There exist alternative ways to choose the tuning parameter $L$. In particular, the physically motivated
choice $L_{cl} = \tfrac{\alpha^2}{\lambda + 2 \mu /d}$ is considered in \cite{KimTchelepiJuanes2011FixedStress}.
Whereas, \cite{MikelicWangWheeler2014} suggests $L_{opt} = \tfrac{\alpha^2}{2 \, (\lambda + 2 \mu /d)}$.
The recent study \cite{Bothetall2017} suggests the numerical evidence on the iteration counts w.r.t. the full range 
of the Lam\'{e} parameters for heterogeneous media. Numerical investigation of the optimality of these
parameters and comparison with physically and mathematically motivated values from the literature was done in
\cite{BothKoecher2018}. The authors demonstrated that their optimal value is  
not only dependent on the mechanical material parameters but on the boundary conditions and material 
parameters associated with the fluid flow problem.
\end{remark}

{  
\begin{remark}
The inequality \eqref{eq:contration-general-2} shows that the sequence $\{{\delta \eta}^{i}\}_{i \in \mathds{N}}$ 
is generated by a contractive operator. Therefore, due to the Banach theorem, it tends to a certain fixed point.
Moreover, since all the terms in the left-hand side of \eqref{eq:contration-general-2} are positive, in practice, 
$\{{\delta\eta}^{i}\}_{i \in \mathds{N}}$ might converge with even better contraction rate than $q = \tfrac{L}{\beta + L}$.
\end{remark}

\vskip 10pt
\begin{corollary}
\label{cl:corollary-contraction}
From Theorem \ref{th:theorem-contraction}, it follows that $\nabla \delta p^{i}$ and $\strain (\delta \vectoru^{i})$ 
in \eqref{eq:deltas} are also converging sequences and satisfy 
\begin{equation*}
\| \nabla \delta p^{i} \|^2_{\tensorK_{\tau}} \leq q \| \delta {\eta}^{i-1} \|^2
%
\quad \mbox{and} \quad 
%
\| \strain (\delta \vectoru^{i}) \|^2_{2 \mu} \leq q^2 \| \delta {\eta}^{i-1} \|^2,
\end{equation*}
respectively.
\end{corollary}
%
Corollary \ref{cl:corollary-contraction} is used in the derivation of the error estimate for the term 
$|\!|\!|\, e_p \,|\!|\!|^2_{p}$. In particular, it yields the following result based on the estimates for the Banach 
contraction mappings (see \cite{Banach1922,RepinDeGruyter2008}).

\begin{lemma}[Estimates for contractive mapping]
\label{eq:lemma-1}
Let the assumptions of Theorem \ref{th:theorem-contraction} hold. Then, we have the estimates
\begin{alignat}{2}
\| \nabla (p - p^{i})\|^2_{\tensorK_{\tau}} 
& \leq \tfrac{q}{(1-q)^2} \|  \delta {\eta}^{i-1} \|^2,
\label{eq:estimates-of-nabla-p-via-sigma} \\
%
\| \strain (\vectoru - \vectoru^{i}) \|_{2 \mu}
& \leq \tfrac{q^2}{(1-q)^2} \|  \delta {\eta}^{i-1} \|^2.
\label{eq:estimates-of-strain-via-sigma}
\end{alignat}
\end{lemma}
\ProofBegin
Consider
\begin{alignat*}{2}
\| \nabla (p^{i + m} - p^{i})\|_{\tensorK_{\tau}}
& \leq  \| \nabla (p^{i + m} - p^{i + m - 1})\|_{\tensorK_{\tau}}
  + \ldots + \| \nabla (p^{i + 1} - p^{i})\|_{\tensorK_{\tau}} \\
& \leq q \, (\| {\eta}^{i + m - 1} - {\eta}^{i + m - 2} \| 
  + \ldots + \| {\eta}^{i} - {\eta}^{i-1} \|) \\
& \leq q \, (q^{m} + \ldots + 1) 
  \| {\eta}^{i} - {\eta}^{i-1}\|.
\end{alignat*}
%
By taking a limit $m \rightarrow \infty$ and noting that in this case
$(q^{m} + q^{m-1} +\ldots + 1) \rightarrow \tfrac{1}{1 - q}$,
we arrive at \eqref{eq:estimates-of-nabla-p-via-sigma}. 
The inequality \eqref{eq:estimates-of-strain-via-sigma} is proved by similar arguments.
%
\ProofEnd
\authors{
If in Lemma \ref{eq:lemma-1}, we consider iteration $i$ and $i - m$, $i > m$ as two subsequent 
iterations, a more general version of estimates \eqref{eq:estimates-of-nabla-p-via-sigma} and
\eqref{eq:estimates-of-strain-via-sigma} can be formulated.
\begin{lemma}[General estimates for contractive mapping]
\label{eq:lemma-2}
Let the assumptions of Theorem \ref{th:theorem-contraction} hold. Then, we have the estimates
\begin{alignat}{2}
\| \nabla (p - p^{i})\|^2_{\tensorK_{\tau}} & \leq \min_{1 \leq m \leq i} \Big\{  \tfrac{q^m}{(1-q^m)^2} \| {\eta}^{i} - {\eta}^{i-m}\|^2 \Big\},
\label{eq:estimates-of-nabla-p-via-sigma-2} \\
%
\| \strain (\vectoru - \vectoru^{i}) \|_{2 \mu}& \leq \min_{1 \leq m \leq i} \Big\{ \tfrac{q^{2m}}{(1-q^m)^2} \| {\eta}^{i} - {\eta}^{i-m}\|^2 \Big\}.
\label{eq:estimates-of-strain-via-sigma-2}
\end{alignat}
\end{lemma}
\ProofBegin Proof follows along the lines of the proof of Lemma \ref{eq:lemma-1} with $m = 1, \ldots, i$. \ProofEnd
}
\sveta{
\begin{remark}
\label{eq:remark-lemma-1}
Estimates in Lemma \ref{eq:lemma-2} improves the value of $\tfrac{q^m}{(1-q^m)^2}$ and $\tfrac{q^{2m}}{(1-q^m)^2}$ if 
$q$ is close to 1. It might look contra intuitive, but the choice $m=i$ is not always most optimal.  
Simultaneously with the decrease of quotient with $q$, the term $\| {\eta}^{i} - {\eta}^{0}\|^2$ grows.
Therefore, the choice of $m$ must be made carefully. Computing majorant on each step of our iterative algorithm 
might be computationally expensive, therefore, in demanding applications, the choice of $m$ must be done 
a priori. Alternatively, with several extra iterations made after reaching the desired convergence in $p$ and 
$\vectoru$, both $\tfrac{q\prime}{(1-q\prime)^2}$ with $q\prime = q^m$ and $\| {\eta}^{i} - {\eta}^{i-m}\|^2$ will 
decrease, impacting the total values of the majorant.
\end{remark}
}

\section{Estimates of errors generated by the discretization}
\label{sec:error-estimates}

Before deriving estimates of the approximation errors appearing in the contractive iterative scheme, 
we need to study discretization
errors encompassed in \eqref{eq:flow-fully-discrete}--\eqref{eq:mechanics-fully-discrete} 
for the $i$th iteration. Henceforth, the pair $(\vectoru, p)^{i} = (\vectoru^i, p^i)$ is 
considered as the exact solution of \eqref{eq:flow-fully-discrete}--\eqref{eq:mechanics-fully-discrete}, 
whereas $(\vectoru, p)_h^{i} = (\vectoru_h^i, p_h^i)$ denotes its approximation computed by a certain 
discretization method. We aim to derive computable and reliable estimates of the error measured 
in the terms 
%
$|\!|\!| e_p^{i} |\!|\!|^2_{p}$ 
\mbox{and}
$|\!|\!| e_{\vectoru}^{i} |\!|\!|^2_{\vectoru}.$
\paragraph{Majorant of the error in the pressure term}
\label{ssec:pressure}

For the first equation \eqref{eq:flow-fully-discrete}, Lemma \ref{lem:pressure-full-norm-majorant} 
presents a computable upper bound of the difference $$e^{i}_{p} := p^{i}- p^{i}_h$$ between the 
exact solution $p^{i} \in W_0$ and its approximation $p^i_h \in W_0$, measured in terms of 
the energy norm $|\!|\!|\, e^{i}_{p} \,|\!|\!|^2_{p}$. 
\begin{lemma}
\label{lem:pressure-full-norm-majorant}
For any $p^{i}_h \in W_0$, any auxiliary vector-valued function 
\begin{equation}
\vectorz^i_h \in H_{\Sigma^p_N}(\Omega, \dvrg) 
:= \big\{ \, \vectorz^i_h \in [L(\Omega)]^d \; | \; 
\dvrg \vectorz^i_h \in \L{2}(\Omega), \; 
\vectorz^i_h \cdot \vectorn \in \L{2}(\Sigma^{p}_N)\,\big\},
\label{eq:flux-iter}
\end{equation}
and any parameter $\zeta \geq 0$, we have the following estimate
\begin{equation*}
\| \, \nabla e^{i}_{p} \, \|^2_{\tensorK_{\tau}}
+ \| \, e^{i}_{p} \|^2_\beta
=: |\!|\!|\,  e^{i}_{p} \,|\!|\!|^2_{p}
\leq \overline{\rm M}^h_{p} (p^i_h, \vectorz^i_h; \zeta),
\label{eq:estimate-pressure}
\end{equation*}
where
\begin{equation}
\overline{\rm M}^h_{p} (p^i_h, \vectorz^i_h; \zeta) := 
(1 + \zeta) \, \| \mathbf{r}_{\rm d} (p^i_h, \vectorz^i_h) \|^2_{\tensorK^{-1}_{\tau}} 
+ (1 + \tfrac{1}{\zeta}) \, C^p_\Omega
\Big( 
\|\mathbf{r}_{\rm eq} (p^i_h, \vectorz^i_h) \|^2_{\Omega} 
+ \overline{\rm M}^h_q (({\vectoru}, p)^{1}_h)
+ \| \vectorz^i_h \cdot \vectorn \|^2_{\Sigma^{p}_N} \Big).
\label{eq:majorant-pressure-energy-norm}
\end{equation}
Here,  
\begin{alignat*}{2}
\mathbf{r}_{\rm d} (p^i_h, \vectorz^i_h) & 
:= \vectorz^i_h - \tensorK_{\tau} \nabla p^i_h, 
\quad 
{\mathbf{r}}_{\rm eq} (p^i_h, \vectorz^i_h)
:= \widetilde{g} - \gamma \, {\eta}^{i-1}_h - (\beta + L) \, p^i_h + \dvrg \vectorz^i_h,
\end{alignat*}
where $\widetilde{g}$ is defined in \eqref{eq:g-n}, and 
\sveta{
$$
\overline{\rm M}^h_q := 
\bigg(C_q \, \Big(\frac{\alpha}{\gamma} \, \overline{\rm M}^{h, \; \rfrac{1}{2}}_{p, L^2}
      + \frac{L}{\gamma} \, \overline{\rm M}^{h, \; \rfrac{1}{2}}_{u, \dvrg} \Big)
      + (C_q + 1) \, \| \,\eta^{0} - \eta^{0}_h \, \|\bigg)^2,$$
where $\overline{\rm M}^h_{p, L^2}(p^1_h, \vectorz^1_h)$ and 
$\overline{\rm M}^h_{u, \dvrg} (({\vectoru}, p)^{1}_h, \vectortau^1_h, \vectorz^1_h)$ defined in 
Corollaries \ref{cor:pressure-weak-term-majorant} and \ref{cor:displacement-div-norm} dependent on 
explicitly given ${\eta}^{0}$. }Constant 
\begin{equation}
(C^p_\Omega)^2
:= \tfrac{1}{(\beta + L)}\Big(1 + \big(C^{\rm tr}_{\Sigma^p_N}\big)^2\Big)
\label{eq:cp-omega-iter}
\end{equation}
is defined via the constant in the trace-type inequality
\begin{equation}
\|\, w \,\|_{\Sigma^p_N} \leq C^{\rm tr}_{\Sigma^p_N} \,\| \, w \, \|_{\Omega}, 
\quad \forall w \in W_0,
\label{eq:trace-iter}
\end{equation}
and positive parameters of the Biot model $\beta$ and $L$. 
\end{lemma}

\ProofBegin
The majorant $\overline{\rm M}_p (p^i_h, \vectorz^i_h; \zeta)$ follows from 
\cite[Section 2]{RepinSauter2006} and \cite[Section 4.2--4.3]{RepinDeGruyter2008}, i.e., 
we consider \eqref{eq:flow-fully-discrete} with subtracted from its left- and right-hand side bilinear form 
$(\tensorK_{\tau} \nabla {p^{i}_h}, \nabla w) + (\beta + L) ({p^{i}}_h, w)$ 
\begin{alignat*}{2}
(\tensorK_{\tau} \nabla e^{i}_{p}, \nabla w) + (\beta + L) \, (e^{i}_{p}, w)
& = (\widetilde{g} - \gamma \, {\eta}^{i-1} - (\beta + L) \,{p^{i}}_h, w)
- (\tensorK_{\tau} \nabla {p^{i}_h}, \nabla w).
\end{alignat*}
Next, we set $w = e^{i}_{p}$ and introduce an auxiliary function 
$\vectorz^i_h \in H_{\Sigma^p_N}(\Omega, \dvrg)$ (cf. \eqref{eq:flux-iter}) satisfying the identity 
\linebreak
$(\dvrg \vectorz^i_h, \vectorw)_\Omega + (\vectorz^i_h, \nabla \vectorw)_\Omega 
= (\vectorz^i_h \cdot \vectorn, \vectorw)_{\Sigma^{p}_N}$, such that
\begin{alignat}{2}
\| \nabla e^{i}_{p} \|^2_{\tensorK_\tau}+ \| e^{i}_{p} \|^2_{\beta + L} 
& = 
(\vectorz^i_h - \tensorK_{\tau} \nabla {p^{i}_h}, \nabla e^{i}_{p}) 
+ (\widetilde{g} - \gamma \, {\eta}^{i-1} - (\beta + L) \, {p^{i}}_h + \dvrg \vectorz^i_h , e^{i}_{p})
-  (\vectorz^i_h \cdot \vectorn, e^{i}_{p})_{\Sigma^{p}_N} \nonumber\\
& = (\mathbf{r}_{\rm d} (p^i_h, \vectorz^i_h), \nabla e^{i}_{p}), \nabla e^{i}_{p}) 
+ (\widetilde{\mathbf{r}}_{\rm eq} (p^i_h, \vectorz^i_h), e^{i}_{p}) , e^{i}_{p})
-  (\vectorz^i_h \cdot \vectorn, e^{i}_{p})_{\Sigma^{p}_N},
\label{eq:fluid-h-minus-vh}
\end{alignat}
where $\widetilde{\mathbf{r}}_{\rm eq} (p^i_h, \vectorz^i_h) := \widetilde{g} - \gamma \, {\eta}^{i-1} - (\beta + L) \, p^i_h + \dvrg \vectorz^i_h.$
Using the H\"older and Young inequalities, the first term on the right-hand side of \eqref{eq:fluid-h-minus-vh} can be estimated as  
\begin{equation}
(\mathbf{r}_{\rm d} (p^i_h, \vectorz^i_h), \nabla e^{i}_{p}) 
\leq \tfrac{1}{2}\, (1 + \zeta)\, \| \mathbf{r}_{\rm d} (p^i_h, \vectorz^i_h)\|^2_{1/\tensorK_{\tau}}
+ \tfrac{1}{2 (1 + \zeta)}\, \| \nabla e^{i}_{p} \|^2_{\tensorK_{\tau}}.
\label{eq:estimate-1}
\end{equation}
The second term on the right-hand side of \eqref{eq:fluid-h-minus-vh} is bounded analogously, i.e., 
\begin{alignat}{2}
(\widetilde{\mathbf{r}}_{\rm eq} (p^i_h, \vectorz^i_h), e^{i}_{p})
-  (\vectorz^i_h \cdot \vectorn, e^{i}_{p})_{\Sigma^{p}_N} & 
\leq 
\tfrac{1}{2} \, (1 + \tfrac{1}{\zeta}) \,(C^p_{\Omega})^2 
(\| \widetilde{\mathbf{r}}_{\rm eq} (p^i_h, \vectorz^i_h) \|^2 + \| \vectorz^i_h \cdot \vectorn \|^2_{\Sigma^{p}_N}) \, 
+ \tfrac{1}{2} \, \tfrac{\zeta}{1 + \zeta} \, \| e^{i}_{p} \|^2_{\beta + L},
\label{eq:estimate-2}
\end{alignat}
where $C^p_{\Omega}$ (cf. \eqref{eq:cp-omega-iter}) is a constant in the inequality 
$$\| w \|^2 + \| w \|^2_{\Sigma^{p}_N} 
\leq (C^p_{\Omega})^2 \, \| w \|^2_{\beta + L}, 
\qquad 
\forall w \in W_0,$$
defined in \eqref{eq:trace-iter}.
By summing up the results of \eqref{eq:estimate-1} and \eqref{eq:estimate-2}, we obtain 
\begin{alignat}{2}
\| \nabla e^{i}_{p} \|^2_{\tensorK_\tau}  + \| e^{i}_{p} \|^2_{\beta} 
& \leq
\| \nabla e^{i}_{p} \|^2_{\tensorK_\tau}  + \| e^{i}_{p} \|^2_{\beta + L} \nonumber\\
& \leq
(1 + \zeta) \, \| {\mathbf{r}}_{\rm d} (p^i_h, \vectorz^i_h)\|^2_{1/\tensorK_{\tau}} 
+ (1 + \tfrac{1}{\zeta}) \, (C^p_{\Omega})^2 (\| \widetilde{\mathbf{r}}_{\rm eq} (p^i_h, \vectorz^i_h) \|^2 + \| \vectorz^i_h \cdot \vectorn \|^2_{\Sigma^{p}_N}).
\end{alignat}
\sveta{
At this point, the term $\| \widetilde{\mathbf{r}}_{\rm eq} (p^i_h, \vectorz^i_h) \|^2$ 
is not fully computable in the usual sense 
of functional majorants since it is defined using $\eta_{i-1}$. However, it can be estimated by 
\begin{alignat}{2}
\| \widetilde{\mathbf{r}}_{\rm eq} (p^i_h, \vectorz^i_h) \|^2
& \leq \| \widetilde{g} - \gamma \, {\eta}^{i-1} - (\beta + L) \, p^i_h + \dvrg \vectorz^i_h \|^2 \nonumber\\
& \leq 2 \| \widetilde{g} - \gamma \, {\eta}^{i-1}_h - (\beta + L) \, p^i_h + \dvrg \vectorz^i_h \|^2 
   + 2 \gamma^2 \|   {\eta}^{i-1} - {\eta}^{i-1}_h  \|^2 \nonumber\\
& \leq  2 \| {\mathbf{r}}_{\rm eq} (p^i_h, \vectorz^i_h) \|^2 
   + 2 \gamma^2 \|   {\eta}^{i-1} - {\eta}^{i-1}_h  \|^2.
\label{eq:estimation-req-diff-etai}
\end{alignat}
Consider the norm $\| {\eta}^{i-1} - {\eta}^{i-1}_h\|$ (without the squares) and apply an approach similar to that was used in the proof of 
Lemma \ref{eq:lemma-1}
\begin{alignat*}{2}
\| {\eta}^{i} - {\eta}^{i}_h\| \leq \sum_{k = 1}^{i-1} \, \| \,\delta \eta^{k} - \delta \eta^{k}_h \, \| +  \| \,\eta^{0} - \eta^{0}_h \, \| 
\; \mnote{Theorem \ref{th:theorem-contraction-appendix}} 
& \leq (q^{i-1} + \ldots + 1) \, \| \,\delta \eta^{1} - \delta \eta^{1}_h \, \| +  \| \,\eta^{0} - \eta^{0}_h \, \| \\
& \leq \Big(\sum_{k=1}^{i-1} q^k + 1\Big) \, \| \,\delta \eta^{1} - \delta \eta^{1}_h \, \| +  \| \,\eta^{0} - \eta^{0}_h \, \| \\
& \leq \Big(\sum_{k=1}^{i-1} q^k + 1\Big) \, \big(\| \,\eta^{1} - \eta^{1}_h \, \| + \| \,\eta^{0} - \eta^{0}_h \, \| \big) +  \| \,\eta^{0} - \eta^{0}_h \, \| \\
& \leq C_q \, \| \,\eta^{1} - \eta^{1}_h \, \| + (C_q + 1) \, \| \,\eta^{0} - \eta^{0}_h \, \|.
\end{alignat*}
Here, $\| \,\eta^{0} - \eta^{0}_h \, \|$ is a computable term, whereas $\| {\eta}^{1} - {\eta}^{1}_h\|^2$ is controlled 
by combination of majornats $\overline{\rm M}^h_{p, L^2}(p^1_h, \vectorz^1_h)$ and 
$\overline{\rm M}^h_{u, \dvrg} (({\vectoru}, p)^{1}_h, \vectortau^1_h, \vectorz^1_h)$ defined in 
Corollaries \ref{cor:pressure-weak-term-majorant} and \ref{cor:displacement-div-norm} containing computable 
$\eta^0 := \tfrac{\alpha}{\gamma} \dvrg \vectoru^{0} - \tfrac{L}{\gamma} \, p^{0}$., i.e., 
\begin{alignat}{2}
\| {\eta}^{1} - {\eta}^{1}_h\|
& \leq \frac{\alpha}{\gamma} \, \| \, \dvrg (\vectoru^{1} - \vectoru^{1}_h) \, \|
   + \frac{L}{\gamma} \, \| \,  p^{1} - p^{1}_h \, \|
   \leq \frac{\alpha}{\gamma} \, \Big(\overline{\rm M}^h_{p, L^2}(p^1_h, \vectorz^1_h; \zeta) \Big)^{\rfrac{1}{2}}
      + \frac{L}{\gamma} \, \Big(\overline{\rm M}^h_{u, \dvrg} (({\vectoru}, p)^{1}_h, \vectortau^1_h, \vectorz^1_h)\Big)^{\rfrac{1}{2}}.
\label{eq:etimate-eta1-etah-1}
\end{alignat}
This yields that
\begin{alignat}{2}
\| {\eta}^{i} - {\eta}^{i}_h\|^2 & \leq \overline{\rm M}^h_q := 
\bigg(C_q \, \Big(\frac{\alpha}{\gamma} \, \overline{\rm M}^{h, \; \rfrac{1}{2}}_{p, L^2}
      + \frac{L}{\gamma} \, \overline{\rm M}^{h, \; \rfrac{1}{2}}_{u, \dvrg} \Big)
      + (C_q + 1) \, \| \,\eta^{0} - \eta^{0}_h \, \|\bigg)^2
\label{eq:majorant-q}
\end{alignat}
The combination of \eqref{eq:estimation-req-diff-etai} and \eqref{eq:majorant-q} yields
\begin{alignat*}{2}
\| \nabla e^{i}_{p} \|^2_{\tensorK_\tau}  + \| e^{i}_{p} \|^2_{\beta} 
\leq
(1 + \zeta) \, \| {\mathbf{r}}_{\rm d} (p^i_h, \vectorz^i_h)\|^2_{1/\tensorK_{\tau}} 
+ (1 + \tfrac{1}{\zeta}) \, (C^p_{\Omega})^2 (  2 \| {\mathbf{r}}_{\rm eq} (p^i_h, \vectorz^i_h) \|^2  
+ \overline{\rm M}^h_q(({\vectoru}, p)^{1}_h)
+ \| \vectorz^i_h \cdot \vectorn \|^2_{\Sigma^{p}_N}).
\end{alignat*}
}
\ProofEnd

\begin{remark}
The numerical reconstruction of the majorant involves several steps. They are motivated by the accuracy 
requirements imposed on the upper bound of the error. To generate guaranteed bounds with the 
realistic efficiency index $\Ieff (\overline{\rm M}_p) := 
\tfrac{\overline{\rm M}_p}{|\!|\!|\, p^{i} - p^i_h \,|\!|\!|^2_{p}}$, 
we can reconstruct $\vectorz^i_h$ from $\nabla p^i_h$ (where $p^i_h$ is approximated by the chosen 
discretization method recovering the exact solution of \eqref{eq:flow-fully-discrete}). 
However, to obtain the sharpest estimate, functional $\overline{\rm M}_p$ must be optimized w.r.t. 
$\vectorz^i_h$ and $\zeta$ iteratively. This generates an auxiliary variational problem w.r.t. vector-valued 
function $\vectorz^i_h$. 

Alternatively, one can consider the mixed formulation of \eqref{eq:flow-fully-discrete} and reconstruct the pair 
$(p^i_h, \vectorz^i_h)$ simultaneously using one of the well-developed mixed methods 
\cite{ArnoldBrezzi1985, Arnold1990}. Then, both variables required for the reconstruction of 
$\overline{\rm M}_{p}$ are directly computable, and no additional post-processing 
(computational overhead) is required. 
\end{remark}

Majorant in Lemma \ref{lem:pressure-full-norm-majorant} yields an estimate of the $e^{i}_p$ 
measured in terms of $\L{2}$-norm. 
%
\begin{corollary}
\label{cor:pressure-weak-term-majorant}
For any $p^{i}_h \in W_0$, any auxiliary functions and parameters defined in Lemma 
\ref{lem:pressure-full-norm-majorant}, the estimate
%
\begin{equation}
\| e^{i}_{p} \|^2
\leq \overline{\rm M}^h_{p, L^2}(p^i_h, \vectorz^i_h; \zeta) := 
\Big(\tau \,\lambda_{\tensorK}\big(C^{\rm F}_{\Sigma^p_D}\big)^{-2} + \beta \Big)^{-1} \,
\overline{\rm M}^h_p (p^i_h, \vectorz^i_h; \zeta)
\label{eq:majorant-pressure-in-l2-norm}
\end{equation}
holds, where $\overline{\rm M}^h_{p}(p^i_h, \vectorz^i_h; \zeta)$ is defined in \eqref{eq:majorant-pressure-energy-norm},
$C^{\rm F}_{\Sigma^p_D}$ is a constant in Friedrichs' inequality (cf. \eqref{eq:Poincare}), and 
$\lambda_\tensorK$ is the minimum eigenvalue of the permeability tensor (cf. \eqref{eq:eigenvalue-K}).
\end{corollary}
%
\ProofBegin
By means of the Friedrichs' inequality and \eqref{eq:eigenvalue-K}, we obtain 
\begin{equation}
|\!|\!| e^{i}_{p}  |\!|\!|^2_{p} 
\geq \tau \, \lambda_{\tensorK} \, \big(C^{\rm F}_{\Sigma^p_D}\big)^{-2} \| e^{i}_p  \|^2 
     + \| e^{i}_{p}  \|^2_{\beta} \nonumber \\
\geq \Big(\tau \,\lambda_{\tensorK}\big(C^{\rm F}_{\Sigma^p_D}\big)^{-2} + \beta \Big) 
  \| \, e^{i}_{p} \, \|^2.
\label{eq:h1-l2-connection}
\end{equation}
By combining \eqref{eq:h1-l2-connection} and \eqref{eq:majorant-pressure-energy-norm}, we arrive at 
\eqref{eq:majorant-pressure-in-l2-norm}. 
\ProofEnd


%
\paragraph{Majorant of the error in the displacement term}
\label{ssec:displacement}

Current section considers estimates for the error 
\begin{equation}
e^{i}_{\vectoru} := {\vectoru}^{i} - {\vectoru}^{i}_h
\label{eq:eiu}
\end{equation}
between the exact solution ${\vectoru}^{i} \in \boldsymbol{V_0}$ and its respective approximation 
${\vectoru}^{i}_h \in \boldsymbol{V_0}$ measured in terms of the energy norm 
$|\!|\!|\, \cdot \,|\!|\!|^2_{\vectoru}$ (cf. \eqref{eq:error-pressure-displacement}). 
Since $p^i_h$ is, in fact, used instead of $p^i$, the original problem \eqref{eq:mechanics-fully-discrete} is 
replaced by 
\begin{equation}
2\mu \, \big(\strain ({\tilde{\vectoru}^{i}}), \strain (\vectorv) \big) 
+ \lambda (\dvrg {\tilde{\vectoru}^{i}}, \dvrg \vectorv)
 = (\boldsymbol{f}^{i} - \alpha  \nabla p^{i}_h, \vectorv), 
 \quad \forall \, \vectorv \in \boldsymbol{V_0},
\label{eq:mechanics-h-perturbed}
\end{equation}
with a perturbed right-hand side. 
Therefore, ${\vectoru}^{i}_h$ is an approximation of $\tilde{\vectoru}^{i}$ instead of ${\vectoru}^{i}$. 
In other words, $e^{i}_{\vectoru}$ is composed of the error arising due to the original problem is 
replaced by \eqref{eq:mechanics-h-perturbed}, i.e., ${\vectoru}^{i} - \tilde{\vectoru}^{i}$, and the error 
$\tilde{\vectoru}^{i} - {\vectoru}^{i}_h$ arising because \eqref{eq:mechanics-h-perturbed} is solved
approximately. By means of the triangle inequality, $e^{i}_{\vectoru}$ can be estimated by above-described 
errors as follows:
\begin{alignat}{2}
& \| \, \strain (e^{i}_{\vectoru}) \, \|^2_{2\mu}
+ \| \, \dvrg (e^{i}_{\vectoru}) \, \|^2_{\lambda} 
=: |\!|\!| e^{i}_{\vectoru}|\!|\!|^2_{\vectoru}
\leq 2 \, |\!|\!| \, {\vectoru}^{i} - \tilde{\vectoru}^{i} |\!|\!|^2_{\vectoru}
+ 2 \, |\!|\!| \, \tilde{\vectoru}^{i} - {\vectoru}^{i}_h \, |\!|\!|^2_{\vectoru}.
\label{eq:comb-error}
\end{alignat}
Here, $|\!|\!| \, \tilde{\vectoru}^{i} - {\vectoru}^{i}_h \, |\!|\!|^2_{\vectoru}$ can be estimated by the functional
majorant for a class of the elasticity problems (see Lemma \ref{lem:displacement-func-estimate}), whereas 
$|\!|\!| \, {\vectoru}^{i} - \tilde{\vectoru}^{i} |\!|\!|^2_{\vectoru}$ is controlled by the bound following from 
the difference of model problems \eqref{eq:mechanics-fully-discrete} and \eqref{eq:mechanics-h-perturbed} 
(see Lemma \ref{lem:displacement-estimate}).

\begin{lemma}
\label{lem:displacement-func-estimate}
For any $\vectoru^i_h \in \boldsymbol{V_0}$ approximating $\tilde{\vectoru}^{i}$ in 
\eqref{eq:mechanics-h-perturbed},
any auxiliary tensor-valued function
\begin{alignat*}{2}
\vectortau^i_h \in [\mathcal{T}_{\rm Div}(\Omega)]^{d \times d} := 
\Big \{\, \vectortau^i_h \in [\L{2}(\Omega)]^{d \times d} \;\big|\; 
	& {\rm Div} {\vectortau^i_h} \in [\L{2}(\Omega)]^d, 
	\; \vectortau^i_h \cdot \vectorn \in \L{2}(\Sigma^{\vectoru}_N)
	\, \Big \},
\end{alignat*}
and any parameter $\xi \geq 0$, we have the estimate 
%
\begin{equation}
\begin{alignedat}{2}
\|\strain(\tilde{\vectoru}^{i} - {\vectoru}^{i}_h) \|^2_{2\mu}
& + \|\dvrg (\tilde{\vectoru}^{i} - {\vectoru}^{i}_h) \|^2_{\lambda} 
=: |\!|\!| \tilde{\vectoru}^{i} - {\vectoru}^{i}_h |\!|\!|^2_{\vectoru}
\leq \overline{\rm M}_{\tilde{u}}(({\vectoru}, {p})^i_h, \vectortau^i_h) \\
 & 
 := (1 + \xi)\, \int_\Omega \mathbf{r}_{\rm d}(p^{i}_h, \vectortau^i_h) \dx 
 + (1 + \tfrac{1}{\xi}) \, C^{\vectoru}_{\Omega}
 \Big( \| \mathbf{r}_{\rm eq}(p^{i}_h, \vectortau^i_h) \|^2_{\Omega} 
 + \| \vectortau^i_h \cdot \vectorn\|^2_{\Sigma^{\vectoru}_N} \Big), 
\end{alignedat}
\label{eq:majorant-displpacement-energy-norm}
\end{equation}
where 
\begin{equation}
\begin{alignedat}{2}
\mathbf{r}_{\rm eq}(p^{i}_h, \vectortau^i_h) 
& := \boldsymbol{f}^{i} - \alpha  \nabla p^{i}_h + {\rm Div} \vectortau^i_h, \\
\mathbf{r}_{{\rm} d, \mu, \lambda}({\vectoru}^{i}_h, \vectortau^i_h) 
 & := 2 \,\mu \,|\strain({\vectoru}^{i}_h)|^2
      + \lambda \,|\dvrg {\vectoru}^{i}_h|^2
      + \tfrac{1}{2\mu} \, (|\vectortau^i_h|^2 - \tfrac{\lambda}{3\lambda+2\mu} |\dvrg \vectortau^i_h|^2)
      - 2 \, \strain({\vectoru}^{i}_h) : \tau^i,
\end{alignedat}
\label{eq:residuals-p-iterative}
\end{equation}
$\beta$, $\alpha$, $\mu$, $\lambda$ are characteristics of the Biot model, and 
\begin{equation}
C^{\vectoru}_{\Omega} :=  C^{\rm K}  \, \big(1 + C^{\rm tr}_{\Sigma^{\vectoru}_N}\big)
\label{eq:cu-omega}
\end{equation}
is defined through the constants $C^{\rm tr}_{\Sigma^{\vectoru}_N}$ and $C^{\rm K}$ in 
trace-type and the Korn first inequalities
%
\begin{equation}
\|\, \vectorw \,\|_{\Sigma^{\vectoru}_N} \leq 
C^{\rm tr}_{\Sigma^{\vectoru}_N} \,\| \vectorw \|_{[H^1(\Omega)]^d} \quad \mbox{and} \quad
\| \vectorw \|_{[H^1(\Omega)]^d} 
\leq C^{\rm K}
\|\strain(\vectorw)\|_{[L^2(\Omega)]^{{d} \times {d}}}
\quad \forall \vectorw \in \boldsymbol{V_0},
\label{eq:trace-inequality}
\end{equation}
%
%
respectively.

\end{lemma}
%
\ProofBegin
For the simplicity of exposition, let us assume the following representation of the elasticity tensor 
\begin{equation}
\tensorL \,\strain (\vectoru) := 2\, \mu \, \strain(\vectoru) + \lambda \, \dvrg (\vectoru).
\label{eq:l-tensor-via-strain-divu}
\end{equation}
Then, the derivation of an a posteriori error estimate for the problem 
\begin{equation}
(\tensorL \,\strain ({\tilde{\vectoru}^{i}}), \strain (\vectorv) \big) 
= (\boldsymbol{f}^{i} - \alpha  \nabla p^{i}_h, \vectorv), 
\quad \forall \, \vectorv \in \boldsymbol{V_0},
\label{eq:mechanics-h-perturbed-strain-tensor}
\end{equation}
follows the lines presented in \cite[Section 5.2]{RepinDeGruyter2008}. 
In particular, considering an approximation $\vectoru^i_h \in \boldsymbol{V_0}$, 
we subtract bilinear form $(\tensorL \strain ({\vectoru}^{i}_h), \strain (\vectorv) \big)$ from the left- and right-hand side 
of \eqref{eq:mechanics-h-perturbed-strain-tensor} and set 
$\vectorv = \tilde{\vectoru}^{i} - {\vectoru}^{i}_h$ to obtain 
\begin{equation}
(\tensorL \,\strain ({\tilde{\vectoru}^{i}} - {\vectoru}^{i}_h), \strain ({\tilde{\vectoru}^{i}} - {\vectoru}^{i}_h) \big) 
 = (\boldsymbol{f}^{i} - \alpha  \nabla p^{i}_h, \vectorv)
 - \big(\tensorL \, \strain ({{\vectoru}}), \strain ({\tilde{\vectoru}^{i}} - {\vectoru}^{i}_h) \big).
\label{eq:mechanics-h-minus-vh}
\end{equation}
Next, we set $\vectorv = {\tilde{\vectoru}^{i}} - {\vectoru}^{i}_h$ and add the divergence of 
the tensor-valued function $\vectortau^i_h \in [\mathcal{T}_{\rm Div}(\Omega)]^{d \times d}$, i.e., 
%
\begin{equation}
({\rm Div} \vectortau^i_h, \vectorv) + (\vectortau^i_h, \strain(\vectorv)) = 
(\vectortau^i_h \cdot \vectorn, \vectorv)_{\Sigma^{\vectoru}_N}, 
\quad \forall \vectorv \in \boldsymbol{V_0}, 
\label{eq:tau-identity}
\end{equation}
into the left- and right-hand side of \eqref{eq:mechanics-h-minus-vh}, which results in the identity
\begin{equation}
(\tensorL \,\strain ({\tilde{\vectoru}^{i}} - {\vectoru}^{i}_h), 
\strain ({\tilde{\vectoru}^{i}} - {\vectoru}^{i}_h) \big) 
 = \big(\mathbf{r}_{{\rm d}, \tensorL}({\vectoru}^{i}_h, \vectortau^i_h) , \strain ({\tilde{\vectoru}^{i}} - {\vectoru}^{i}_h) \big)
 + (\mathbf{r}_{\rm eq}(p^{i}_h, \vectortau^i_h), {\tilde{\vectoru}^{i}} - {\vectoru}^{i}_h)
 - (\vectortau^i_h \cdot \vectorn, \tilde{\vectoru}^{i} - {\vectoru}^{i}_h)_{\Sigma^{\vectoru}_N},
\label{eq:mechanics-h-minus-vh-2}
\end{equation}
where 
\begin{equation*}
\mathbf{r}_{\rm d, \tensorL}({\vectoru}^{i}_h, \vectortau^i_h) 
 := \vectortau^i_h - \tensorL \, \strain (\vectoru^{i}_h)
\end{equation*}
and $\mathbf{r}_{\rm eq}(p^{i}_h, \vectortau^i_h)$ is defined \eqref{eq:residuals-p-iterative}. 
By means of the H\"older and Young inequalities, the first term 
on the right-hand side of \eqref{eq:mechanics-h-minus-vh} can be estimated as  
\begin{equation*}
\big(\mathbf{r}_{\rm d, \tensorL}({\vectoru}^{i}_h, \vectortau^i_h), 
           \strain ({\tilde{\vectoru}^{i}} - {\vectoru}^{i}_h) \big)
\leq 
\big \| \mathbf{r}_{\rm d, \tensorL}({\vectoru}^{i}_h, \vectortau^i_h)  \|_{\tensorL^{-1}}
      \| \strain ({\tilde{\vectoru}^{i}} - {\vectoru}^{i}_h) \|_{\tensorL}
\leq \tfrac{\alpha_1}{2} \, \| \mathbf{r}_{\rm d}({\vectoru}^{i}_h, \vectortau^i_h) \|^2_{\tensorL^{-1}}
+ \tfrac{1}{2 \, \alpha_1} \| \strain (\tilde{\vectoru}^{i} - {\vectoru}^{i}_h) \|^2_{\tensorL}
%
\end{equation*}
The second and the third terms are combined and estimated as follows
\begin{equation*}
(\mathbf{r}_{\rm eq}(p^{i}_h, \vectortau^i_h), {\tilde{\vectoru}^{i}} - {\vectoru}^{i}_h)
 - (\vectortau^i_h \cdot \vectorn, \tilde{\vectoru}^{i} - {\vectoru}^{i}_h)_{\Sigma^{\vectoru}_N}
 \leq 
\tfrac{\alpha_2}{2} \, (C^{\vectoru}_{\Omega})^2 \, 
(\| \mathbf{r}_{\rm eq}({p}^{i}_h, \vectortau^i_h) \|^2 
+ \| \vectortau^i_h \cdot \vectorn \|^2_{{\Sigma^{\vectoru}_N}} )
+ \tfrac{1}{2 \, \alpha_2} \| \strain (\tilde{\vectoru}^{i} - {\vectoru}^{i}_h) \|^2_{\tensorL} 
\end{equation*}
where $C^{\vectoru}_{\Omega}$ (cf. \eqref{eq:cu-omega}) is a constant in 
$$\| {\tilde{\vectoru}^{i}} - {\vectoru}^{i}_h \|^2 
+ \| \tilde{\vectoru}^{i} - {\vectoru}^{i}_h\|^2_{\Sigma^{\vectoru}_N} 
\leq (C^{\vectoru}_{\Omega})^2 \| \strain (\tilde{\vectoru}^{i} - {\vectoru}^{i}_h) \|^2_{\tensorL}$$ 
defined through constants
$C^{\rm K}$ and $C^{\rm tr}_{\Sigma^{\vectoru}_N}$  in the Korn and trace inequalities
defined in \eqref{eq:korn} and \eqref{eq:trace}, respectively.
By choosing parameters $\alpha_1 = (\xi +1)$, $\alpha_2 = (1 + \tfrac{1}{\xi})$, 
where $\xi >0$, we arrive at
\begin{equation}
\| \strain ({\tilde{\vectoru}^{i}} - {\vectoru}^{i}_h) \|_{\tensorL}
 \leq 
 (1 + \xi)\, \| \mathbf{r}_{\rm d, \tensorL}({\vectoru}^{i}_h, \vectortau^i_h) \|^2_{\tensorL^{-1}}
 + (1 + \tfrac{1}{\xi}) \, (C^{\vectoru}_{\Omega})^2
 \big(\| \mathbf{r}_{\rm eq}(p^{i}_h, \vectortau^i_h)\|^2 
 + \| \vectortau^i_h \cdot \vectorn\|^2_{\Sigma^{\vectoru}_N}\big). 
\label{eq:majorant-L-tensor}
\end{equation}

Consider now \eqref{eq:l-tensor-via-strain-divu} and the representation of tensor 
$\tensorL^{-1} \vectortau^i_h$ through Lame parameters, i.e., 
$$\tensorL^{-1} \vectortau^i_h := 
\tfrac{1}{2\mu} \big(\vectortau^i_h - \tfrac{\lambda}{3\lambda+2\mu} \dvrg \vectortau^i_h \tensorI\big).$$
Then, the first term on the right-hand side of \eqref{eq:majorant-L-tensor} can be rewritten as
\begin{alignat}{2}
\tensorL^{-1} \mathbf{r}_{\rm d, \tensorL}({\vectoru}^{i}_h, \vectortau^i_h) &  :
\mathbf{r}_{\rm d, \tensorL}({\vectoru}^{i}_h, \vectortau^i_h)
= \big(\tensorL^{-1} \vectortau^i_h - \strain ({{\vectoru}^{i}_h})\big) :
   (\vectortau^i_h - \tensorL \strain ({{\vectoru}^{i}_h})\big) \nonumber\\
& = 2 \,\mu \,|\strain({\vectoru}^{i}_h)|^2
      + \lambda \,|\dvrg {\vectoru}^{i}_h|^2
      + \tfrac{1}{2\mu} \, (|\vectortau^i_h|^2 - \tfrac{\lambda}{3\lambda+2\mu} |\dvrg \vectortau^i_h|^2)
      - 2 \, \strain({\vectoru}^{i}_h) : \tau^i 
   =: {\rm r}_{\rm d, \mu, \lambda}.
\label{eq:residuals-p-iterative-lame}     
\end{alignat}
Taking the latter into account, we arrive at the alternative estimate
\begin{alignat}{2}
\|\strain({\tilde{\vectoru}^{i}} - {\vectoru}^{i}_h) \|^2_{2\mu}
& + \|\dvrg ({\tilde{\vectoru}^{i}} - {\vectoru}^{i}_h) \|^2_{\lambda} \nonumber\\ 
& \leq  (1 + \xi)\, \int_\Omega {\rm r}_{\rm d, \mu, \lambda}({\vectoru}^{i}_h, \vectortau^i_h) \dx 
 + (1 + \tfrac{1}{\xi}) \, (C^{u}_{\Omega})^2
 \big(\| \mathbf{r}_{\rm eq}(p^{i}_h, \vectortau^i_h)\|^2 
 + \| \vectortau^i_h \cdot \vectorn\|^2_{\Sigma^{\vectoru}_N}\big), 
\label{eq:majorant-stain-divu-tensor}
\end{alignat}
where $\mathbf{r}_{\rm d, \mu, \lambda}({\vectoru}^{i}_h, \vectortau^i_h)$ in defined in 
\eqref{eq:residuals-p-iterative-lame} and $\vectortau^i_h = \vectortau^i_h$ is an auxiliary stress approximating function 
reconstructed in the correspondence with $\vectoru^i_h$.
\ProofEnd

\begin{remark}
We note the choice of the auxiliary tensor-function providing the optimal values of the error estimate is 
$\vectortau^\star_h 
:= \tensorL \strain(\tilde{\vectoru}^{i}) 
:= 2\, \mu \strain(\tilde{\vectoru}^{i}) + \lambda \dvrg (\tilde{\vectoru}^{i}) \tensorI$. 
In this case, we can show that the equilibration residual $\mathbf{r}_{\rm eq}(p^{i}_h, \vectortau^i_h)$ vanishes, 
and the dual one provides the exact representation of the error.
\end{remark}

Lemma \ref{lem:displacement-estimate} proceeds with estimation of 
$e^{i}_{\vectoru}$ (cf. \eqref{eq:eiu}), accounting for the error arising if 
\eqref{eq:mechanics-fully-discrete} is replaced by \eqref{eq:mechanics-h-perturbed}.

\begin{lemma}
\label{lem:displacement-estimate}
For any $p^i_h \in W_0$, 
any $\vectoru^i_h \in \boldsymbol{V_0}$ approximating $\tilde{\vectoru}^i$ in \eqref{eq:mechanics-h-perturbed},
and any $\vectorz^i_h \in H_{\Sigma^p_N}(\Omega, \dvrg)$ and \linebreak
$\vectortau^i_h \in [\mathcal{T}_{\rm Div}(\Omega)]^{d \times d}$, the estimate
%
\begin{equation}
|\!|\!| \, e^{i}_{\vectoru} \, |\!|\!|^2_{\vectoru} 
\leq {\overline{\rm M}^h_{u}(({\vectoru}, p)^{i}_h, (\vectortau, \vectorz)^i_h)} := 
\tfrac{2 \, \lambda \, \eta^2 \, \alpha^2}{2\, \chi \, \lambda - 1}\,
\overline{\rm M}^h_{p, L^2}(p^i_h, \vectorz^i_h)
+ 2 \, \overline{\rm M}_{\tilde{u}} (({\vectoru}, {p})^i_h, \vectortau^i_h) 
\label{eq:majorant-displpacement-energy-norm-h}
\end{equation}
holds, where $\zeta \geq 0$ and $\chi \in \big[\tfrac{1}{2\,\lambda}, +\infty)$. Here, 
$\overline{\rm M}_{p, L^2}$ and $\overline{\rm M}_{\tilde{\vectoru}}$ are 
defined in \eqref{eq:majorant-pressure-in-l2-norm} and \eqref{eq:majorant-displpacement-energy-norm}, 
respectively, and $\alpha$ and $\lambda$ are characteristics of the Biot model.
\end{lemma}

\ProofBegin
As it was noted in \eqref{eq:comb-error}, the error is two-folded and composed from  
$ \| \, {\vectoru}^{i} - \tilde{\vectoru}^{i} \|^2_{\vectoru}$ and 
$\| \, \tilde{\vectoru}^{i} - {\vectoru}^{i}_h \, \|^2_{\vectoru}$, where the second term is controlled by 
\eqref{eq:majorant-displpacement-energy-norm} in
Lemma \ref{lem:displacement-func-estimate}.
The estimate of the first term is derived by considering 
the difference of \eqref{eq:mechanics-fully-discrete} and \eqref{eq:mechanics-h-perturbed}, i.e., 
%
\begin{equation*}
2 \mu\, (\strain({\vectoru}^{i} - \tilde{\vectoru}^{i}), \strain(\vectorv)) 
+ \lambda (\dvrg ({\vectoru}^{i} - \tilde{\vectoru}^{i}), \dvrg \vectorv) 
= - \alpha (p^{i} - p^i_h, \dvrg \vectorv).
\end{equation*}
By choosing $\vectorv = {\vectoru}^{i} - \tilde{\vectoru}^{i}$, we obtain the identity
\begin{equation*}
\| \strain({\vectoru}^{i} - \tilde{\vectoru}^{i}) \|^2_{2\mu} 
+ \|\, \dvrg ({\vectoru}^{i} - \tilde{\vectoru}^{i})\, \|^2_{\lambda} 
= - \alpha (p^{i} - p^i_h, \dvrg ({\vectoru}^{i} - \tilde{\vectoru}^{i})).
\end{equation*}
The latter one can be estimated from above by the Cauchy inequality, which yields
\begin{alignat*}{2}
\| \strain(e^{i}_{\vectoru}) \|^2_{2\mu}
+ \|\, \dvrg (e^{i}_{\vectoru})\, \|^2_{\lambda}
\leq \alpha \| e_{p}^{i} \| 
\| \dvrg (e^{i}_{\vectoru}) \|.
\end{alignat*}
%
By using the Young inequality with $\chi \geq \tfrac{1}{2\lambda}$, we arrive at
\begin{equation}
\| \strain({\vectoru}^{i} - \tilde{\vectoru}^{i}) \|^2_{2\mu} 
+ (\lambda - \tfrac{1}{2\, \chi}) \,
 \|\, \dvrg ({\vectoru}^{i} - \tilde{\vectoru}^{i})\, \|^2
\leq \tfrac{\chi}{2} \alpha^2 \| p^{i} - p^i_h\|^2.
\label{eq:majorant-energy-norm-2}
\end{equation}
%
According to Lemma \ref{cor:pressure-weak-term-majorant}, the linear combination in 
\eqref{eq:majorant-energy-norm-2} can be estimated as
\begin{equation*}
\| \strain ({\vectoru}^{i} - \tilde{\vectoru}^{i}) \|^2_{2\mu} + (\lambda - \tfrac{1}{2\, \chi}) \,
 \|\, \dvrg ({\vectoru}^{i} - \tilde{\vectoru}^{i})\, \|^2
\leq  \tfrac{\chi \, \alpha^2}{2} \,\overline{\rm M}_{p, L^2}(p^i_h).
\end{equation*}
%
By using 
\begin{equation*}
\| \strain({\vectoru}^{i} - \tilde{\vectoru}^{i}) \|^2_{2\mu} 
+ \|\, \dvrg ({\vectoru}^{i} - \tilde{\vectoru}^{i})\, \|^2_{\lambda} \\
\leq \tfrac{2\, \chi \, \lambda}{2\, \chi \, \lambda - 1} \, 
\Big(\| \strain({\vectoru}^{i} - \tilde{\vectoru}^{i}) \|^2_{2\mu} 
+ (\lambda - \tfrac{1}{2\, \chi}) \, 
 \|\, \dvrg ({\vectoru}^{i} - \tilde{\vectoru}^{i})\, \|^2\Big),
\end{equation*}
we obtain 
\begin{alignat*}{2}
|\!|\!| {\vectoru}^{i} - \tilde{\vectoru}^{i} |\!|\!|^2_{\vectoru} &
\leq \tfrac{\lambda \, \chi^2 \, \alpha^2}{2\, \chi \, \lambda - 1}\,\overline{\rm M}_{p, L^2}(p^i_h).
\label{eq:last-res-u-estimated-M}
\end{alignat*}
Combining \eqref{eq:majorant-displpacement-energy-norm} and \eqref{eq:last-res-u-estimated-M}, 
we arrive at 
\begin{alignat}{2}
 \| \, \strain (e^{i}_{\vectoru}) \, \|^2_{2\mu}
+ \| \, \dvrg (e^{i}_{\vectoru}) \, \|^2_{\lambda}
\leq \tfrac{2 \, \lambda \, \chi^2 \, \alpha^2}{2\, \chi \, \lambda - 1}\,\overline{\rm M}_{p, L^2}(p^i_h)
       + 2 \, \overline{\rm M}_{\tilde{u}}({\vectoru}^{i}_h, {p}^i_h).
\end{alignat}

\ProofEnd
%
%
In addition to \eqref{eq:majorant-displpacement-energy-norm}, we can obtain the estimate 
for the error measured in terms of $\| \dvrg \cdot \|^2$-norm.

\begin{corollary}
\label{cor:displacement-div-norm}
%
For any $p^i_h \in W_0$, 
any $\vectoru^i_h \in \boldsymbol{V_0}$ approximating $\tilde{\vectoru}^i$ in \eqref{eq:mechanics-h-perturbed},
as well as any parameters and function defined 
in Lemma \ref{lem:displacement-estimate}, we have
\begin{equation}
\|\, \dvrg (e_{\vectoru}^{i})\, \|^2
\leq \overline{\rm M}^h_{u, \dvrg} (({\vectoru}, p)^{i}_h, \vectortau^i_h, \vectorz^i_h) := 
\tfrac{1}{(2\mu \, d + \lambda)} \,
\Big( 
\tfrac{2 \, \lambda \, \eta^2 \, \alpha^2}{2\, \eta \, \lambda - 1}\,\overline{\rm M}_{p, L^2}(p^i_h, \vectorz^i_h)
+ \overline{\rm M}_{\tilde{u}}(({\vectoru}, {p})^i_h, \vectortau^i_h)
\Big), 
\label{eq:majorant-dvrg-u}
\end{equation}
where $\overline{\rm M}_{p}(p^i_h, \vectorz^i_h)$ and 
$\overline{\rm M}_{\tilde{u}}(({\vectoru}, {p})^i_h, \vectortau^i_h)$ are defined in 
\eqref{eq:majorant-pressure-in-l2-norm} and \eqref{eq:majorant-displpacement-energy-norm}
for any $\vectorz^i_h \in H_{\Sigma^p_N}(\Omega, \dvrg)$ and 
$\vectortau^i_h \in [\mathcal{T}_{\rm Div}(\Omega)]^{d \times d}$, respectively, 
and $\mu$ is characteristic of the Biot model.
\end{corollary}

\ProofBegin
By using inequality \eqref{eq:inequality-divu-strain} and substituting it in \eqref{eq:majorant-displpacement-energy-norm}, 
we arrive at 
\begin{equation*}
(2\mu \, d + \lambda) \|\, \dvrg (e_{\vectoru}^{i})\, \|^2
\leq {\overline{\rm M}_{u, \dvrg}(p^i_h, \vectorz^i_h; \zeta)} 
= \tfrac{2 \, \lambda \, \eta^2 \, \alpha^2}{2\, \eta \, \lambda - 1}\,\overline{\rm M}_{p, L^2}(p^i_h)
+ 2 \, \overline{\rm M}_{\tilde{u}}(({\vectoru}, {p})^i_h, \vectortau^i_h).
\end{equation*}
%
\ProofEnd


\section{Estimates of errors generated by the iteration method}
\label{sec:iterative-estimates}

Next, we consider guaranteed bounds of errors arising in the process
of contractive iterations \eqref{eq:flow-fully-discrete}--\eqref{eq:mechanics-fully-discrete} applied
to the system \eqref{eq:mechanics-semi-discrete}--\eqref{eq:flow-semi-discrete}.

\paragraph{Error estimates for pressure term} First, we
prove the following result for the error in the flow equation, which can be done in two different ways.
For Lemma \ref{lem:iterative-pressure-estimate}, we consider functions $p^i_h, p^{i-1}_h \in W_0$ as approximations
of two consequent pressures associated with the iterations $i$ and $i-1$, whereas 
$\vectoru^i_h, \vectoru^{i-1}_h \in \boldsymbol{V_0}$ are approximations of $\tilde{u}^i$ and 
$\tilde{u}^{i-1} \in \boldsymbol{V_0}$ in \eqref{eq:mechanics-h-perturbed}, respectively. From now on, 
when we refer to both dual variables $\vectorz^i_h$ and $\vectortau^i_h$ corresponding to the $i$th iteration step, 
we refer to them as a pair $(\vectortau, \vectorz)^i_h$.

%
%



\begin{lemma}
\label{lem:iterative-pressure-estimate}
\sveta{For $p^i, p^{i-1} \in W_0$ approximating $p \in W_0$ in \eqref{eq:flow-semi-discrete},}
the estimate of the error incorporated in the pressure term on the $i$th iteration step has the following 
form
%
\begin{equation}
\begin{alignedat}{2}
|\!|\!|\, p - p^i \,|\!|\!|^2_{p} & \leq \overline{\rm M}^{\,i}_{p} 
\big((\vectoru, p)^{i-1}_h, (\vectortau, \vectorz)^{i-1}_h, (\vectoru, p)^{i}_h, (\vectortau, \vectorz)^i_h \big)
:= 
\; \tfrac{3 \, q}{1-q^2} \Bigg(\tfrac{\big(C^{\rm F}_{\Sigma^p_D}\big)^2\, \beta}{\lambda_\tensorK \, \tau} + 1 \Bigg) 
   \bigg(\| {\eta}^{i}_h - {\eta}^{i-1}_h  \|^2 \\ 
& \qquad \qquad \qquad \qquad \qquad \qquad
+ \tfrac{\lambda}{2} \,\big({\overline{\rm M}}^h_{u, \dvrg} (({\vectoru}, p)^{i}_h, (\vectortau, \vectorz)^i_h)
	   + {\overline{\rm M}}^h_{u, \dvrg}(({\vectoru}, p)^{i-1}_h, (\vectortau, \vectorz)^{i-1}_h) \big) \\ 
& \qquad \qquad \qquad \qquad \qquad \qquad
+ \tfrac{L}{4} \big({\overline{\rm M}}^h_{p, L^2} (p^i_h, \vectorz^i_h)  
                            + {\overline{\rm M}}^h_{p, L^2}(p^{i-1}_{h}, \vectorz^{i-1}) \big) \bigg),
\end{alignedat}
\label{eq:estimate-p-iterative}
\end{equation}
%
where ${\overline{\rm M}}^h_{u, \dvrg}$ and ${\overline{\rm M}}^h_{p, L^2}$ are defined 
in Corollaries \ref{cor:pressure-weak-term-majorant} and \ref{cor:displacement-div-norm} 
with $\vectortau^i_h \in [\mathcal{T}_{\rm Div}(\Omega)]^{d \times d}$ and $\vectorz^i_h \in H_{\Sigma^p_N}(\Omega, \dvrg)$, respectively, $q = \tfrac{L}{\beta + L}$, and 
$${\eta}^{i}_h = \tfrac{\alpha}{\gamma} \dvrg \vectoru^{i}_h - \tfrac{L}{\gamma} \, p^{i}_h, 
\qquad L  \geq  \tfrac{\alpha^2}{2 \,\lambda}, \qquad \forall 
p^i_h \in W_0, \;\vectoru^i_h \in \boldsymbol{V_0}.$$ 
Parameters $\alpha$, $\beta$, $\lambda$, $\mu_f$, $C^{\rm F}_{\Sigma^p_D}$, 
$\lambda_\tensorK$, and $\tau$ are characteristics of the semi-discrete Biot model 
\eqref{eq:mechanics-iterative}--\eqref{eq:flow-iterative}.
\end{lemma}

\ProofBegin 
We begin by noting that for the error $p - p^{i}$ caused by the iterative scheme
\begin{alignat}{2}
|\!|\!| \, p - p^i \,|\!|\!|^2_{p} 
& = \| \, p - p^i \, \|^2_{\beta}
+  \| \, \nabla (p - p^i) \, \|^2_{\tensorK_{\tau}}\; 
{\mnote{\eqref{eq:Poincare}}} \leq \bigg(\tfrac{(C^{\rm F}_{\Sigma^p_D})^2 \, \beta \,}{\lambda_\tensorK \, \tau} + 1\bigg) \, 
       \| \, \nabla (p - p^i) \, \|^2_{\tensorK_{\tau}}.
\label{eq:contraction-in-norm-P}
\end{alignat}
The estimate of $\|\,\nabla (p - p^i) \,\|^2_{\tensorK_{\tau}}$ 
follows from \eqref{eq:estimates-of-nabla-p-via-sigma}.
%
To proceed forward, we need to estimate the right-hand side of \eqref{eq:estimates-of-nabla-p-via-sigma}, namely 
$\| {\eta}^{i} - {\eta}^{i-1}\|^2$. By adding and extracting the discretized approximations 
${\eta}^{i-1}_h$ and ${\eta}^{i}_h$, we obtain
\begin{equation}
\| {\eta}^{i} - {\eta}^{i-1}\|^2 
\leq 3 \, \big(\| {\eta}^{i}_h - {\eta}^{i-1}_h \|^2 
   + \| {\eta}^{i}   - {\eta}^{i}_h \|^2 
	 + \| {\eta}^{i-1} - {\eta}^{i-1}_h \|^2\big). 
\label{eq:eta-i-i-1-etimate}
\end{equation}
%
Here, the first term $\| {\eta}^{i}_h - {\eta}^{i-1}_h\|^2$ is fully computable, and by means 
of relation
$${\eta}^{i} = \tfrac{1}{\gamma} (\alpha \, \dvrg \vectoru^{i} - L p^{i}),$$ 
we obtain the estimate for the second and third terms:
\begin{equation*}
\| {\eta}^{i} - {\eta}^{i}_h\|^2 
\leq \tfrac{1}{2\gamma^2} \big(\alpha^2 \| \dvrg (e_{\vectoru}^{i})\|^2 
   + L^2 \| e_p^{i} \|^2\big)
\; \mnote{\eqref{eq:majorant-pressure-in-l2-norm}, \eqref{eq:majorant-dvrg-u}} 
\leq \tfrac{1}{2\gamma^2} \Big(\alpha^2 \, {\overline{\rm M}}^h_{u, \dvrg} (p^{i}_h)
   + L^2 \, {\overline{\rm M}}^h_{p, L^2}(p^{i}_h)\Big).
\end{equation*}
%
For simplicity, we exclude parameter $\gamma$ by substituting $\gamma^2 = 2 \, L$: 
\begin{equation}
\| {\eta}^{i-1} - {\eta}^{i-1}_h \|^2
\leq \tfrac{1}{4\, L} \big(\alpha^2 \, {\overline{\rm M}}^h_{u, \dvrg}(p^{i-1}_h) 
                            + L^2 \, {\overline{\rm M}}^h_{p, L^2}(p^{i-1}_h)\big).
\label{eq:estimate-of-p}
\end{equation}
%
Therefore, the estimate of $|\!|\!|\, p - p^i \,|\!|\!|^2_{p}$ can be 
represented as follows
\begin{alignat}{2}
|\!|\!|\, p - p^i \,|\!|\!|^2_{p} 
& \leq \overline{\rm M}^{\,i}_P 
:=  \Big(\tfrac{(\CF)^2\, \beta}{\lambda_\tensorK \, \tau} + 1\Big)\, \tfrac{3\,q}{1-q^2}
   \bigg\{ \| {\eta}^{i}_h - {\eta}^{i-1}_h \|^2 \nonumber\\
& \qquad \quad
+ \tfrac{1}{4\, L} 
	\Big( \alpha^2 \, \big({\overline{\rm M}}^h_{u, \dvrg}(p^{i}_h) 
	                + {\overline{\rm M}}^h_{u, \dvrg}(p^{i-1}_h) \big) 
+ L^2 \big({\overline{\rm M}}^h_{p, L^2}(p^{i}_h) + {\overline{\rm M}}^h_{p, L^2}(p^{i-1}_h) \big) \Big)
\bigg\}.
\label{eq:estimate-p-iterative-pre}
\end{alignat}
Finally, by substituting $\lambda = \tfrac{\alpha^2}{2\, L}$ in \eqref{eq:estimate-p-iterative-pre}, 
we arrive at \eqref{eq:estimate-p-iterative}.
\ProofEnd
\sveta{
\begin{corollary}
\label{cor:iterative-pressure-estimate}
Let conditions of Lemma \ref{eq:lemma-2} and Lemma \ref{lem:iterative-pressure-estimate} hold.
Then, for $1 \leq m \leq i$, the estimate of the error incorporated in the pressure term on the $i$th iteration step has
the alternative form
%
\begin{equation}
\begin{alignedat}{2}
|\!|\!|\, p - p^i \,|\!|\!|^2_{p} & \leq \overline{\rm M}^{\,i, m}_{p}
\big((\vectoru, p)^{i-m}_h, (\vectortau, \vectorz)^{i-m}_h, (\vectoru, p)^{i}_h, (\vectortau, \vectorz)^i_h \big)
:= 
\; \tfrac{3 \, q^m}{1-q^{2m}}
\Bigg(\tfrac{\big(C^{\rm F}_{\Sigma^p_D}\big)^2\, \beta}{\lambda_\tensorK \, \tau} + 1 \Bigg)
\bigg(\| {\eta}^{i}_h - {\eta}^{i-m}_h  \|^2 \\
& \qquad \qquad \qquad \qquad \qquad \qquad
+ \tfrac{\lambda}{2} \,\big({\overline{\rm M}}^h_{u, \dvrg} (({\vectoru}, p)^{i}_h, (\vectortau, \vectorz)^i_h)
+ {\overline{\rm M}}^h_{u, \dvrg}(({\vectoru}, p)^{i-m}_h, (\vectortau, \vectorz)^{i-m}_h) \big) \\
& \qquad \qquad \qquad \qquad \qquad \qquad
+ \tfrac{L}{4} \big({\overline{\rm M}}^h_{p, L^2} (p^i_h, \vectorz^i_h)
+ {\overline{\rm M}}^h_{p, L^2}(p^{i-m}_{h}, \vectorz^{i-m}) \big) \bigg).
\end{alignedat}
\label{eq:estimate-p-iterative-m}
\end{equation}
\end{corollary}
}

\authors{
The norm in the left-hand side of \eqref{eq:eta-i-i-1-etimate} can be alternatively estimated by the contractive 
estimate, which results into alternative error bound independent of the functional majorants.  This result is 
presented below.
\begin{lemma}
\label{lem:iterative-pressure-estimate-2}
For any $p^i \in W_0$ approximating $p \in W_0$ in \eqref{eq:flow-semi-discrete},
the estimate of the error incorporated in the pressure term on the $i$th iteration step has the following 
form
%
\begin{equation}
\begin{alignedat}{2}
|\!|\!|\, p - p^i \,|\!|\!|^2_{p} & \leq \reallywidetilde{\rm M}^{\,i}_{p} 
\big({\eta}^{1}_h, {\eta}^{0}_h \big)
:=  \bigg(\tfrac{\big(C^{\rm F}_{\Sigma^p_D}\big)^2 \, \beta \,}{\lambda_\tensorK \, \tau} + 1\bigg) \, 
      \tfrac{3 \, q^{2i - 1}}{1-q^2} \, \\
& \qquad \qquad   \qquad \qquad    \qquad \qquad     
\Big( 
      \frac{\alpha^2}{\gamma^2} \overline{\rm M}^h_{p, L^2}(p^1_h) 
      + \frac{L^2}{\gamma^2} \, \overline{\rm M}^h_{u, \dvrg} (({\vectoru}, p)^{1}_h)
      + \| {\eta}^{0} - {\eta}^{0}_h\|^2 + \| {\eta}^{1}_h - {\eta}^{0}_h\|^2\Big),
\end{alignedat}
\label{eq:estimate-p-iterative-2}
\end{equation}
%
where 
$$q = \tfrac{L}{\beta + L}, \qquad {\eta}^{i}_h = \tfrac{\alpha}{\gamma} \dvrg \vectoru^{i}_h - \tfrac{L}{\gamma} \, p^{i}_h, 
\qquad L \geq \tfrac{\alpha^2}{2 \,\lambda}, \qquad \forall 
p^i_h \in W_{0h}, \;\vectoru^i_h \in \boldsymbol{V_{0h}},$$ 
and ${\overline{\rm M}}^h_{u, \dvrg}$ and ${\overline{\rm M}}^h_{p, L^2}$ are defined 
in Corollaries \ref{cor:pressure-weak-term-majorant} and \ref{cor:displacement-div-norm}.
Parameters $\alpha$, $\beta$, $\lambda$, $\mu_f$, $C^{\rm F}_{\Sigma^p_D}$, 
$\lambda_\tensorK$, and $\tau$ are characteristics of the semi-discrete Biot model 
\eqref{eq:flow-iterative} -- \eqref{eq:mechanics-iterative}.
\end{lemma}

\ProofBegin
We use inequality \eqref{eq:contraction-in-norm-P}, to obtain
\begin{alignat}{2}
|\!|\!| \, p - p^i \,|\!|\!|^2_{p} 
& \leq \bigg(\tfrac{(C^{\rm F}_{\Sigma^p_D})^2 \, \beta \,}{\lambda_\tensorK \, \tau} + 1\bigg) \, 
      \tfrac{q}{1-q^2} \, \| {\eta}^{i} - {\eta}^{i-1}\|^2.
\label{eq:contraction-in-norm-P-2}
\end{alignat}
}
\authors{
However, in this proof, we use contraction properties of sequence $\{{\delta \eta}^{i}\}$, to obtain
\begin{alignat}{2}
\| {\eta}^{i} - {\eta}^{i-1}\|^2 
 \leq q^{2(i-1)} \, \| {\eta}^{1} - {\eta}^{0}\|^2
& \leq q^{2(i-1)} \, (3 \| {\eta}^{1} - {\eta}^{1}_h\|^2 + 3 \| {\eta}^{0} - {\eta}^{0}\|^2 + 3 \| {\eta}^{1}_h - {\eta}^{0}_h\|^2) 
\label{eq:contraction-eta-0-and-1}
\end{alignat}
where all terms on the right-hand side except  $\| {\eta}^{1} - {\eta}^{1}_h\|^2$ are computable. 
As in Lemma \ref{lem:pressure-full-norm-majorant}, $\| {\eta}^{1} - {\eta}^{1}_h\|^2$ is controlled by combination of 
$\overline{\rm M}^h_{p, L^2}(p^1_h, \vectorz^1_h; \zeta)$ and 
$\overline{\rm M}^h_{u, \dvrg} (({\vectoru}, p)^{1}_h, \vectortau^1_h, \vectorz^1_h)$, i.e., 
\begin{alignat}{2}
\| {\eta}^{1} - {\eta}^{1}_h\|^2 
\leq \frac{\alpha^2}{\gamma^2} \, \| \, \dvrg (\vectoru^{1} - \vectoru^{1}_h) \, \|^2
   + \frac{L^2}{\gamma^2} \, \| \,  p^{1} - p^{1}_h \, \|^2 
\leq \frac{\alpha^2}{\gamma^2} \overline{\rm M}^h_{p, L^2}(p^1_h, \vectorz^1_h; \zeta) 
      + \frac{L^2}{\gamma^2} \, \overline{\rm M}^h_{u, \dvrg} (({\vectoru}, p)^{1}_h, \vectortau^1_h, \vectorz^1_h),
\label{eq:etimate-eta1-etah-1-square}
\end{alignat}
where right-hand side contains computable 
$\eta^0 := \tfrac{\alpha}{\gamma} \dvrg \vectoru^{0} - \tfrac{L}{\gamma} \, p^{0}$.
Combining \eqref{eq:contraction-in-norm-P-2},  \eqref{eq:contraction-eta-0-and-1}, and \eqref{eq:etimate-eta1-etah-1-square}, 
we arrive at
\begin{alignat*}{2}
|\!|\!| \, p - p^i \,|\!|\!|^2_{p} 
& \leq \bigg(\tfrac{(C^{\rm F}_{\Sigma^p_D})^2 \, \beta \,}{\lambda_\tensorK \, \tau} + 1\bigg) \, 
      \tfrac{3 \, q^{2i - 1}}{1-q^2} \, \Big( 
      \frac{\alpha^2}{\gamma^2} \overline{\rm M}^h_{p, L^2}(p^1_h) 
      + \frac{L^2}{\gamma^2} \, \overline{\rm M}^h_{u, \dvrg} (({\vectoru}, p)^{1}_h)
      + \, \| {\eta}^{0} - {\eta}^{0}_h\|^2 
      + \| {\eta}^{1}_h - {\eta}^{0}_h\|^2\Big).
\end{alignat*}
\ProofEnd
}

\paragraph{Error estimates for displacement term}
To derive the upper bound for the error $\vectoru - \vectoru^{i}$ measured
in terms of
\begin{equation}
|\!|\!|\, \vectoru - \vectoru^{i} \,|\!|\!|^2_{\vectoru} :=
\| \, \strain (\vectoru - \vectoru^{i}) \, \|^2_{2\mu} 
+ \| \, \dvrg (\vectoru - \vectoru^{i}) \, \|^2_{\lambda}, 
\label{eq:dispacement-error-k-iteration}
\end{equation}
we exploit the idea analogous to the one used to estimate the error in the pressure term. %
\begin{lemma}
\label{lem:iterative-displacement-estimate}
%
For any ${\vectoru}^i, {\vectoru}^{i-1} \in \boldsymbol{V_0}$ approximating ${\vectoru} \in \boldsymbol{V_0}$ in 
\eqref{eq:mechanics-semi-discrete}, the error in the displacement 
on the $i$th iteration step has the following form 
%
%
\begin{equation}
\begin{alignedat}{2}
|\!|\!|\, {\vectoru} - {\vectoru}^i \,|\!|\!|^2_{\vectoru} & \leq {\overline{\rm M}}^{i}_{u} 
\big((\vectoru, p)^{i-1}_h, (\vectortau, \vectorz)^{i-1}_h, (\vectoru, p)^{i}_h, (\vectortau, \vectorz)^i_h \big)
 := (1 + \tfrac{d \, \lambda}{2\mu}) \tfrac{3 \, q^2}{1-q^2} 
\Big(\| {\eta}^{i}_h - {\eta}^{i-1}_h \| \\
& \qquad \qquad \qquad \qquad \qquad \qquad
+ \tfrac{\lambda}{2} \, \big({\overline{\rm M}}^h_{u, \dvrg} ( ({\vectoru}, p)^{i}_h, (\vectortau, \vectorz)^i_h)
	   + {\overline{\rm M}}^h_{u, \dvrg}(({\vectoru}, p)^{i-1}_h, (\vectortau, \vectorz)^{i-1}_h) \big) \\
& \qquad \qquad \qquad \qquad \qquad \qquad
+ \tfrac{L}{4}  \big({\overline{\rm M}}^h_{p, L^2}(p^{i}_h,  \vectorz^i_h) 
                             + {\overline{\rm M}}^h_{p, L^2}(p^{i-1}_h,  \vectorz^{i-1})\big) \Big),
\label{eq:estimate-u-iterative}
\end{alignedat}
\end{equation}
where ${\overline{\rm M}}^h_{u, \dvrg}$ and ${\overline{\rm M}}^h_{p, L^2}$ are defined 
in Corollaries \ref{cor:pressure-weak-term-majorant} and \ref{cor:displacement-div-norm} for 
$\vectorz^i_h, \vectorz^{i-1}_h \in H_{\Sigma^p_N}(\Omega, \dvrg)$ and 
$\vectortau^i_h, \vectortau^{i-1}_h \in [\mathcal{T}_{\rm Div}(\Omega)]^{d \times d}$, respectively, 
$q = \tfrac{L}{\beta + L}$, and 
$${\eta}^{i}_h = \tfrac{\alpha}{\gamma} \dvrg \vectoru^{i}_h - \tfrac{L}{\gamma} \, p^{i}_h, 
\qquad L = \tfrac{\alpha^2}{2 \,\lambda}, \qquad \forall 
p^i_h \in W_0, \, 
\vectoru^i_h \in \boldsymbol{V_0}.$$ 
Parameters $\lambda$, $\mu$, $\alpha$ are characteristics of the Biot model.
\end{lemma}
%
\ProofBegin
We consider 
\begin{equation}
|\!|\!|\, {\vectoru} - {\vectoru}^i \,|\!|\!|^2_{\vectoru} \;
\leq 
\big( 2\mu + d \, \lambda \big) \, 
\| \, \strain (e_{\vectoru}) \, \|^2,
\label{eq:contrative-estimates-for-displacement-norm}
\end{equation}
where the right-hand side is controlled by the contractive term $\|{\eta} - {\eta}^i \|^2$, 
which follows from \eqref{eq:estimates-of-strain-via-sigma}.
Therefore, we obtain 
%
\begin{alignat}{2}
|\!|\!|\, {\vectoru} - {\vectoru}^i \,|\!|\!|^2_{\vectoru}
& \leq \tfrac{2\mu + d \, \lambda}{2\mu} \tfrac{q^2}{1 - q^2}
\| {\eta}^{i} - {\eta}^{i-1} \|^2 \nonumber \\
& \leq \tfrac{2\mu+ d \, \lambda}{2\mu} \tfrac{3 \, q^2}{1-q^2} 
\big(\|{\eta}^{i} - {\eta}^{i}_h\|^2 
+ \|{\eta}^{i-1} - {\eta}^{i-1}_h\|^2 
+ \|{\eta}^{i}_h - {\eta}^{i-1}_h\|^2\big).
\label{eq:contractive-estimates-for-displacement-norm-via-stress}
\end{alignat}
%
Analogously to the proof of Lemma \ref{lem:iterative-pressure-estimate} (cf. \eqref{eq:estimate-of-p}), 
the estimate for the second term in \eqref{eq:contractive-estimates-for-displacement-norm-via-stress} 
results into 
\begin{alignat}{2}
|\!|\!|\, {\vectoru} - {\vectoru}^i \,|\!|\!|^2_{\vectoru} &
\leq {\overline{\rm M}}^{i}_{u} 
:= (1 + \tfrac{d \, \lambda}{2\mu})\,  \tfrac{3 \, q^2}{1-q^2} 
\Big(\| {\eta}^{i}_h - {\eta}^{i-1}_h \| \nonumber\\
& + \tfrac{1}{4\, L}
	\Big(\alpha^2 \, \big({\overline{\rm M}}^h_{u, \dvrg}(p^{i}_h, \vectorz^i_h) 
                          + {\overline{\rm M}}^h_{u, \dvrg}(p^{i-1}_h,  \vectorz^{i-1})\big)
+ L^2 \big({\overline{\rm M}}^h_{p, L^2}(p^{i}_h,  \vectorz^i_h) 
                             + {\overline{\rm M}}^h_{p, L^2}(p^{i-1}_h,  \vectorz^{i-1})\big) \Big) \Big).
\label{eq:estimate-u-iterative-pre}
\end{alignat}
Again, by substituting $\lambda = \tfrac{\alpha^2}{2\,L}$ in \eqref{eq:estimate-u-iterative-pre}, 
we arrive at \eqref{eq:estimate-u-iterative}.
\ProofEnd
\sveta{
\begin{corollary}
\label{cor:iterative-displacement-estimate}
Let conditions of Lemma \ref{eq:lemma-2} and Lemma \ref{lem:iterative-displacement-estimate} hold.
Then, for $1 \leq m \leq i$, the error in the displacement on the $i$th iteration step has the alternative form
%
%
\begin{equation}
\begin{alignedat}{2}
|\!|\!|\, {\vectoru} - {\vectoru}^i \,|\!|\!|^2_{\vectoru} & \leq {\overline{\rm M}}^{i, m}_{u}
\big((\vectoru, p)^{i-m}_h, (\vectortau, \vectorz)^{i-m}_h, (\vectoru, p)^{i}_h, (\vectortau, \vectorz)^i_h \big)
:= (1 + \tfrac{d \, \lambda}{2\mu}) \tfrac{3 \, q^{2m}}{1-q^{2m}}
\Big(\| {\eta}^{i}_h - {\eta}^{i-m}_h \| \\
& \qquad \qquad \qquad \qquad \qquad \qquad
+ \tfrac{\lambda}{2} \, \big({\overline{\rm M}}^h_{u, \dvrg} ( ({\vectoru}, p)^{i}_h, (\vectortau, \vectorz)^i_h)
+ {\overline{\rm M}}^h_{u, \dvrg}(({\vectoru}, p)^{i-m}_h, (\vectortau, \vectorz)^{i-m}_h) \big) \\
& \qquad \qquad \qquad \qquad \qquad \qquad
+ \tfrac{L}{4}  \big({\overline{\rm M}}^h_{p, L^2}(p^{i}_h,  \vectorz^i_h)
+ {\overline{\rm M}}^h_{p, L^2}(p^{i-m}_h,  \vectorz^{i-m})\big) \Big).
\label{eq:estimate-u-iterative}
\end{alignedat}
\end{equation}
\end{corollary}
}
\authors{
Analogously, contraction properties can be exploited to derive an alternative bound of the error
$\vectoru - \vectoru^{i}$.
\begin{lemma}
\label{lem:iterative-displacement-estimate-2}
%
For any ${\vectoru}^i \in \boldsymbol{V_0}$ approximating ${\vectoru} \in \boldsymbol{V_0}$ in 
\eqref{eq:mechanics-semi-discrete}, the error in the displacement 
on the $i$th iteration step has the following form 
%
%
\begin{equation}
\begin{alignedat}{2}
|\!|\!|\, {\vectoru} - {\vectoru}^i \,|\!|\!|^2_{\vectoru} &
\leq {\reallywidetilde{\rm M}}^{i}_{u} \big({\eta}^{1}_h, {\eta}^{0}_h\big)
:= \tfrac{2\mu + d \, \lambda}{2\mu} \tfrac{3\, q^{2i}}{1 - q^2} \,
\Big( \frac{\alpha^2}{\gamma^2} \overline{\rm M}^h_{p, L^2}(p^1_h) 
      + \frac{L^2}{\gamma^2} \, \overline{\rm M}^h_{u, \dvrg} (({\vectoru}, p)^{1}_h) 
      + \, \| {\eta}^{0} - {\eta}^{0}_h\|^2 + \| {\eta}^{1}_h - {\eta}^{0}_h\|^2\Big),
\label{eq:estimate-u-iterative-2}
\end{alignedat}
\end{equation}
where 
$$q = \tfrac{L}{\beta + L}, \qquad {\eta}^{i}_h = \tfrac{\alpha}{\gamma} \dvrg \vectoru^{i}_h - \tfrac{L}{\gamma} \, p^{i}_h, 
\qquad L = \tfrac{\alpha^2}{2 \,\lambda}, \qquad \forall 
p^i_h \in W_0, \, 
\vectoru^i_h \in \boldsymbol{V_0},$$ 
and ${\overline{\rm M}}^h_{u, \dvrg}$ and ${\overline{\rm M}}^h_{p, L^2}$ are defined 
in Corollaries \ref{cor:pressure-weak-term-majorant} and \ref{cor:displacement-div-norm}.
Parameters $\lambda$, $\mu$, $\alpha$ are characteristics of the Biot model.
\end{lemma}
\ProofBegin
The derivation comprises of  combining 
\begin{alignat}{2}
|\!|\!|\, {\vectoru} - {\vectoru}^i \,|\!|\!|^2_{\vectoru}
& \leq \tfrac{2\mu + d \, \lambda}{2\mu} \tfrac{q^2}{1 - q^2}
\| {\eta}^{i} - {\eta}^{i-1} \|^2,
\label{eq:contractive-estimates-for-displacement-norm-via-stress-2}
\end{alignat}
\eqref{eq:contraction-eta-0-and-1}, and \eqref{eq:etimate-eta1-etah-1-square}.
\ProofEnd
}



\paragraph{General estimate for the error in the iterative coupling scheme}
Derivation of the reliable estimate for the error in the pressure approximation reconstructed on the $i$th
iteration is based on the combination of two different approaches, i.e., estimates for contractive mapping 
used for iterative methods and functional error estimates. 
Theorem \ref{th:pressure-displacement-estimate} presents an upper bound of the error in the pressure term, 
which combines the above-mentioned techniques. 

\begin{theorem}
\label{th:pressure-displacement-estimate}
%
For any $p^i_h \in W_0$ and $\vectoru^{i}_h \in \boldsymbol{V_0}$ that form a contraction relative to 
$p^{i-1}_h \in W_0$ and $\vectoru^{i-1}_h \in \boldsymbol{V_0}$, we have the estimates 
\authors{
\begin{alignat*}{2}
|\!|\!|\, {e}_p \,|\!|\!|^2_{p} &
\leq \overline{\rm M}_p 
:= 2\, \bigg( \overline{\rm M}^{h}_{p} (p^i_h, \vectorz^i_h) 
 + \min \Big\{
\overline{\rm M}^{i}_{p}, \overline{\rm M}^{i, m}_{p}, \reallywidetilde{\rm M}^{i}_{p} \Big\} \bigg), \\
%
|\!|\!|\, {e}_{\vectoru} \,|\!|\!|^2_{\vectoru} & \leq 
\overline{\rm M}_u 
:= 2\, \bigg({\overline{\rm M}^h_{u} \big(({\vectoru}, p)^{i}_h, (\vectortau, \vectorz)^i_h\big)}
+
\min \Big\{
{\overline{\rm M}}^{i}_{u}, {\overline{\rm M}}^{i, m}_{u}, \reallywidetilde{\rm M}^{i}_{u}
\Big\}\bigg).
\end{alignat*}
Here, $\overline{\rm M}^{h}_{p}$, $\overline{\rm M}^{i}_{p}$, $\overline{\rm M}^{i, m}_{p}$, and $\reallywidetilde{\rm M}^{i}_{p}$
are defined in
Lemmas \ref{lem:pressure-full-norm-majorant},
\ref{lem:iterative-pressure-estimate},
Corollary \ref{cor:iterative-pressure-estimate}, and Lemma
\ref{lem:iterative-pressure-estimate-2},
whereas ${\overline{\rm M}}^{h}_{u}$, ${\overline{\rm M}}^{i}_{u}$, and $\reallywidetilde{\rm M}^{i}_{u}$ are presented by 
Lemmas \ref{lem:displacement-estimate}, \ref{lem:iterative-displacement-estimate},
Corollary \ref{cor:iterative-displacement-estimate}, and
Lemma \ref{lem:iterative-displacement-estimate-2}, respectively,
}
where
$p^{i-1}_h \in W_0$, $\vectoru^{i-1}_h \in \boldsymbol{V_0}$, 
$\vectorz^i_h, \vectorz^{i-1} \in H_{\Sigma^p_N}(\Omega, \dvrg)$, 
$\vectortau^i_h, \vectortau^{i-1} \in [\mathcal{T}_{\rm Div}(\Omega)]^{d \times d}$,
and parameter $\zeta \geq 0$.
\end{theorem}

\ProofBegin
To decompose the error  $|\!|\!|\, {e}_p \,|\!|\!|^2_{p}$ in two parts, we apply 
the triangle inequality
\begin{equation}
	|\!|\!|\, {e}_p \,|\!|\!|^2_{p}  = |\!|\!|\, p - p_h^{i} \,|\!|\!|^2_{p}
	\leq 
	2\, \big( |\!|\!|\, p - p^{i} \,|\!|\!|^2_{p} 
	+ |\!|\!|\, p^{i} - p_h^{i} \,|\!|\!|^2_{p} \big).
	\label{eq:estimate-energy-pressure}
\end{equation}
The first term in the right-hand side of \eqref{eq:estimate-energy-pressure} is bounded by 
\eqref{eq:estimate-p-iterative} from Lemma \ref{lem:iterative-pressure-estimate}, whereas 
the second term is controlled by \eqref{eq:majorant-pressure-energy-norm} from Lemma \ref{lem:pressure-full-norm-majorant}.
%
Analogously, by using the triangle rule, we obtain
\begin{equation}
|\!|\!|\, {e}_{\vectoru} \,|\!|\!|^2_{\vectoru} = 
|\!|\!|\, \vectoru - \vectoru_h^{i} \,|\!|\!|^2_{\vectoru} \leq 
2 \,\big( |\!|\!|\, \vectoru - \vectoru^{i} \,|\!|\!|^2_{\vectoru} + |\!|\!|\, \vectoru^i - \vectoru_h^{i} \,|\!|\!|^2_{\vectoru} \big).
\label{eq:error-inthree-parts}
\end{equation}
%
The first term in the right-hand side of \eqref{eq:error-inthree-parts} is controlled by \eqref{eq:estimate-u-iterative}, 
whereas the estimate of the second term follows from 
\eqref{eq:majorant-displpacement-energy-norm}.
\ProofEnd

\begin{theorem}
\label{th:total-estimate}
For any $p^i_h \in W_0$ and $\vectoru^{i}_h \in \boldsymbol{V_0}$, 
we have the estimates 
\begin{alignat*}{2}
\big|\!\big[ {e}_p, {e}_{\vectoru} \big]\!\big| \leq \overline{\rm {M}} :=
\overline{\rm M}_p + \overline{\rm M}_u,
\end{alignat*}
where $\overline{\rm M}_p$ and $\overline{\rm M}_u$ are defined in 
Theorem \ref{th:pressure-displacement-estimate}.
\end{theorem}

\section{Numerical example}
\label{sec:numerical-example}

\authors{
The numerical properties of above estimates are explained in the series of the following numerical examples. We take two tests where we begin with a manufactured solution for the pressure and displacement. These tests pertain to different values of the contraction coefficient $q$  in the fixed stress scheme and study the efficiency of the error estimators.  For each of the two tests, we consider both parameters that are academic as well as repeat them for more realistic parameters of the material properties (permeability, Lame coefficients etc). We solve the Biot model at each time step using the fixed stress scheme by choosing different mesh sizes and  time steps. Moreover, we also vary the finite element pairs used for solving the problem. All these variations, including material properties, discretization parameters, and the finite element pairs allow us to study the performance of the estimators and show their robustness. 

Our approach is to compute the approximations of the pressure and displacement for the test cases and study the convergence of this to the exact solution. The errors are computed by employing different norms including the combination of pressure and displacement as well the individual ones in both $L^2$ and the energy norms.  We choose a time and a spatial mesh size and at this discrete time, we fix the number of fixed stress iterations and study the convergence.  The majorants $\overline{\rm M}^h_{p}$ and $\overline{\rm M}^h_{u}$ are obtained via minimisations w.r.t. to the auxiliary functions. To characterise the efficiency of all the above mentioned error estimatators, we use the efficiency index defined as 
$\Ieff := \overline{\rm M} / |\!|\!| e |\!|\!|^2$, where $\overline{\rm M}$ is a chosen majorant and $|\!|\!| e |\!|\!|^2$ 
is error measured in the corresponding error norm. We study the individual contributions of the different terms in the majorants, namely the dual term and  reliability/equilibration term. The first one allows us to study the efficiency as the error indicators and the small contribution of the second one ensures the reliability of the total error bounds. 
This is studied for the different majorants, the one for the discretization error and for the iterative ones and for the different choices of error norms.  
}
\begin{example}
\label{ex:1} 
\rm
{\bf Simplified parameters}.
Let $\Omega := (0, 1)^2 \in \Rtwo$, $T = 10$. 
The exact solution 
of \eqref{eq:poroelastic-system} is defined as
$$\vectoru(x, y, t) := t \, x \, (1 - x)\, y \, (1 - y)\, [1, 1]^{\rm T} \quad 
\mbox{and} \quad
p(x, y, t) := t \, x \, (1 - x)\, y \, (1 - y).$$ 
The Lame parameters are  $\mu = \lambda_{\rm vol} = 1$, which leads to 
$\lambda_{\rm plane} = \tfrac{2 \, \lambda_{\rm vol} \, \mu}{\lambda_{\rm vol} + 2 \, \mu} = \tfrac{2}{3}$. 
We set $\alpha = \beta = 1$, $C_{\rm F} = \tfrac{1}{\sqrt{2}\, \pi}$, and $\tensorK$ $[\rm mD/cP]$ is a unit tensor. From 
the parameters above, it follows that $L =\tfrac{\alpha^2}{2 \, (\lambda + 2 \mu /d)} =  0.3$ and 
$q = \tfrac{L}{\beta + L} = \tfrac{3}{13} \approx 0.23077$. 
Taking into account that the error estimates depend on the term $\tfrac{q^2}{(1 - q^2)}$, which goes 
to infinity as $q$ goes to $1$, the efficiency of the resulting 
$\overline{\rm M}_u$ and $\overline{\rm M}_p$ drastically depends on the value of $q$. 
In this particular example, the ratio $\tfrac{q^2}{(1 - q^2)} = 0.05625$ is relatively small, which 
prevents ${\overline{\rm M}}^{i}_{u}$ and ${\overline{\rm M}}^{i}_{p}$ overestimating the errors. 

\begin{table}[!ht]
\footnotesize
\scriptsize
\centering
\newcolumntype{g}{>{\columncolor{gainsboro}}c} 	
\newcolumntype{k}{>{\columncolor{lightgray}}c} 	
\newcolumntype{s}{>{\columncolor{silver}}c} 
\newcolumntype{a}{>{\columncolor{ashgrey}}c}
\newcolumntype{b}{>{\columncolor{battleshipgrey}}c}
\begin{tabular}{c|cc|cc|cc|cc}
\parbox[c]{1.6cm}{\centering $i = 1, \ldots, I$ } & 
\parbox[c]{1.4cm}{\centering $|\!|\!| e_p |\!|\!|^2$ }  & 	  
\parbox[c]{1.4cm}{\centering $\overline{\rm M}^h_{p}$ } &    
\parbox[c]{1.0cm}{\centering $\| e_p \|^2_{\beta}$ } & 
\parbox[c]{1.4cm}{\centering $\overline{\rm M}^h_{p, L^2}$ } & 
\parbox[c]{1.4cm}{\centering $|\!|\!| e_u |\!|\!|^2$ } & 
\parbox[c]{1.4cm}{\centering $\overline{\rm M}^h_{u}$  }&    
\parbox[c]{1.0cm}{\centering $\| \dvrg e_u \|^2_{\lambda}$ }   & 	
\parbox[c]{1.4cm}{\centering $\overline{\rm M}^h_{u, \dvrg}$ } \\
\midrule
\multicolumn{9}{l}{(a) $N = 10$, \; \quad $[t_{9}, t_{10}]$, $\tau = 1.0$} \\
\midrule
     1 & 1.2868e-02  &   $-$ &     2.4070e-04  &   $-$ & 2.0722e-04  &   $-$  &     5.1044e-05   &  $-$ \\
     2 & 1.8733e-04  &   1.8732e-04 &     4.7759e-09  &   8.9032e-06 & 1.8687e-04  &   4.9549e-04  &     4.4544e-05   &  1.0618e-04 \\
     3 & 1.8731e-04  &   $-$ &     2.0666e-09  &   $-$ & 1.8687e-04  &   $-$  &     4.4544e-05  &   $-$ \\ 
     4 & 1.8731e-04  &  1.8738e-04 &     2.0481e-09 &    8.9060e-06 & 1.8687e-04  &   4.9549e-04  &     4.4544e-05   &  1.0618e-04 \\
     5 & 1.8731e-04  &   1.8738e-04 &     2.0479e-09  &   8.9060e-06 & 1.8687e-04  &   4.9549e-04  &     4.4544e-05  &   1.0618e-04 \\
\midrule
\multicolumn{9}{l}{(b) $N = 10^2$, \; \quad $[t_{99}, t_{100}]$, $\tau = 0.1$} \\
\midrule
     1 & 7.0040e-01  &   $-$ &     1.1032e-01   &  $-$ & 1.5274e-03   &   $-$  &     4.6795e-04  &   $-$ \\
     2 & 1.7443e-04  &    1.2920e-04 &     1.0388e-05   &  3.9462e-05 & 1.8690e-04   &   5.0505e-04  &     4.4558e-05  &   1.0823e-04 \\
     3 & 1.3128e-04  &   $-$ &     8.8865e-08   &  $-$ & 1.8687e-04   &   $-$  &     4.4544e-05  &  $-$ \\
     4 & 1.3113e-04  &    1.3111e-04 &     1.7743e-08   &   4.0047e-05 & 1.8687e-04  &    4.9275e-04  &     4.4544e-05  &   1.0559e-04 \\
     5 & 1.3113e-04  &   1.3112e-04 &     1.5374e-08   &   4.0049e-05 & 1.8687e-04 &    4.9276e-04  &     4.4544e-05 &    1.0559e-04 \\
\end{tabular}
\caption{{\it Example \ref{ex:1}}. Errors and majorants w.r.t. the iteration steps for $h = \tfrac{1}{64}$ and $I = 5$
(both values are measured relative to the increment in $|\!|\!| p |\!|\!|^2_p$ and $|\!|\!| \vectoru |\!|\!|^2_{u}$ on the $N$th time step).}
\label{tab:example-1-iteration-convergence-relative}
\end{table}
First, we generate approximations employing $10$ time steps (where the total number of steps in time is denoted by
$N$) with corresponding size of the step $\tau$ = 1.0 and spatial mesh-size $h = \tfrac{1}{64}$ using standard $\mathds{P}_1$
(polynomial/lagrangian first-order) finite elements (see \eqref{eq:p0-p1}). Let $I$ denote the number of iterations
to solve \eqref{eq:flow-fully-discrete}--\eqref{eq:mechanics-fully-discrete} on each time-step; we consider $I = 5$
for this discretization.
Table \ref{tab:example-1-iteration-convergence-relative} illustrates the convergence of the errors in $\vectoru$ and
$p$ w.r.t. the iteration steps $i = 1, \ldots, I$ for the time-interval $[t_9, t_{10}] = [9.0, 10.0]$.

The values of $\overline{\rm M}^h_{p}$ and $\overline{\rm M}^h_{u}$ are obtained by minimization of each functional
w.r.t. to auxiliary functions. For instance, the functional estimate to control the error in variable $p$
\begin{equation}
\overline{\rm M}^h_{p} (p^i_h, \vectorz^i_h; \zeta) := 
(1 + \zeta) \, \| \mathbf{r}_{\rm d} (p^i_h, \vectorz^i_h) \|^2_{\tensorK^{-1}_{\tau}} 
+ (1 + \tfrac{1}{\zeta}) \, C^p_\Omega
\Big( 2\, \|\mathbf{r}_{\rm eq} (p^i_h, \vectorz^i_h) \|^2_{\Omega}
+ 2\, \gamma^2 \, q^{2(i-1)} \, \|{\eta}^{0} - {\eta}^{0}_h\|^2
+ \| \vectorz^i_h \cdot \vectorn \|^2_{\Sigma^{p}_N} \Big)
\end{equation}
is minimized w.r.t. the function $\vectorz^i_h$ and parameter $\zeta$. The results of this optimization procedure are
listed in Table \ref{tab:example-1-majorant-optimization}. We see that only two iterations are enough to achieve good
efficiency $\Ieff := \overline{\rm M}^h_{p} / |\!|\!| e_p |\!|\!|^2 = 1.0023$.
In Table \ref{tab:example-1-majorant-optimization}, we see the decrease of the dual term
$\overline{\rm m}^2_{\rm d}
:= \| \mathbf{r}_{\rm d} (p^i_h, \vectorz^i_h) \|^2_{\tensorK^{-1}_{\tau}} 
= \| \vectorz^i_h - \tensorK_{\tau} \nabla p^i_h \|^2_{\tensorK^{-1}_{\tau}}$
as well as the reliability/equilibration term of the majorant
$\overline{\rm m}^2_{\rm eq}
:= \| \mathbf{r}_{\rm eq} (p^i_h, \vectorz^i_h) \|^2_{\Omega}
= \| \widetilde{g} - \gamma \, {\eta}^{i-1}_h - (\beta + L) \, p^i_h + \dvrg \vectorz^i_h\|^2_{\Omega}$.
In order to achieve the desired efficiency in the error estimate, the reliability term must be several orders of
magnitude less than the dual one, which is satisfied in this case.

\begin{table}[!t]
\footnotesize
\scriptsize
\centering
\newcolumntype{g}{>{\columncolor{gainsboro}}c}
\newcolumntype{k}{>{\columncolor{lightgray}}c}
\newcolumntype{s}{>{\columncolor{silver}}c}
\newcolumntype{a}{>{\columncolor{ashgrey}}c}
\newcolumntype{b}{>{\columncolor{battleshipgrey}}c}
\begin{tabular}{c|c|cg|cc|k}
\parbox[c]{2.4cm}{\centering optimization step} &
\parbox[c]{0.8cm}{\centering $\zeta$} &
\parbox[c]{1.4cm}{\centering $|\!|\!| e_p |\!|\!|^2$ }  &
\parbox[c]{1.4cm}{\centering $\overline{\rm M}^h_{p}$ } &
\parbox[c]{1.4cm}{\centering $\overline{\rm m}^2_{\rm d}$ } &
\parbox[c]{1.4cm}{\centering $\overline{\rm m}^2_{\rm eq}$ } &
\parbox[c]{2.4cm}{\centering $\Ieff := \overline{\rm M}^h_{p} / |\!|\!| e_p |\!|\!|^2$ } \\
\midrule
0 & 0.4359 & 1.87e-04 & 4.55e-04 & 1.87e-04 & 8.06e-05 & 1.5589 \\
1 & 0.0045 & 1.87e-04 & 2.69e-04 & 1.87e-04 & 1.53e-08 & 1.1983 \\
2 & 0.0002 & 1.87e-04 & 1.88e-04 & 1.87e-04 & 1.67e-08 & 1.0023 \\
\end{tabular}
\caption{{\it Example \ref{ex:1}}. Decrease of the values of the parameter $\zeta$ and the majorant $\overline{\rm M}^h_{p}$,
as well as the terms $\overline{\rm m}^2_{\rm d}$ and $\overline{\rm m}^2_{\rm eq}$ it contains,
w.r.t. the number of optimization cycles.}
\label{tab:example-1-majorant-optimization}
\end{table}

\begin{figure}[!t]
\centering
\subfloat[Error $|\!|\!| e_p |\!|\!|^2$ and indicator  $\overline{\rm m}^2_{\rm d, \tensorK}$]{
\distr{
\includegraphics[width=8cm, trim={0.2cm 0.7cm 1.5cm 1.2cm}, clip]{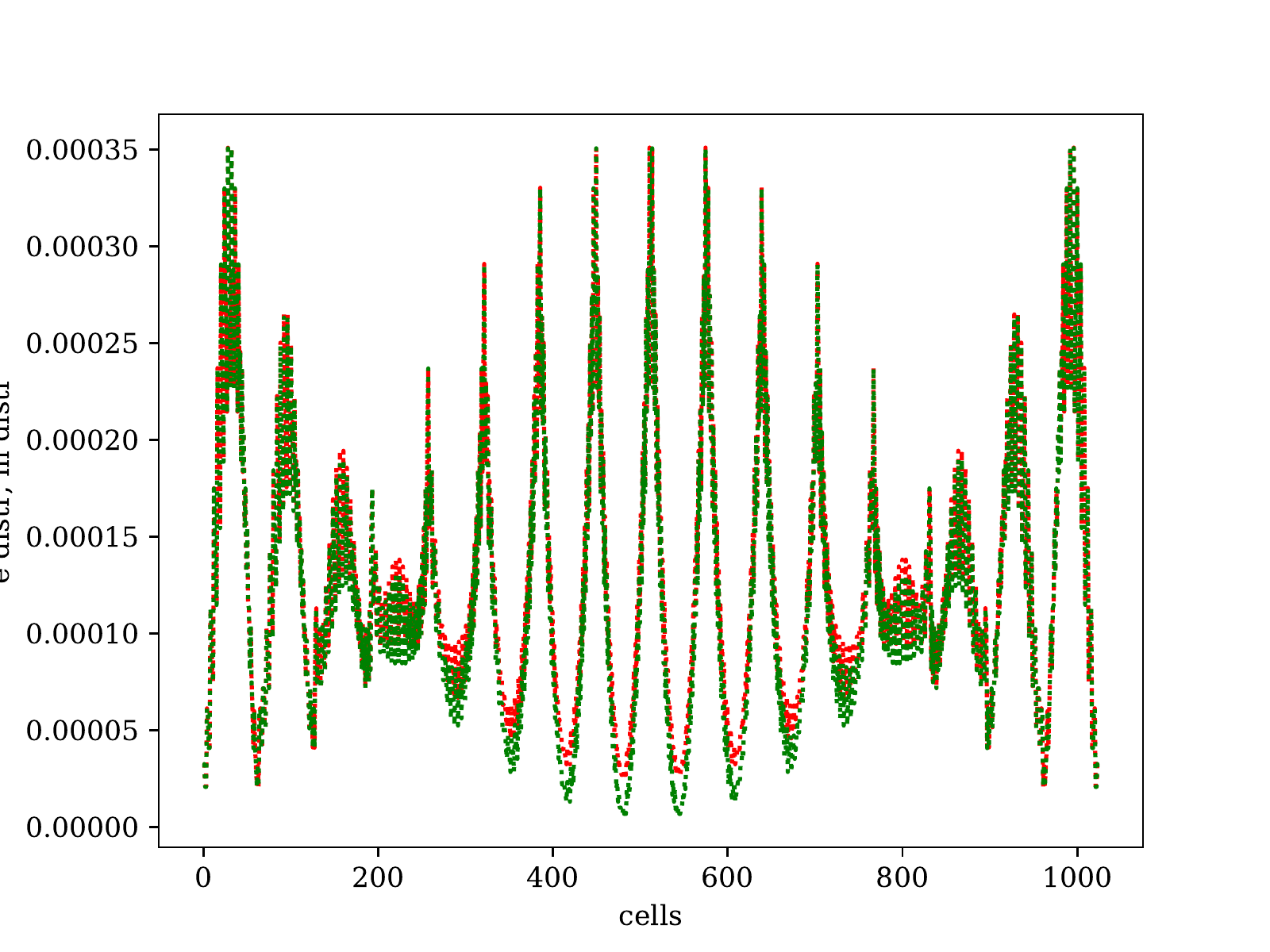}
}
}
\subfloat[Error $|\!|\!| e_{\vectoru} |\!|\!|^2$ and indicator  $\overline{\rm m}^2_{\rm d, \mu, \lambda}$]{
\distr{
\includegraphics[width=8cm, trim={0.2cm 0.7cm 1.5cm 1.2cm}, clip]{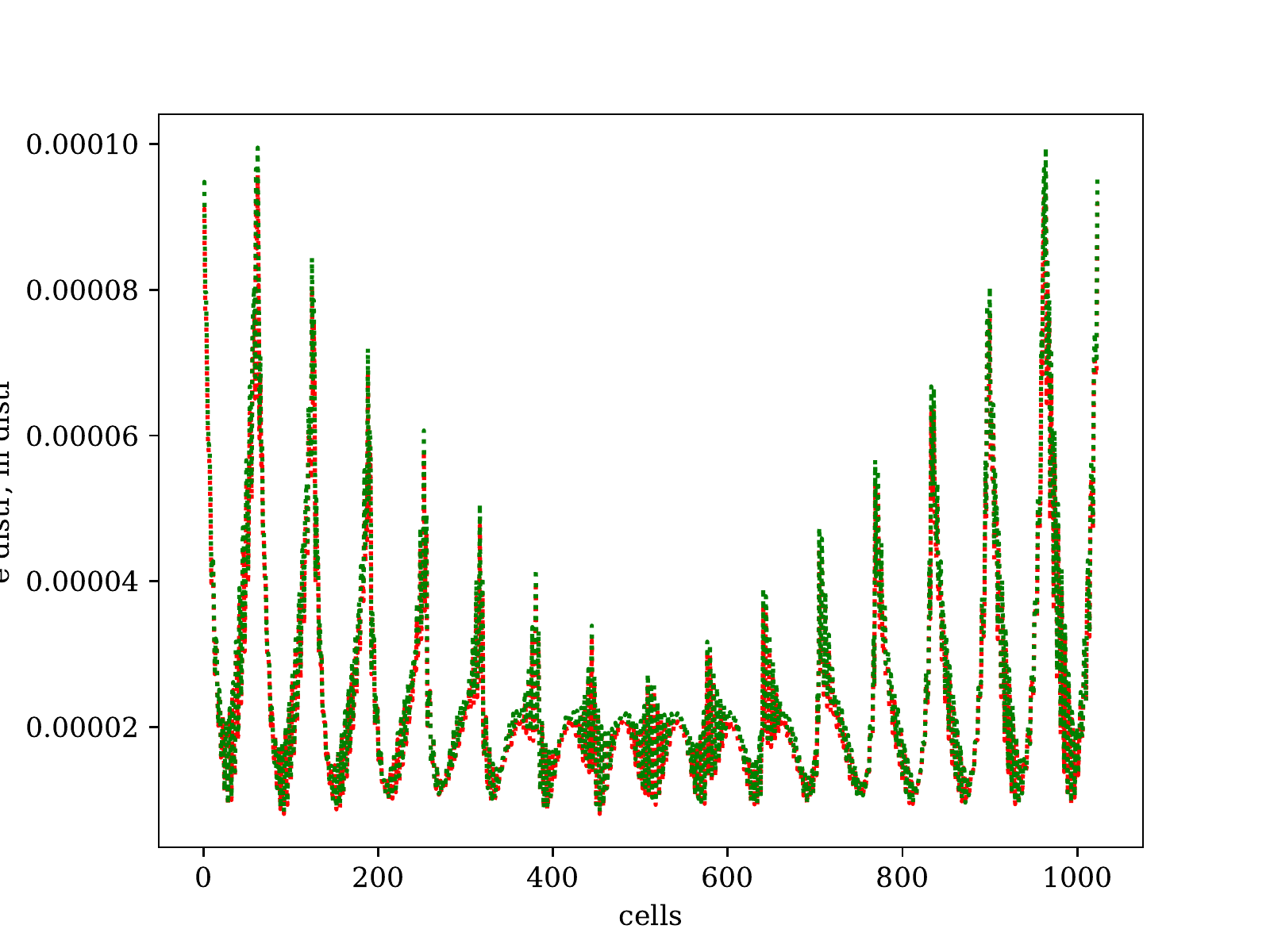}
}
}\\
\subfloat[Error $|\!|\!| e_p |\!|\!|^2$ and indicator  $\overline{\rm m}^2_{\rm d, \tensorK}$]{\distr{
\includegraphics[width=8cm, trim={0.3cm 0.7cm 1.5cm 1.2cm}, clip]{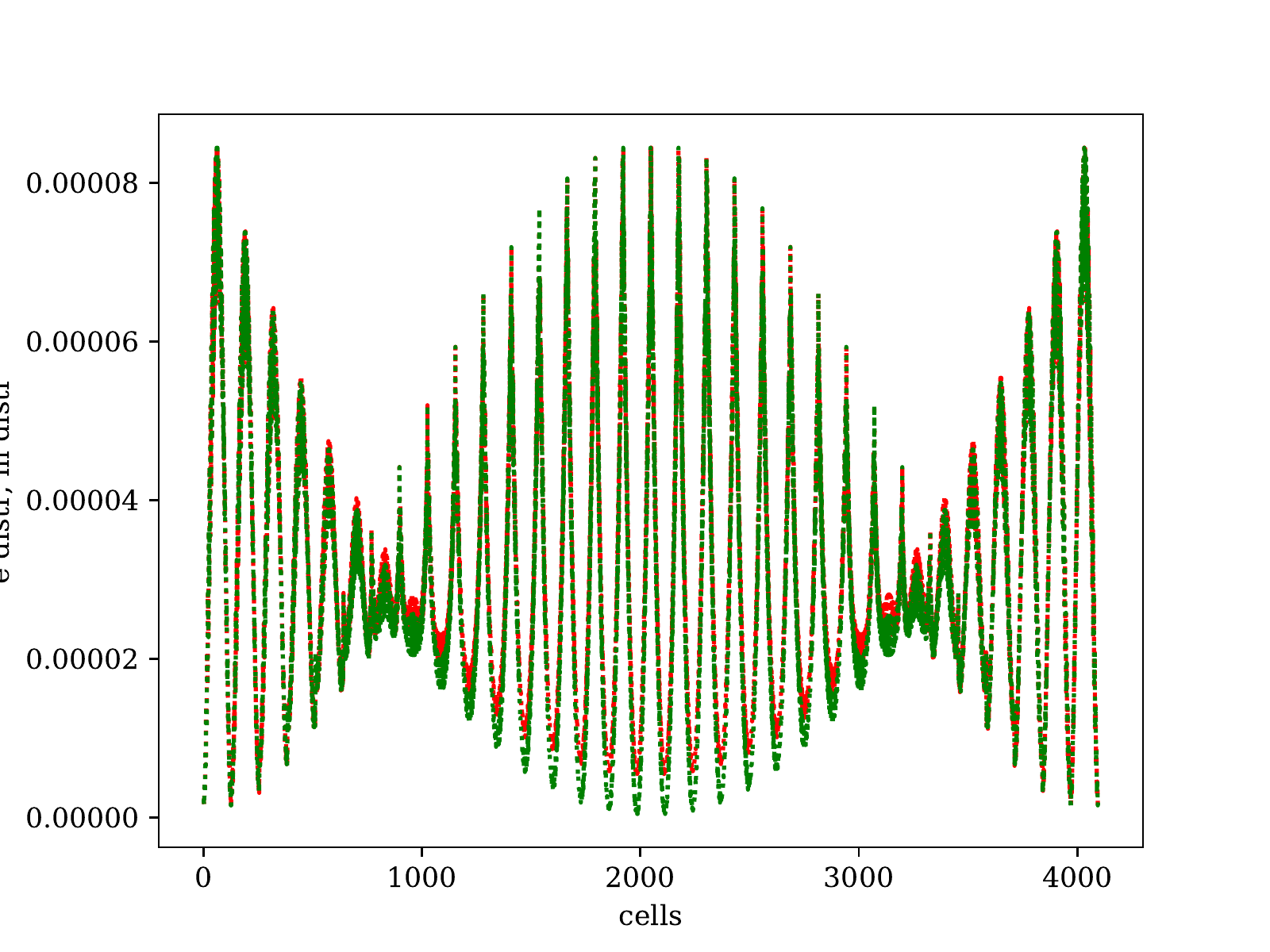}
}
}
\subfloat[Error $|\!|\!| e_{\vectoru} |\!|\!|^2$ and indicator  $\overline{\rm m}^2_{\rm d, \mu, \lambda}$]{
\distr{
\includegraphics[width=8cm, trim={0.1cm 0.7cm 1.5cm 1.2cm}, clip]{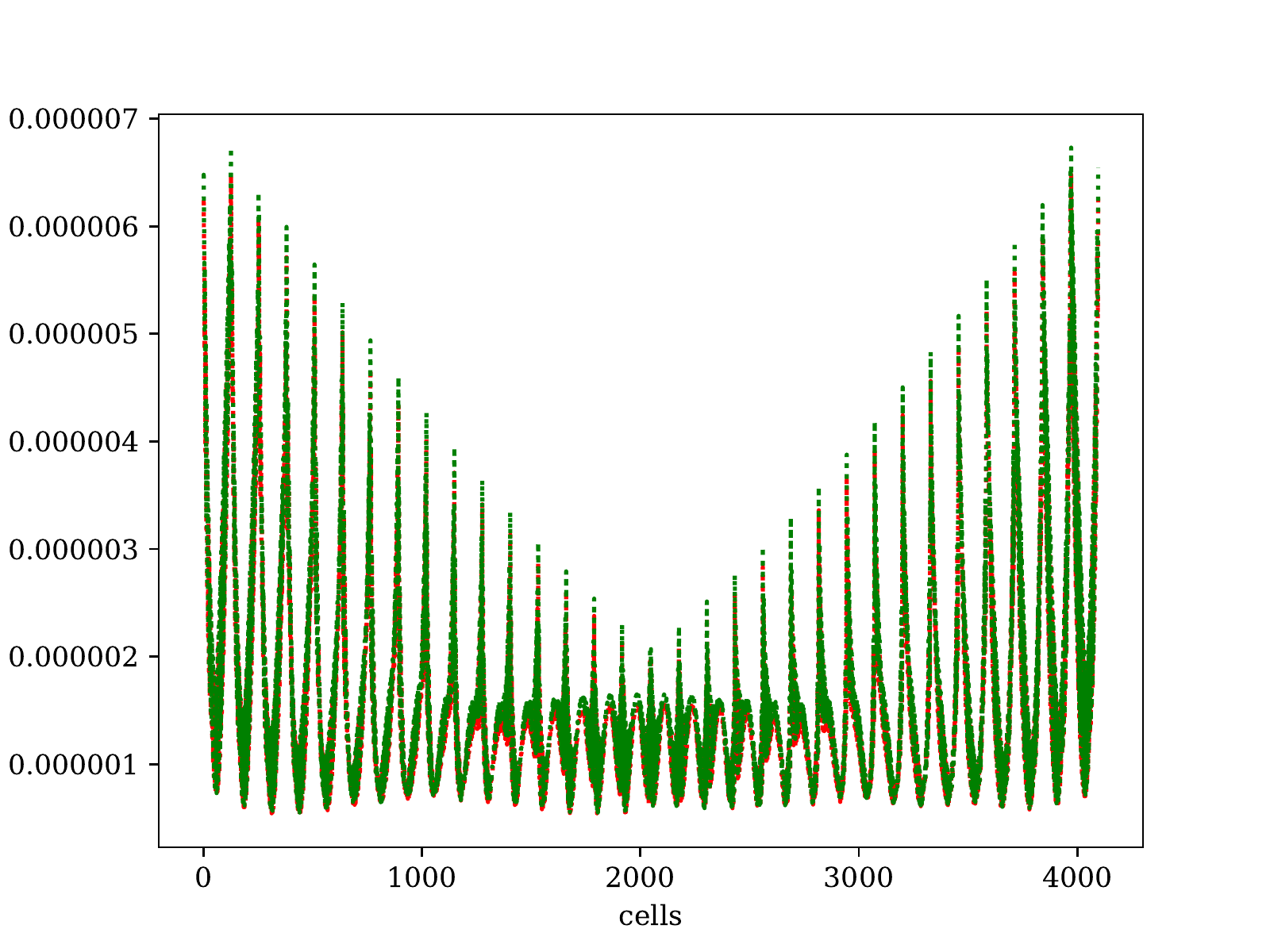}
}
}
\caption{{\it Example \ref{ex:1}}. (a),(c) Local distribution of error $|\!|\!| e_p |\!|\!|^2$ in red and error indicator
$\overline{\rm m}^2_{\rm d, \tensorK}$ generated by the majorant $\overline{\rm M}^h_{p}$ in green and
(b),(d) local distribution of error $|\!|\!| e_{\vectoru} |\!|\!|^2$ in red and error indicator
$\overline{\rm m}^2_{\rm d, \mu, \lambda}$ generated by the majorant $\overline{\rm M}^h_{\vectoru}$ in green.
The distributions are depicted w.r.t. the numbered finite element cells of the uniform meshes with
(a)--(b)  $h = \rfrac{1}{32}$ and (c)--(d) $h = \rfrac{1}{64}$.}
\label{fig:example-1-adaptive-mesh}
\end{figure}
Analogous minimization is performed for the majorant
$\overline{\rm M}^h_{u}$. After minimization, the terms $\overline{\rm m}^2_{\rm d, \tensorK}$ and
$\overline{\rm m}^2_{\rm d, \lambda, \mu}$ become very close to the true errors $|\!|\!| e_p |\!|\!|^2$ and
$|\!|\!| e_u |\!|\!|^2$, respectively, and can be used as error indicators for the local error distribution over
the computational domain. To confirm that, Figure \ref{fig:example-1-adaptive-mesh} presents the distribution of the errors
and indicators (generated by majorants) w.r.t. to the numbered finite element cells. On the right side of
Figure \ref{fig:example-1-adaptive-mesh}, we depict $|\!|\!| e_p |\!|\!|^2$ and
$\overline{\rm m}^2_{\rm d, \tensorK} := \| \mathbf{r}_{{\rm} d, \tensorK}({p}^{i}_h, \vectorz^i_h)\|^2$
and on the left side the error in $\vectoru$ measured in terms of $|\!|\!| e_u |\!|\!|^2$ as well as the indicator
$\overline{\rm m}^2_{\rm d, \mu, \lambda} := \| \mathbf{r}_{{\rm} d, \mu, \lambda}({\vectoru}^{i}_h, \vectortau^i_h)\|^2$
w.r.t. to the numbered finite element cells. To depict the error green marker is used, whereas for the majorant the red
color is chosen. We see that the indicator produces quantitatively efficient error indicator.
First, let us assume that the error bounds $\overline{\rm M}^h_{p}$ and $\overline{\rm M}^h_{p, L^2}$, as well as 
$\overline{\rm M}^h_{u}$ and $\overline{\rm M}^h_{u, \dvrg}$ are reconstructed only on the 
$4$th and $5$th steps to calculate {$\overline{\rm M}^{\,i}_{\vectoru}$} and {$\overline{\rm M}^{\,i}_{p}$} 
for $i=5$ (using \eqref{eq:estimate-p-iterative} and \eqref{eq:estimate-u-iterative}).
Let $|\!|\!| e^{(n)}_p |\!|\!|^2$ and $|\!|\!| e^{(n)}_u |\!|\!|^2$, $n = 1, \ldots, N$ correspond to the error
increments corresponding to the $n$th time-step. The contribution of the majorants of discretization errors and
the iterative majorant on the $N$th time step are of the same magnitude, i.e.,
$${\overline{\rm M}}^{h, (N)}_{p} = \mbox{1.8736e-04},  \quad
\overline{\rm M}^{\,i, (N)}_{p} = \mbox{3.0550e-04},  \quad
\overline{\rm M}^{h, (N)}_{u} = \mbox{4.9549e-04},  \quad
\overline{\rm M}^{\,i, (N)}_{\vectoru} = \mbox{3.2025e-05}.$$
Then, $N$th increment of both total errors and corresponding majorants
for the pressure as well as the displacement respectively are as follows:
\begin{alignat*}{2}
|\!|\!| e^{(N)}_p |\!|\!|^2 & = \mbox{4.8652e-05}, \quad
\overline{\rm M}^{(N)}_p = \mbox{2.5603e-04}, \quad
\Ieff(\overline{\rm M}^{(N)}_p) = \mbox{2.29} \quad \mbox{and}\quad \\
|\!|\!| e^{(N)}_u |\!|\!|^2 & = \mbox{4.8536e-05}, \quad
\overline{\rm M}^{(N)}_u = \mbox{2.7403e-04}, \quad
\Ieff(\overline{\rm M}^{(N)}_u) = 2.38,
\end{alignat*}
making the general error, majorant, and corresponding efficiency index contributions the following  
$$[(e^{(N)}_u, e^{(N)}_p)]^2 = \mbox{4.8559e-04}, \quad
{\overline{\rm M}^{(N)}} = \mbox{9.9944e-05}, \quad
\Ieff({\overline{\rm M}^{(N)}}) = 2.36.$$

If instead of $4$th and $5$th iterations we use $2$nd and $5$th to compute {$\overline{\rm M}^{\,i, m, (N)}_{\vectoru}$} and
{$\overline{\rm M}^{\,i, m, (N)}_{p}$} (using $m = 3$ in \eqref{eq:estimates-of-nabla-p-via-sigma-2} and
\eqref{eq:estimates-of-strain-via-sigma-2}). Then contribution of iterative majorants is
$$\overline{\rm M}^{\,i, m, (N)}_{p} = \mbox{5.0348e-06} \quad \mbox{and}
\quad \overline{\rm M}^{\,i, m, (N)}_{\vectoru} = \mbox{ 2.8107e-08}.$$
Then, the accumulated overall time steps values of the relative errors and majorants are

\begin{alignat*}{4}
|\!|\!| e_p |\!|\!|^2 = \sum_{n=1}^N |\!|\!| e^{(n)}_p |\!|\!|^2 & = \mbox{1.8731e-04}, \quad &
|\!|\!| e_u |\!|\!|^2 = \sum_{n=1}^N |\!|\!| e^{(n)}_u |\!|\!|^2 & = \mbox{1.8687e-04}, \quad \mbox{and} \quad \nonumber\\
\overline{\rm M}_p = \sum_{n=1}^N \overline{\rm M}^{(n)}_p & = \mbox{3.8479e-04}, \quad &
\overline{\rm M}_u = \sum_{n=1}^N \overline{\rm M}^{(n)}_u & = \mbox{9.9104e-04},
\end{alignat*}
respectively, with corresponding efficiency indices over the whole time interval are summarised in Table
\ref{tab:example-1-convergence-wrt-mesh-size-and-tau}(a). The results are presented corresponding to meshes with
two different mesh-sizes $h$ and $\tau$. As the caption highlights, all these values are obtained for the auxiliary
functions are reconstructed by the following finite elements $\boldsymbol{z}^i_h \in {\rm \mathds{RT}_1}$ and
$\vectortau^i_h \in [{\rm P_2}]^{2 \times 2}$, which is an order higher then usual choice of the finite element
approximation spaces for the fluxes and stresses in mixed formulations. However, it is important to obtain better
accuracy of the error bounds. To emphasize this importance, Table \ref{eq:example-1-different-z-and-tau} lists the
results obtained for cases, when different finite elements pairs for auxiliary functions are used, i.e., 
\authors{
(a) $\boldsymbol{z}^i_h \in {\rm \mathds{RT}_0}$ and $\vectortau^i_h \in [{\rm P_2}]^{2 \times 2}$ or 
(b) $\boldsymbol{z}^i_h \in {\rm \mathds{RT}_1}$ and $\vectortau^i_h \in [{\rm \mathds{P}_1}]^{2 \times 2}$.}

\begin{table}[!t]
\footnotesize
\centering
\newcolumntype{g}{>{\columncolor{gainsboro}}c} 	
\newcolumntype{k}{>{\columncolor{lightgray}}c} 	
\newcolumntype{s}{>{\columncolor{silver}}c} 
\newcolumntype{a}{>{\columncolor{ashgrey}}c}
\newcolumntype{b}{>{\columncolor{battleshipgrey}}c}
\begin{tabular}{c|rs|rs|rs|k}
\parbox[c]{0.8cm}{\centering $h$} & 
\parbox[c]{1.1cm}{\centering  $|\!|\!| e_p |\!|\!|^2 $}   & 	  
\parbox[c]{1.1cm}{\centering $\overline{\rm M}_{p}$  }&    
\parbox[c]{1.1cm}{\centering $|\!|\!| e_u |\!|\!|^2$ } & 
\parbox[c]{1.1cm}{\centering $\overline{\rm M}_{\vectoru}$  }&
\parbox[c]{1.1cm}{\centering $[(e_u, e_p)]^2$  }&  
\parbox[c]{1.1cm}{\centering ${\overline{\rm M}}$ }  & 
\parbox[c]{1.1cm}{\centering $\Ieff ({\overline{\rm M}})$ } \\
\midrule
\multicolumn{8}{l}{(a) $N = 10$, \; \quad $\tau = 1.0$} \\
\midrule
$\rfrac{1}{4}$  &	   4.8163e-02 &	   1.0299e-01 &	   4.8035e-02 &	   2.6402e-01 &	   6.1827e-02 &	  2.9352e-01 &	 2.18 \\
$\rfrac{1}{8}$  &	   1.2006e-02 &	   2.4971e-02 &	   1.1977e-02 &	   6.4119e-02 &	   1.5415e-02 &	  7.1270e-02 &	 2.15 \\
$\rfrac{1}{16}$  &	   2.9981e-03 &	   6.1776e-03 &	   2.9909e-03 &	   1.5899e-02 &	   3.8495e-03 &	  1.7668e-02 &	 2.14 \\
$\rfrac{1}{32}$  &	   7.4930e-04 &	   1.5398e-03 &	   7.4752e-04 &	   3.9663e-03 &	   9.6209e-04 &	  4.4072e-03 &	 2.14 \\
$\rfrac{1}{64}$  &	   1.8731e-04 &	   3.8479e-04 &	   1.8687e-04 &	   9.9104e-04 &	   2.4050e-04 &	  1.1012e-03 &	 2.14 \\
\midrule
\multicolumn{8}{l}{(b) $N = 10^2$, \quad $\tau = 0.1$} \\
\midrule
$\rfrac{1}{4}$  &	   3.4606e-02 &	   6.9224e-02 &	   4.8035e-03 &	   2.6245e-02 &	   6.2192e-03 &	  2.9077e-02 &	 2.16 \\
$\rfrac{1}{8}$  &	   8.4602e-03 &	   1.7171e-02 &	   1.1977e-03 &	   6.3753e-03 &	   1.5438e-03 &	  7.0778e-03 &	 2.14 \\
$\rfrac{1}{16}$ &	   2.1022e-03 &	   4.2867e-03 &	   2.9909e-04 &	   1.5809e-03 &	   3.8509e-04 &	  1.7563e-03 &	 2.14 \\
$\rfrac{1}{32}$  &	   5.2472e-04 &	   1.0754e-03 &	   7.4752e-05 &	   3.9441e-04 &	   9.6218e-05 &	  4.3841e-04 &	 2.13 \\
$\rfrac{1}{64}$  &	   1.3113e-04 &	   2.7329e-04 &	   1.8687e-05 &	   9.8551e-05 &	   2.4051e-05 &	  1.0973e-04 &	 2.14 \\
\midrule
\multicolumn{8}{l}{(c) $N = 10^3$, \quad $\tau = 0.01$} \\
\midrule
$\rfrac{1}{4}$ &	   1.0880e-02 &	   3.5205e-02 &	   4.8035e-04 &	   2.6104e-03 &	   6.5839e-04 &	  3.1865e-03 &	 2.20 \\
$\rfrac{1}{8}$ &	   2.2543e-03 &	   8.7416e-03 &	   1.1977e-04 &	   6.3403e-04 &	   1.5666e-04 &	  7.7708e-04 &	 2.23 \\
$\rfrac{1}{16}$ &	   5.3424e-04 &	   2.1800e-03 &	   2.9909e-05 &	   1.5722e-04 &	   3.8652e-05 &	  1.9289e-04 &	 2.23 \\
$\rfrac{1}{32}$ &	   1.3172e-04 &	   5.4464e-04 &	   7.4752e-06 &	   3.9223e-05 &	   9.6307e-06 &	  4.8135e-05 &	 2.24 \\
$\rfrac{1}{64}$ &	   3.2816e-05 &	   1.3621e-04 &	   1.8687e-06 &	   9.8005e-06 &	   2.4056e-06 &	  1.2029e-05 &	 2.24 \\
\end{tabular}
\caption{{\it Example \ref{ex:1}}. 
Convergence of errors and majorants w.r.t. the different choice of spatial mesh sizes $h$ and time steps $\tau$
(measured relatively w.r.t. $|\!|\!| p |\!|\!|^2_p$ and $|\!|\!| \vectoru |\!|\!|^2_{u}$). For all the cases,
the auxialiary functions are reconstrcted by the following finite elements
$\boldsymbol{z}^i_h \in {\rm \mathds{RT}_1}$ and $\vectortau^i_h \in [{\rm P_2}]^{2 \times 2}$.}
\label{tab:example-1-convergence-wrt-mesh-size-and-tau}
\end{table} 

\begin{table}[!t]
\footnotesize
\centering
\newcolumntype{g}{>{\columncolor{gainsboro}}c}
\newcolumntype{k}{>{\columncolor{lightgray}}c}
\newcolumntype{s}{>{\columncolor{silver}}c}
\newcolumntype{a}{>{\columncolor{ashgrey}}c}
\newcolumntype{b}{>{\columncolor{battleshipgrey}}c}
\begin{tabular}{c|rs|rs|rs|k}
\parbox[c]{0.8cm}{\centering $h$} &
\parbox[c]{1.1cm}{\centering  $|\!|\!| e_p |\!|\!|^2 $}   &
\parbox[c]{1.1cm}{\centering $\overline{\rm M}_{p}$  }&
\parbox[c]{1.1cm}{\centering $|\!|\!| e_u |\!|\!|^2$ } &
\parbox[c]{1.1cm}{\centering $\overline{\rm M}_{\vectoru}$  }&
\parbox[c]{1.1cm}{\centering $[(e_u, e_p)]^2$  }&
\parbox[c]{1.1cm}{\centering ${\overline{\rm M}}$ }  &
\parbox[c]{1.1cm}{\centering $\Ieff ({\overline{\rm M}})$ } \\
\bottomrule
\multicolumn{8}{l}{ \rule{0pt}{4ex}
(a) $\boldsymbol{z}^i_h \in {\rm \mathds{RT}_0}$, $\vectortau^i_h \in [{\rm P_2}]^{2 \times 2}$} \\[3pt]
\toprule
$\rfrac{1}{4}$ &	   1.5184e-02 &	   9.9777e-02 &	   1.2934e-02 &	   7.3430e-02 &	   2.8119e-02 &	  1.7321e-01 &	 2.48 \\
$\rfrac{1}{8}$ &	   3.7854e-03 &	   2.5695e-02 &	   3.2249e-03 &	   1.7884e-02 &	   7.0102e-03 &	  4.3579e-02 &	 2.49 \\
$\rfrac{1}{16}$ &	   9.4525e-04 &	   6.4694e-03 &	   8.0535e-04 &	   4.4381e-03 &	   1.7506e-03 &	  1.0907e-02 &	 2.50 \\
$\rfrac{1}{32}$ &	   2.3624e-04 &	   1.6202e-03 &	   2.0128e-04 &	   1.1074e-03 &	   4.3752e-04 &	  2.7276e-03 &	 2.50 \\
$\rfrac{1}{64}$ &	   5.9055e-05 &	   4.0529e-04 &	   5.0316e-05 &	   2.7672e-04 &	   1.0937e-04 &	  6.8201e-04 &	 2.50 \\
\bottomrule
\multicolumn{8}{l}{ \rule{0pt}{4ex}
(b) $\boldsymbol{z}^i_h \in {\rm \mathds{RT}_1}$, $\vectortau^i_h \in [{\rm \mathds{P}_1}]^{2 \times 2}$} \\[3pt]
\toprule
$\rfrac{1}{4}$ &	   1.5184e-02 &	   2.3151e-01 &	   1.2934e-02 &	   2.9586e-01 &	   2.8119e-02 &	  5.2737e-01 &	 4.33 \\
$\rfrac{1}{8}$ &	   3.7854e-03 &	   5.9349e-02 &	   3.2249e-03 &	   7.6502e-02 &	   7.0102e-03 &	  1.3585e-01 &	 4.40 \\
$\rfrac{1}{16}$ &	   9.4525e-04 &	   1.4935e-02 &	   8.0535e-04 &	   1.9301e-02 &	   1.7506e-03 &	  3.4236e-02 &	 4.42 \\
$\rfrac{1}{32}$ &	   2.3624e-04 &	   3.7402e-03 &	   2.0128e-04 &	   4.8369e-03 &	   4.3752e-04 &	  8.5771e-03 &	 4.43 \\
$\rfrac{1}{64}$ &	   5.9055e-05 &	   9.3545e-04 &	   5.0316e-05 &	   1.2100e-03 &	   1.0937e-04 &	  2.1454e-03 &	 4.43 \\
\end{tabular}
\caption{{\it Example \ref{ex:1}}. 
Convergence of errors and majorants w.r.t. the different choice of approximation spaces for $\boldsymbol{z}^i_h$ and
$\vectortau^i_h$ (measured relatively w.r.t. $|\!|\!| p |\!|\!|^2_p$ and $|\!|\!| \vectoru |\!|\!|^2_{u}$).
For all the cases, simulation is performed for the discretization with $\tau = 1.0$.
}
\label{eq:example-1-different-z-and-tau}
\end{table}

Assume now that $N=10^3$. Analogously, we consider the last time-step $[t_{999}, t_{1000}]$ with mesh-size
$h = \tfrac{1}{64}$. Table \ref{tab:example-1-iteration-convergence-relative-N-1000} illustrates
the convergence of the errors in $\vectoru$ and $p$ w.r.t. iteration steps (in this case, we consider $I = 10$).
For considered $\tau=0.01$, the contribution of the majorants of the discretization errors and the iterative majorant
on the pressure variable differ in magnitude, i.e.,
$${\overline{\rm M}}^{h,(N)}_{p} =  \mbox{6.6685e-05}, \quad
{\overline{\rm M}}^{i,(N)}_{p} = \mbox{3.0410e-02}, \quad
{\overline{\rm M}}^{h,(N)}_{u} = \mbox{4.9003e-04}, \quad
{\overline{\rm M}}^{i,(N)}_{u} = \mbox{3.1551e-05}.$$
\begin{table}[!ht]
\footnotesize
\centering
\newcolumntype{g}{>{\columncolor{gainsboro}}c}
\newcolumntype{k}{>{\columncolor{lightgray}}c}
\newcolumntype{s}{>{\columncolor{silver}}c}
\newcolumntype{a}{>{\columncolor{ashgrey}}c}
\newcolumntype{b}{>{\columncolor{battleshipgrey}}c}
\begin{tabular}{c|ck|ck|ck|ck}
\parbox[c]{1.6cm}{\centering $i = 1, \ldots, I$ } &
\parbox[c]{1.4cm}{\centering $|\!|\!| e_p |\!|\!|^2$ }  &
\parbox[c]{1.4cm}{\centering $\overline{\rm M}^h_{p}$ } &
\parbox[c]{1.0cm}{\centering $\| e_p \|^2_{\beta}$ } &
\parbox[c]{1.4cm}{\centering $\overline{\rm M}^h_{p, L^2}$ } &
\parbox[c]{1.4cm}{\centering $|\!|\!| e_u |\!|\!|^2$ } &
\parbox[c]{1.4cm}{\centering $\overline{\rm M}^h_{u}$  }&
\parbox[c]{1.0cm}{\centering $\| \dvrg e_u \|^2_{\lambda}$ }   &
\parbox[c]{1.4cm}{\centering $\overline{\rm M}^h_{u, \dvrg}$ } \\
\midrule
1 & 5.8753e+00 &    $-$ &     3.5255e+00  &   $-$ & 1.7140e-02 &    $-$  &     5.2771e-03   &  $-$ \\
2 & 3.0103e-03 &    $-$ &     1.3101e-03  &   $-$ & 1.9054e-04 &    $-$  &     4.5713e-05   &  $-$ \\
3 & 5.5125e-05 &    8.6791e-05 &     1.6877e-05  &   5.7961e-05 & 1.8688e-04 &    4.9901e-04  &     4.4550e-05   &  1.0693e-04 \\
4 & 3.3497e-05 &    $-$ &     5.9589e-07  &   $-$ & 1.8687e-04 &    $-$  &     4.4544e-05   &  $-$ \\
$...$ & $...$ &    $...$ &     $...$ &   $...$ & $...$ &    $...$  &     $...$   &  $...$ \\
8 & 3.2817e-05 &    $-$ &     3.8002e-08  &   $-$ & 1.8687e-04 &    $-$  &     4.4544e-05   &  $-$ \\
9 & 3.2817e-05 &    6.6685e-05 &     3.7904e-08  &   4.4534e-05 & 1.8687e-04 &    4.9003e-04  &     4.4544e-05   &  1.0501e-04 \\
10 & 3.2817e-05 &    6.6685e-05 &     3.7889e-08  &   4.4534e-05 & 1.8687e-04 &    4.9003e-04  &     4.4544e-05   &  1.0501e-04 \\
%
\end{tabular}
\caption{{\it Example \ref{ex:1}}. 
Errors errors and majorants w.r.t. the iteration steps for $I = 7$, $h = \tfrac{1}{64}$, $N = 10^3$, \; \quad $[t_{999}, t_{1000}]$, $\tau = 0.0 1$
(both values are measured relative to the increment in $|\!|\!| p |\!|\!|^2_p$ and $|\!|\!| \vectoru |\!|\!|^2_{u}$ on the $N$th time step).}
\label{tab:example-1-iteration-convergence-relative-N-1000}
\end{table}
These are the results obtained considering two last subsequent iteration steps with $q = 0.2307$ that
correspond to $\tfrac{3 \, q}{1-q^2} 
\Bigg(\tfrac{\big(C^{\rm F}_{\Sigma^p_D}\big)^2\, \beta}{\lambda_\tensorK \, \tau} + 1 \Bigg) = 7.0973$.
Such a difference in ${\overline{\rm M}}^{h, (N)}_{p}$ and ${\overline{\rm M}}^{i, (N)}_{p}$ results in the error
majorant with the efficiency index $\Ieff(\overline{\rm M}^{(N)}_p) = \mbox{43.10}$.
However, if we assume, say, the third and tenth steps to be subsequent ones, it provides a better contractive
parameter $\tilde{q} := q^7$ = 3.4853e-05 and improves the efficiency index of error estimates corresponding
to pressure and displacement by approximately 21.12 and 1.03 times, respectively, i.e.,
$${\overline{\rm M}}^{i,m,(N)}_{p} = \mbox{1.4154e-06} \quad \mbox{and} \quad
{\overline{\rm M}}^{i,m,(N)}_{u} = \mbox{2.2179e-13}$$
This yields
\begin{alignat*}{2}
|\!|\!| e^{(N)}_p |\!|\!|^2 & =  \mbox{9.8302e-08},
\quad \overline{\rm M}^{(N)}_p = \mbox{4.0799e-07}, \quad
\Ieff(\overline{\rm M}^{(N)}_p) = \mbox{2.04} \quad \mbox{and}\quad \\
|\!|\!| e^{(N)}_u |\!|\!|^2 & =  \mbox{5.5976e-09},
\quad \overline{\rm M}^{(N)}_u = \mbox{2.9358e-08},
\quad \Ieff(\overline{\rm M}^{(N)}_u) = 2.29.
\end{alignat*}
The accumulated values of the pressure error, majorant, and corresponding efficiency index over the whole 
time interval are as follows:
$$|\!|\!| e_p |\!|\!|^2 = \mbox{3.2816e-05}, \quad
\overline{\rm M}_p = \mbox{1.3621e-04}, \quad \mbox{and} \quad
\Ieff(\overline{\rm M}_p) = 2.04.$$
For the displacement, we observe  
$$|\!|\!| e_u |\!|\!|^2 = \mbox{1.8687e-06}, \quad
\overline{\rm M}_u = \mbox{9.8005e-06}, \quad \mbox{and} \quad
\Ieff(\overline{\rm M}_u) = 2.29.$$
This yields the total values $\big|\!\big[ ({e}_{\vectoru}, {e}_p) \big]\!\big| = \mbox{2.4056e-06}$ and
$\overline{\rm M}$ = 1.2029e-05 with efficiency index 
$\Ieff := \tfrac{\overline{\rm M}}{\big|\!\big[ ({e}_{\vectoru}, {e}_p) \big]\!\big| } = 2.24$.
The latter values are included in Table \ref{tab:example-1-convergence-wrt-mesh-size-and-tau}(c). They 
illustrate the accumulated values of the errors and majorant over the whole time interval.
We see that even with the decreasing $\tau$, which scales the permeability tensor $\tensorK$, 
the efficiency of total error estimates stays rather robust.

We note here, that each increment of $\overline{\rm M}^{(N)}_p $ and $\overline{\rm M}^{(N)}_u$ can be computed using
$\reallywidetilde{\rm M}^{\,i}_{p}$ and $\reallywidetilde{\rm M}^{\,i}_{u}$, which does not require computation of the
majorant of discretization errors ${\overline{\rm M}}^{h,(N)}_{p}$ and ${\overline{\rm M}}^{h,(N)}_{u}$ on iteration
steps $i-1$ or $i-m$ and, as consequence,
their minimization w.r.t. the auxiliary functions $\vectorz^i_h$ and $\vectortau^i_h$, respectively. Moreover,
the values $\reallywidetilde{\rm M}^{\,i}_{p}$ and $\reallywidetilde{\rm M}^{\,i}_{u}$ are orders of magnitude smaller
than ${\overline{\rm M}}^{i,m,(N)}_{p}$ and ${\overline{\rm M}}^{i,m,(N)}_{u}$, i.e.,
$\reallywidetilde{\rm M}^{\,i}_{p}$ = {\rm 7.0735e-10} {\rm and}
$\reallywidetilde{\rm M}^{\,i}_{u}$ = {\rm 7.3390e-13}. Such values automatically minimize the contribution of iterative estimates 
into ${\overline{\rm M}}^{h,(N)}_{p}$ and ${\overline{\rm M}}^{h,(N)}_{u}$ and even further improves their efficiency. 

\vskip 10pt
\noindent
{\em Realistic parameters.}
Let us assume now more realistic parameters in the above example (similar to those taken from \cite{Bothetall2017} and \cite{MikelicWangWheeler2014}). 
Let exact pressure be scaled as follows: $p(x, y, t) = 10^8 \, x \, (1 - x) \, y \, (1 - y) \, t$. The permeability tensor divided by fluid viscosity 
is taken as $\tensorK = 100 \tensorI$ $[\rm mD/cP]$, fluid compressibility fixed to $4.7 \cdot 10^{-7}$ $[\rm psi^{-1}]$, and 
initial porosity $\phi_0$ is assumed  to be 0.2.
The Biot and bulk modulus are $M = 1.65 \cdot 10^{10}$ [Pa] and $E = 0.594 \cdot 10^{9}$ [Pa], 
respectively. That results into $\beta = \tfrac{1}{M} + c_f \, \phi_0 = 9.406 \cdot 10^{-8}$, which, in turn, yields considerable small 
$L =\tfrac{\alpha^2}{2 \, (\lambda + 2 \mu /d)} = 1.35 \cdot 10^{-9}$. Such a tuning parameter generates 
instantaneous convergence of iterative scheme with contractive parameter $q = \tfrac{L}{\beta + L} = 6.73 \cdot 10^{-12}$.
The resulting errors and corresponding estimates are summarised in Table \ref{tab:example-1-2-confergence-wrt-mesh-size-and-tau}
(for $\tau = 0.1$ and difference spacial mesh-sizes $h$). 
\sveta{For such $q$, even one iteration is enough for convergence. 
However, one can consider two to five iterations to improve the efficiency of the majorant. 
The efficiency indices obtained here confirm the efficiency of obtained majorants also for parameters closed to 
those used in engineering applications.}
\begin{table}[!t]
\footnotesize
\centering
\newcolumntype{g}{>{\columncolor{gainsboro}}c} 	
\newcolumntype{k}{>{\columncolor{lightgray}}c} 	
\newcolumntype{s}{>{\columncolor{silver}}c} 
\newcolumntype{a}{>{\columncolor{ashgrey}}c}
\newcolumntype{b}{>{\columncolor{battleshipgrey}}c}
\begin{tabular}{c|rs|rs|rs|k}
\parbox[c]{0.8cm}{\centering $h$} & 
\parbox[c]{1.1cm}{\centering  $|\!|\!| e_p |\!|\!|^2 $}   & 	  
\parbox[c]{1.1cm}{\centering $\overline{\rm M}_{p}$  }&    
\parbox[c]{1.1cm}{\centering $|\!|\!| e_u |\!|\!|^2$ } & 
\parbox[c]{1.1cm}{\centering $\overline{\rm M}_{\vectoru}$  }&
\parbox[c]{1.1cm}{\centering $[(e_u, e_p)]^2$  }&  
\parbox[c]{1.1cm}{\centering ${\overline{\rm M}}$ }  & 
\parbox[c]{1.1cm}{\centering $\Ieff ({\overline{\rm M}})$ } \\
%
\midrule
$\rfrac{1}{4}$  & 5.0067e-02 & 3.9082e-01 & 4.7948e-03 & 8.3004e-02 & 5.0067e-02 & 3.9082e-01 & 2.79 \\
$\rfrac{1}{8}$  & 1.2576e-02 & 9.6727e-02 & 1.2024e-03 & 2.0992e-02 & 1.2576e-02 & 9.6727e-02 & 2.77 \\
$\rfrac{1}{16}$  & 3.1461e-03 & 2.3977e-02 & 3.0070e-04 & 5.2587e-03 & 3.1461e-03 & 2.3977e-02 & 2.76 \\
$\rfrac{1}{32}$  & 7.8664e-04 & 5.9649e-03 & 7.5179e-05 & 1.3152e-03 & 7.8664e-04 & 5.9649e-03 & 2.75 \\
$\rfrac{1}{64}$  & 1.9667e-04 & 1.4874e-03 & 1.8795e-04 & 3.2882e-03 & 1.9667e-04 & 1.4874e-03 & 2.75 \\
\end{tabular}
\caption{{\it Example \ref{ex:1}}. 
Convergence of errors and majorants w.r.t. the different choice of spatial mesh sizes 
(measured relatively w.r.t. $|\!|\!| p^i_h |\!|\!|^2_p$ and $|\!|\!| \vectoru^i_h |\!|\!|^2_{u}$) for
$N = 100$ and $\tau = 0.1$.}
\label{tab:example-1-2-confergence-wrt-mesh-size-and-tau}
\end{table} 

\end{example}

\begin{example} 
\label{ex:2}
\rm
{\bf Simplified parameters}.
Let $\Omega := (0, 1)^2 \in \Rtwo$, $T = 10$. 
The exact solution of \eqref{eq:poroelastic-system} is defined as
$$\vectoru(x, y, t) := 
\begin{bmatrix} 
t \, (x^2 + y^2) \\ 
t \, (x + y) 
\end{bmatrix} 
\quad  \mbox{and} \quad
p(x, y, t) := t \, x \, (1 - x)\, y \, (1 - y).$$ 
%
We fix Poisson ratio to be $\nu = 0.2$ and Young modulus $E = 0.594$, which yields the Lame parameters $\mu = 0.2475$ and 
$\lambda = 0.1238$. We set $\alpha = 1$ and $\beta = \tfrac{1}{M} + c_f \phi_0 = 0.1165$, where $M = 1.65 \cdot 10^{10}$, 
$K_u = K + \tfrac{\alpha^2}{c_0} = 10.2887$, $c_f = \tfrac{1}{c_0} \lambda_\tensorK \tfrac{K + 4/3 \mu}{K_u + 4/3 \mu} = 0.5827$
$C_{\rm F} = \tfrac{1}{\sqrt{2}\, \pi}$, and $\tensorK = I$, where $\tensorI$ is a unit tensor. From 
parameters above, it follows that $L =\tfrac{\alpha^2}{2 \, (\lambda + 2 \mu /d)} =  1.3468$ and 
$q = \tfrac{L}{\beta + L} = 0.9203$. With such $q$ the ratio $\tfrac{q^2}{(1 - q)^2}$ is 133.3343, which might influence 
the quantitative quality of the majorant. 
%
\begin{table}[!t]
\footnotesize
\centering
\newcolumntype{g}{>{\columncolor{gainsboro}}c}
\newcolumntype{k}{>{\columncolor{lightgray}}c}
\newcolumntype{s}{>{\columncolor{silver}}c}
\newcolumntype{a}{>{\columncolor{ashgrey}}c}
\newcolumntype{b}{>{\columncolor{battleshipgrey}}c}
\begin{tabular}{c|ck|ck|ck|ck}
\parbox[c]{0.8cm}{\centering \# it. } & 
\parbox[c]{1.4cm}{\centering $|\!|\!| e_p |\!|\!|^2$ }  & 	  
\parbox[c]{1.4cm}{\centering $\overline{\rm M}^h_{p}$ } &    
\parbox[c]{1.0cm}{\centering $\| e_p \|^2_{\beta}$ } & 
\parbox[c]{1.4cm}{\centering $\overline{\rm M}^h_{p, L^2}$ } & 
\parbox[c]{1.4cm}{\centering $|\!|\!| e_u |\!|\!|^2$ } & 
\parbox[c]{1.4cm}{\centering $\overline{\rm M}^h_{u}$  }&    
\parbox[c]{1.0cm}{\centering $\| \dvrg e_u \|^2_{\lambda}$ }   & 	
\parbox[c]{1.4cm}{\centering $\overline{\rm M}^h_{u, \dvrg}$ } \\
\midrule
\multicolumn{9}{l}{(a) $N = 10$, \; \quad $\tau = 1.0$} \\
\midrule
     1 & 5.3962e+00 &    $-$ &     2.9479e-02  &   $-$ & 8.4773e-04 &    $-$  &     1.2996e-03   &  $-$ \\
    $...$ & $...$ &    $...$ &     $...$ &   $...$ & $...$ &    $...$  &     $...$   &  $...$ \\
     5 & 1.9553e-04 &    $-$ &     2.4696e-10  &   $-$ & 8.5024e-06 &    $-$  &     9.8152e-06   &  $-$ \\
     6 & 1.9553e-04 &    1.9559e-04 &     2.6065e-10  &   9.2250e-06 & 8.5024e-06 &    2.6115e-05  &     9.8152e-06   &  2.3448e-05 \\
$...$ & $...$ &    $...$ &     $...$ &   $...$ & $...$ &    $...$  &     $...$   &  $...$ \\
    10 & 1.9553e-04 &    $-$ &     2.6127e-10  &   $-$ & 8.5024e-06 &    $-$  &     9.8152e-06   &  $-$ \\
    11 & 1.9553e-04 &    1.9559e-04 &     2.6127e-10  &   9.2250e-06 & 8.5024e-06 &    2.6115e-05  &     9.8152e-06   &  2.3448e-05 \\
    12 & 1.9553e-04 &    1.9559e-04 &     2.6127e-10  &   9.2250e-06 & 8.5024e-06 &    2.6115e-05  &     9.8152e-06   &  2.3448e-05 \\
\end{tabular}
\caption{{\it Example \ref{ex:2}}. Errors errors and majorants w.r.t. the iteration steps for $h = \tfrac{1}{64}$ and $I = 12$
(both values are measured relative to the increment in $|\!|\!| p |\!|\!|^2_p$ and $|\!|\!| \vectoru |\!|\!|^2_{u}$ on the $N$th time step).}
\label{tab:example-2-iteration-convergence-last}
\end{table}

We consider $10$ time-steps of the length $\tau$ = 1.0 discretizing considered interval and spatial mesh-size $h = \tfrac{1}{64}$ 
using standard $\mathds{P}_1$ finite elements (see \eqref{eq:p0-p1}). The number of iterations to solve the problem on each time-step is 
set to 12. For the time-interval $[t_9, t_{10}] = [9.0, 10.0]$, the convergence of the errors in $\vectoru$ and $p$ is 
presented by Table \ref{tab:example-2-iteration-convergence-last}. From one side, we can consider iterations $I$ and $I-1$
as subsequent ones with a contraction parameter $q = 0.9203$. However, by using Lemma \ref{eq:lemma-2} instead, let $m = 6$ so 
that $6$th and $12$th iterations be treated as two consecutive steps but with $q' = q^6 = 0.6027$. Then, the constants dependent from 
the ratio $\tfrac{q^2}{(1 - q^2)} = 2.3012$ attain more acceptable values. 
Table \ref{tab:example-2-confergence-wrt-mesh-size-and-tau} illustrates the improvement of error majorants efficiency indices as 
the explained-above method is employed. 
The local distribution of the error in $p$ and $\vectoru$ on each cell of the finite-element discretization is presented for 
mesh sizes $h = \rfrac{1}{4}$ and $h= \rfrac{1}{8}$ in Figure \ref{fig:example-2-adaptive-mesh}. One can see the 
resemblance in the local error distribution for $p$ since the exact solutions for both examples are the same. The local values of 
error and indicator (generated by the majorant) in $\vectoru$ has more uniform distribution.
%
\begin{table}[!t]
\footnotesize
\centering
\newcolumntype{g}{>{\columncolor{gainsboro}}c} 	
\newcolumntype{k}{>{\columncolor{lightgray}}c} 	
\newcolumntype{s}{>{\columncolor{silver}}c} 
\newcolumntype{a}{>{\columncolor{ashgrey}}c}
\newcolumntype{b}{>{\columncolor{battleshipgrey}}c}
\begin{tabular}{c|rs|rs|rs|k}
\parbox[c]{0.8cm}{\centering $h$} & 
\parbox[c]{1.4cm}{\centering  $|\!|\!| e_p |\!|\!|^2 $}   & 	  
\parbox[c]{1.4cm}{\centering $|\!|\!| e_u |\!|\!|^2$ } &
\parbox[c]{1.4cm}{\centering $\big|\!\big[ ({e}_{\vectoru}, {e}_p) \big]\!\big|$ } &
\parbox[c]{1.4cm}{\centering $\Ieff(\overline{\rm M}_{p})$} &    
\parbox[c]{1.4cm}{\centering $\Ieff(\overline{\rm M}_{u})$} &    
\parbox[c]{1.4cm}{\centering $\Ieff(\overline{\rm M})$} \\
\midrule
\multicolumn{7}{l}{(a) $N = 10$, \; \quad $\tau = 1.0$} \\
\midrule
$\rfrac{1}{4}$  & 5.0147e-02 & 2.1767e-03 &	2.5823e-03 & 40.48 $\rightarrow$ 4.94 & 20.61 $\rightarrow$ 3.15 & 24.81 $\rightarrow$ 3.49 \\
$\rfrac{1}{8}$  & 1.2525e-02 & 5.4416e-04 &	6.4546e-04 & 40.46 $\rightarrow$ 4.93 & 20.59 $\rightarrow$ 3.15 & 24.78 $\rightarrow$ 3.49 \\
$\rfrac{1}{16}$ & 3.1293e-03 & 1.3604e-04 & 1.6135e-04 & 40.46 $\rightarrow$ 4.93 & 20.58 $\rightarrow$ 3.15 & 24.78 $\rightarrow$ 3.49 \\
$\rfrac{1}{32}$ & 7.8217e-04 & 3.4010e-05  & 4.0335e-05 & 40.46 $\rightarrow$ 4.93 & 20.58 $\rightarrow$ 3.15 & 24.77 $\rightarrow$  3.49 \\
$\rfrac{1}{64}$ & 1.9553e-04 & 8.5024e-06 & 1.0084e-05 & 40.46 $\rightarrow$ 4.93 & 20.58 $\rightarrow$ 3.15 & 24.77 $\rightarrow$  3.49 \\
\midrule
\multicolumn{7}{l}{(b) $N = 10^2$, \quad $\tau = 0.1$} \\
\midrule
$\rfrac{1}{4}$  & 4.7811e-02 & 4.7811e-02 &	2.5836e-04 & 129.34 $\rightarrow$ 15.18 & 20.34 $\rightarrow$ 3.12 & 54.62 $\rightarrow$ 6.67 \\
$\rfrac{1}{8}$  & 1.1914e-02 & 5.4416e-05 &	6.4554e-05 & 129.53 $\rightarrow$ 15.27 & 20.33 $\rightarrow$ 3.12 & 54.62 $\rightarrow$ 6.70 \\
$\rfrac{1}{16}$  & 2.9748e-03 & 1.3604e-05 &	1.6135e-05 & 129.60 $\rightarrow$ 15.27 & 20.33 $\rightarrow$ 3.14 & 54.62 $\rightarrow$ 6.81 \\
$\rfrac{1}{32}$  & 7.4345e-04 & 3.4010e-06 &	4.0336e-06  & 129.62 $\rightarrow$ 16.66 & 20.33 $\rightarrow$ 3.23 & 54.62 $\rightarrow$ 7.21 \\
$\rfrac{1}{64}$  & 1.8585e-04 & 8.5024e-07 &	1.0084e-06  & 129.62 $\rightarrow$ 20.45 & 20.33 $\rightarrow$ 3.56 & 54.62 $\rightarrow$ 8.74 \\
\end{tabular}
\caption{{\it Example \ref{ex:2}}. 
Convergence of errors and majorants w.r.t. the different choice of spatial mesh sizes and time steps 
(measured relatively w.r.t. $|\!|\!| p^i_h |\!|\!|^2_p$ and $|\!|\!| \vectoru^i_h |\!|\!|^2_{u}$).}
\label{tab:example-2-confergence-wrt-mesh-size-and-tau}
\end{table}

\begin{figure}[!t]
	\centering
	\subfloat[$|\!|\!| e_p |\!|\!|^2$ and $\overline{\rm m}^2_{\rm d, \tensorK}(p)$]{
	\distr{
	\includegraphics[width=8cm, trim={0.5cm 0.7cm 1.5cm 1.2cm}, clip]{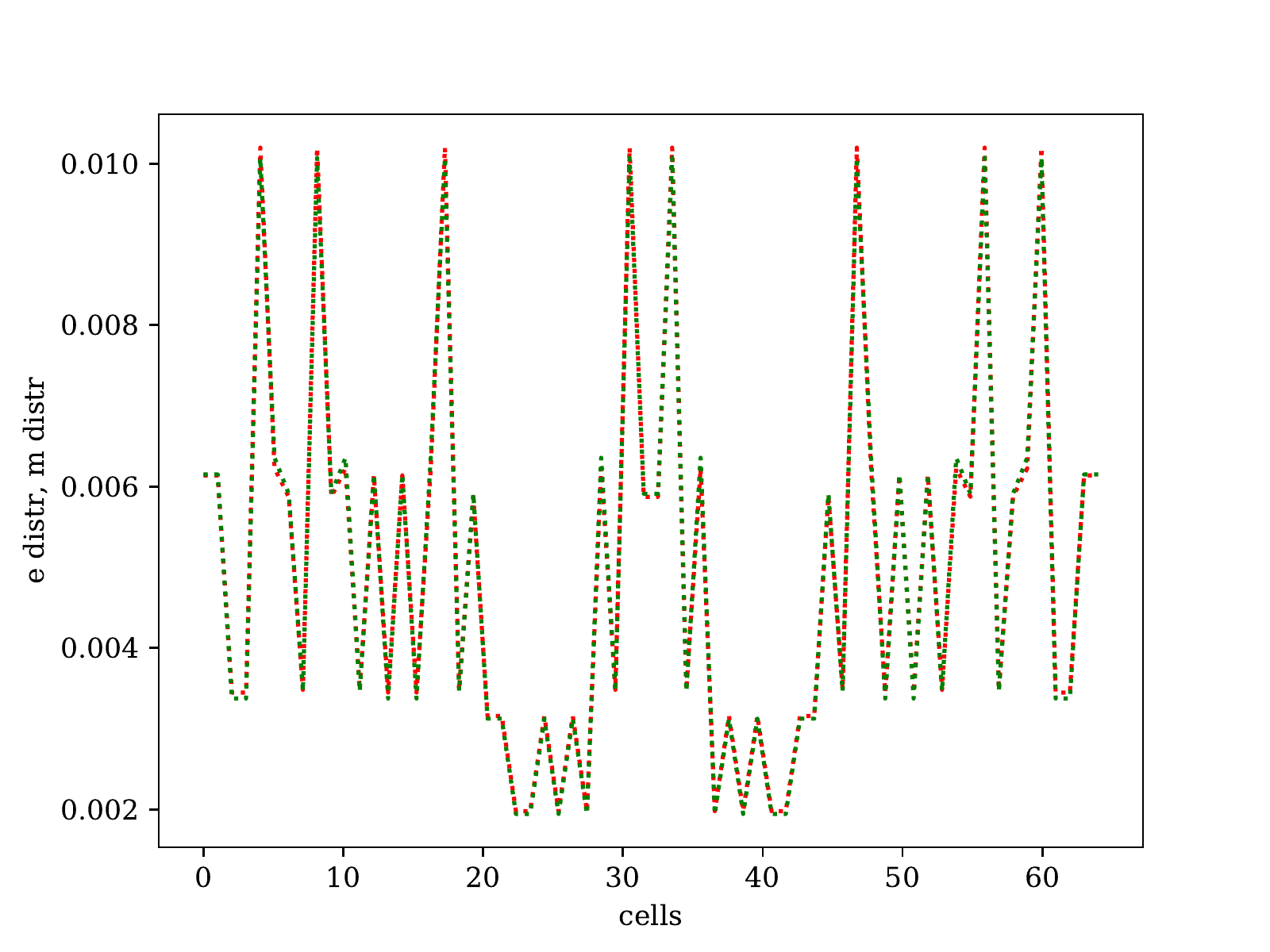}}
	} 
	\subfloat[$|\!|\!| e_{\vectoru} |\!|\!|^2$ and $\overline{\rm m}^2_{\rm d, \mu, \lambda}$]{
	\distr{
	\includegraphics[width=8cm, trim={0.5cm 0.7cm 1.5cm 1.2cm}, clip]{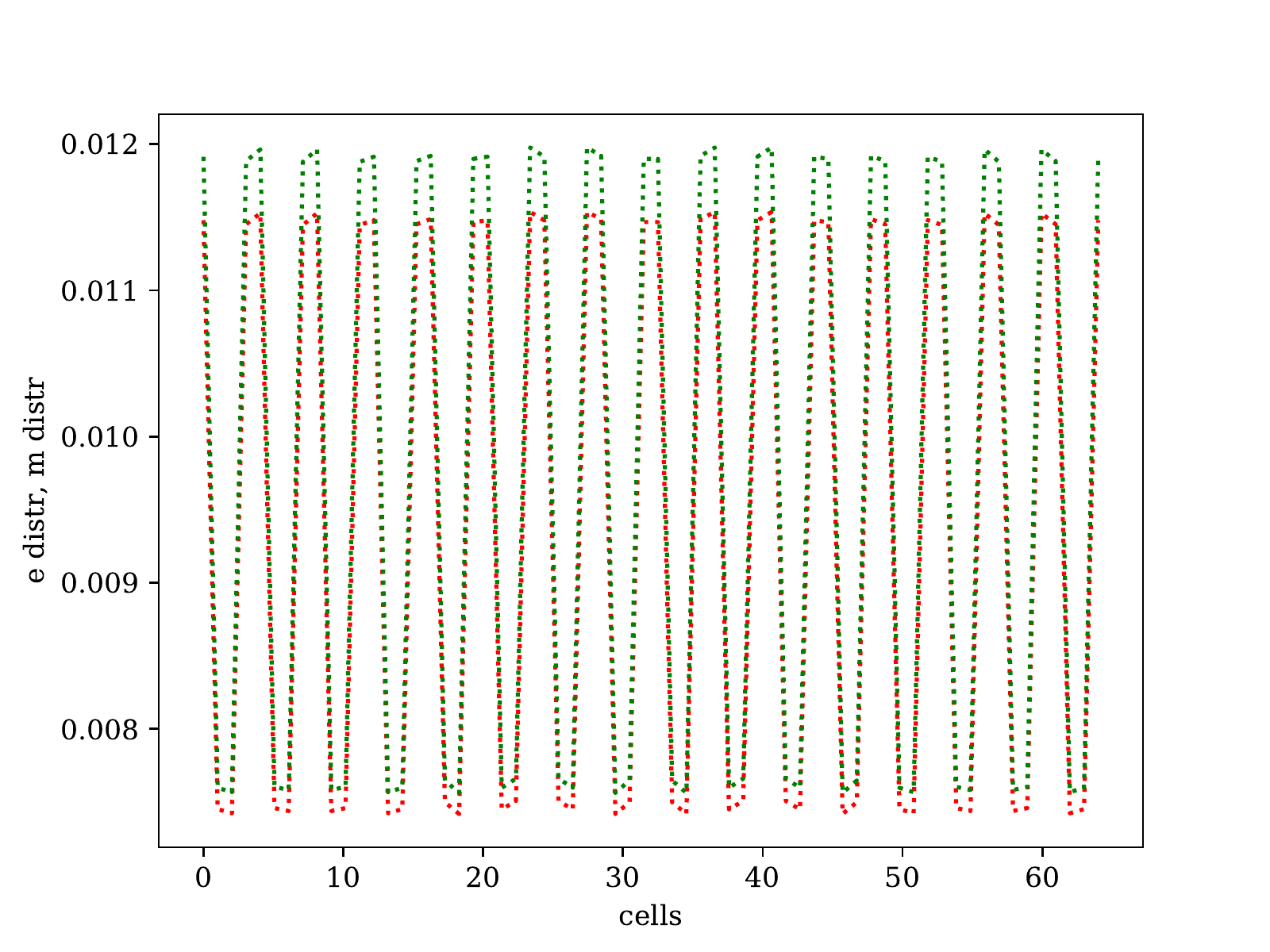}}
	}
	\\
	\subfloat[Error $|\!|\!| e_p |\!|\!|^2$ and indicator $\overline{\rm m}^2_{\rm d, \tensorK}(p)$]{
	\distr{
	\includegraphics[width=8cm, trim={0.5cm 0.7cm 1.5cm 1.2cm}, clip]{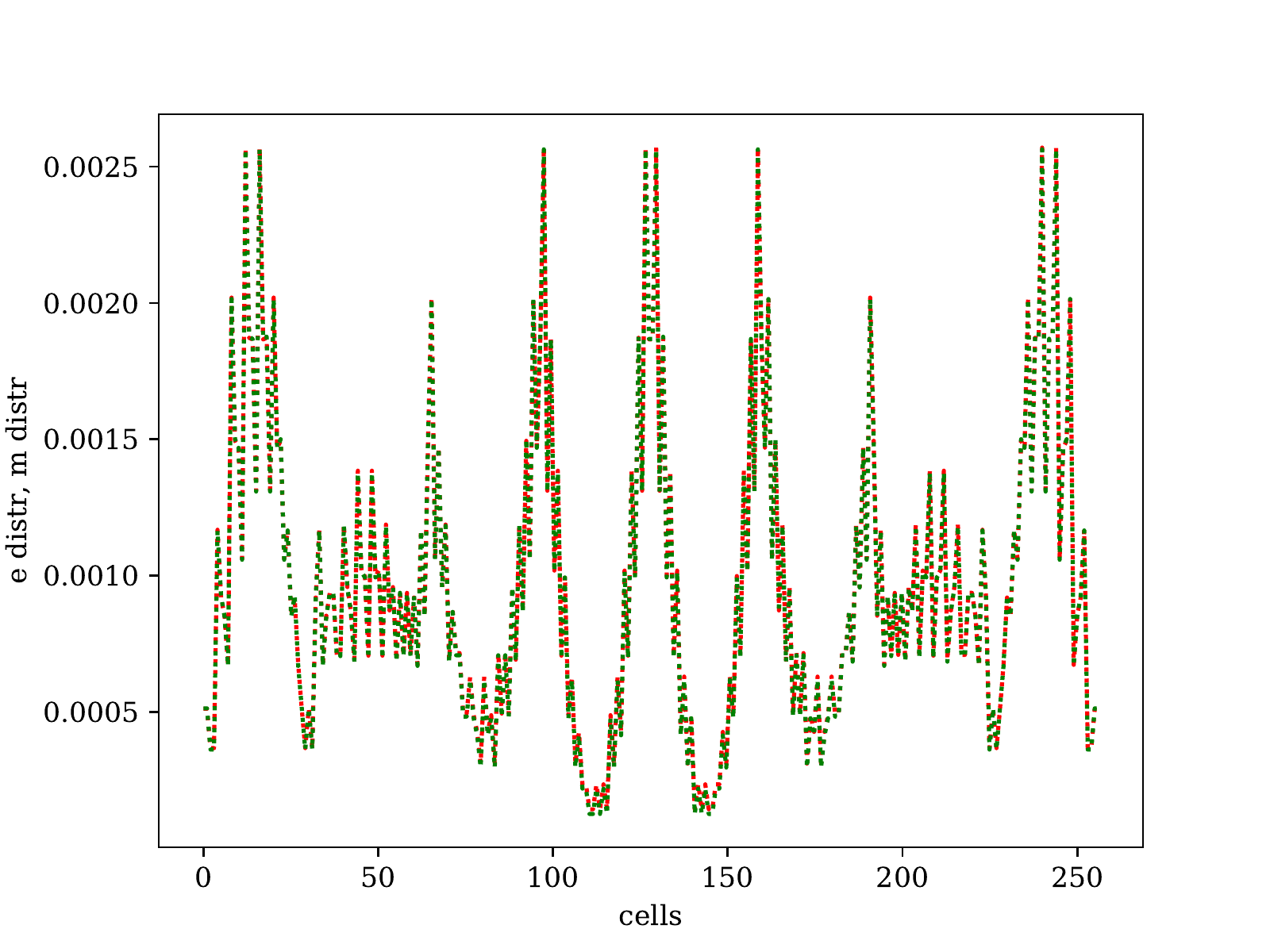}} 
	}
	\subfloat[Error $|\!|\!| e_{\vectoru} |\!|\!|^2$ and indicator  $\overline{\rm m}^2_{\rm d, \mu, \lambda}$]{
	\distr{
	\includegraphics[width=8cm, trim={0.5cm 0.7cm 1.5cm 1.2cm}, clip]{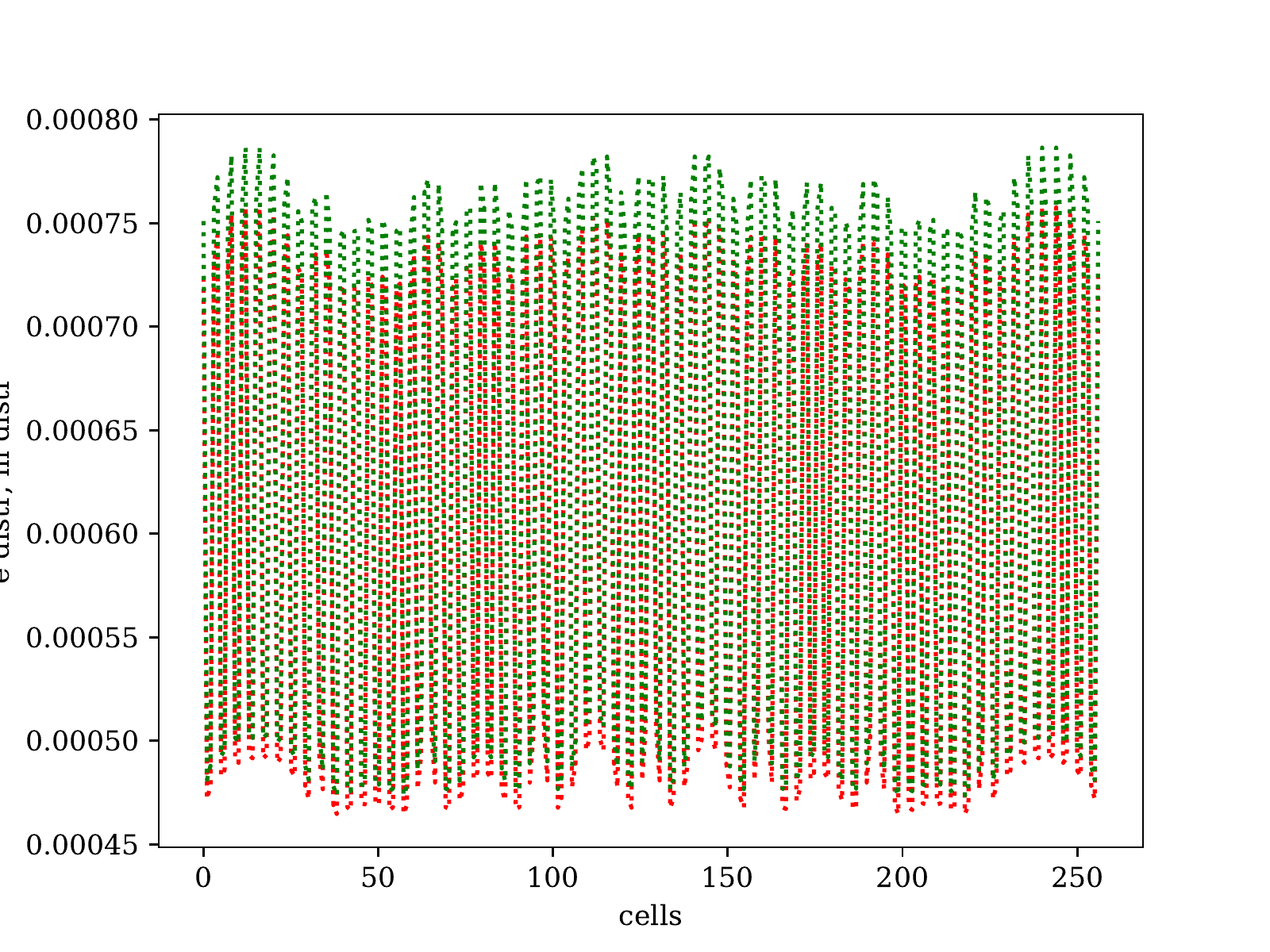}
	}
	}
	\caption{{\it Example \ref{ex:2}}. Errors $|\!|\!| e_p |\!|\!|^2$ and $|\!|\!| e_{\vectoru} |\!|\!|^2$ (in red)  and 
	error indicators $\overline{\rm m}^2_{\rm d, \tensorK}$ and $\overline{\rm m}^2_{\rm d, \mu, \lambda}$ (in green) 
	generated by the majorants $\overline{\rm M}^h_{p}$ and $\overline{\rm M}^h_{\vectoru}$, respectively, 
	distributions w.r.t. to the numbered finite element cells. }
	\label{fig:example-2-adaptive-mesh}
\end{figure}

{\bf Realistic parameters}.
Let us consider more realistic parameters similar to those considered in \cite{Giraultetal2015}.
Again, let exact pressure be $p(x, y, t) = 10^8 \, x \, (1 - x) \, y \, (1 - y) \, t$.
Mechanical parameters such as Biot and bulk modulus are chosen as follows:
$M = 1.4514 \cdot 10^{4}$ $[\rm Pa]$ and $E = 7 \cdot 10^{7}$ $[\rm Pa]$, 
respectively. 
The permeability tensor is$\tensorK =  {\rm diag} \{50, 200\} \cdot 10^{-10}$ $[\rm mD/cP]$, whereas 
the fluid viscosity is $\mu_f = 10^{-3}$.
Thus, we obtain $\beta = \tfrac{1}{M} + c_f \phi_0 = 6.8909 \cdot 10^{-5}$, which, in turn, yields considerable small 
$L =\tfrac{\alpha^2}{2 \, (\lambda + 2 \mu /d)} =1.1428 \cdot 10^{-8}$. 
\sveta{Such parameter generates contractive parameter $q = \tfrac{L}{\beta + L} = 1.6584 \cdot 10^{-4}$.}
The latter means that the ratio $\tfrac{q^2}{(1 - q)^2}$ is of order $10^{-8}$, which does not degenerate the efficiency 
of the majorant.
We summarise the convergence results in Table \ref{tab:example-2-2-confergence-wrt-mesh-size-and-tau}. It confirms that even 
for more realistic parameters (common for engineering applications), the value of the estimates remains rather efficient.
\begin{table}[!t]
\footnotesize
\centering
\newcolumntype{g}{>{\columncolor{gainsboro}}c} 	
\newcolumntype{k}{>{\columncolor{lightgray}}c} 	
\newcolumntype{s}{>{\columncolor{silver}}c} 
\newcolumntype{a}{>{\columncolor{ashgrey}}c}
\newcolumntype{b}{>{\columncolor{battleshipgrey}}c}
\begin{tabular}{c|rs|rs|rs|k}
\parbox[c]{0.8cm}{\centering $h$} & 
\parbox[c]{1.1cm}{\centering  $|\!|\!| e_p |\!|\!|^2 $}   & 	  
\parbox[c]{1.1cm}{\centering $\overline{\rm M}_{p}$  }&    
\parbox[c]{1.1cm}{\centering $|\!|\!| e_u |\!|\!|^2$ } & 
\parbox[c]{1.1cm}{\centering $\overline{\rm M}_{\vectoru}$  }&
\parbox[c]{1.1cm}{\centering $[(e_u, e_p)]^2$  }&  
\parbox[c]{1.1cm}{\centering ${\overline{\rm M}}$ }  & 
\parbox[c]{1.1cm}{\centering $\Ieff ({\overline{\rm M}})$ } \\
\midrule
\multicolumn{8}{l}{
\qquad (a) $N = 10$, \; \quad $\tau = 1.0$} \\
\midrule
$\rfrac{1}{4}$ &	   3.9463e-02 &	   2.0213e-01 &	   2.1767e-03 &	   1.3970e-02 &	   3.9663e-02 &	  2.0342e-01 &	 2.26 \\
$\rfrac{1}{8}$ &	   9.6879e-03 &	   4.9714e-02 &	   5.4416e-04 &	   3.4872e-03 &	   9.7380e-03 &	  5.0034e-02 &	 2.27 \\
$\rfrac{1}{16}$ &	   2.4097e-03 &	   1.2333e-02 &	   1.3604e-04 &	   8.7104e-04 &	   2.4222e-03 &	  1.2413e-02 &	 2.26 \\
$\rfrac{1}{32}$ &	   6.0164e-04 &	   3.0727e-03 &	   3.4010e-05 &	   2.1767e-04 &	   6.0477e-04 &	  3.0927e-03 &	 2.26 \\
$\rfrac{1}{64}$ &	   1.5036e-04 &	   7.6693e-04 &	   8.5024e-06 &	   5.4407e-05 &	   1.5114e-04 &	  7.7193e-04 &	 2.26 \\
\midrule
\multicolumn{8}{l}{
\qquad (a) $N = 100$, \; \quad $\tau = 0.1$} \\
\midrule
$\rfrac{1}{4}$ &	   1.5307e-02 &	   5.5755e-02 &	   2.1767e-04 &	   1.2482e-03 &	   1.5375e-02 &	  5.6145e-02 &	 1.91 \\
$\rfrac{1}{8}$ &	   3.4085e-03 &	   1.3741e-02 &	   5.4416e-05 &	   3.1213e-04 &	   3.4255e-03 &	  1.3839e-02 &	 2.01 \\
$\rfrac{1}{16}$ &	   8.2572e-04 &	   3.4190e-03 &	   1.3604e-05 &	   7.8032e-05 &	   8.2996e-04 &	  3.4433e-03 &	 2.04 \\
$\rfrac{1}{32}$ &	   2.0478e-04 &	   8.5343e-04 &	   3.4010e-06 &	   1.9508e-05 &	   2.0584e-04 &	  8.5952e-04 &	 2.04 \\
$\rfrac{1}{64}$ &	   5.1091e-05 &	   2.1324e-04 &	   8.5024e-07 &	   4.8769e-06 &	   5.1356e-05 &	  2.1477e-04 &	 2.04 \\
\midrule
\multicolumn{8}{l}{
\qquad (a) $N = 10^3$, \; \quad $\tau = 0.01$} \\
\midrule
$\rfrac{1}{4}$ &	   4.5193e-03 &	   7.0435e-03 &	   2.1767e-05 &	   1.2204e-04 &	   4.5282e-03 &	  7.0936e-03 &	 1.25 \\
$\rfrac{1}{8}$ &	   6.0419e-04 &	   1.7367e-03 &	   5.4416e-06 &	   3.0530e-05 &	   6.0642e-04 &	  1.7493e-03 &	 1.70 \\
$\rfrac{1}{16}$ &	   1.1831e-04 &	   4.3238e-04 &	   1.3604e-06 &	   7.6332e-06 &	   1.1887e-04 &	  4.3552e-04 &	 1.91 \\
$\rfrac{1}{32}$ &	   2.7537e-05 &	   1.0798e-04 &	   3.4010e-07 &	   1.9084e-06 &	   2.7677e-05 &	  1.0876e-04 &	 1.98 \\
\end{tabular}
\caption{{\it Example \ref{ex:2}}. Convergence of errors and majorants w.r.t. the different choice of spatial mesh sizes and time steps 
(measured relatively w.r.t. $|\!|\!| p^i_h |\!|\!|^2_p$ and $|\!|\!| \vectoru^i_h |\!|\!|^2_{u}$).}
\label{tab:example-2-2-confergence-wrt-mesh-size-and-tau}
\end{table} 

\end{example}

\section*{Conclusion}

We analyze semi-discrete approximations of the Biot poroelastic
problem and deduce guaranteed and fully computable bounds of corresponding errors.
%
The derivation combines estimates for contraction mappings and 
functional a posteriori error majorants for the elliptic problems. 
An extended version of the paper that 
contains an overview of the related works and derivation of the model as well as 
 detailed proofs can be found in \cite{Kundanetal2018}.

\section*{Acknowledgment}
{The work is funded by the SIU Grant CPRU-2015/10040. The second author wants to acknowledge the support
of the Austrian Science Fund (FWF), in particular, the NFN S117-03 project. KK and JMN would like to acknowledge 
Norwegian Research Council project Toppforsk 250223 for funding.  }

\vskip 20pt
\input{biot-cam-paper-kknr.bbl}

\section*{Appendix}

\paragraph{Notation and definitions of spaces}
 %
We use standard Lebesgue space of 
square-measurable functions $\L{2}(\Omega)$ equipped with the norm 
$\|\,v\,\|_{\Omega} := \|\,v\,\|_{L^2(\Omega)} := (v, v)_\Omega^{\rfrac{1}{2}}$
%
for all $u, v \in L^2(\Omega)$. 
Let ${\mathds{M}}^{{d} \times {d}}$ denote the space of real $d$-dimensional tensors. Product for 
vector-valued functions $\vectorv, \vectorw \in \mathbb{R}^{\rm d}$ and tensor-valued functions 
${\boldsymbol \tau}, {\boldsymbol \stress} \in {\mathds{M}}^{{d} \times {d}}$ are defined by the 
relations 
$$(\vectorv, \vectorw)_\Omega := \int_\Omega \vectorv \cdot \vectorw \dx \quad \mbox{and} \quad 
({\boldsymbol \tau}, {\boldsymbol \stress})_\Omega 
:= \int_\Omega {\boldsymbol \tau} : {\boldsymbol \stress} \dx,$$
where $\vectorv \cdot \vectorw := v_i \, w_i$ and 
${\boldsymbol \tau} : {\boldsymbol \stress} := \tau_{ij} \, \stress_{ij}$, respectively.
Next, ${\mathbb{\bf A}} (x) \in {\mathds{M}}^{d \times d}$, 
$x \in \Omega$ denotes a symmetric uniformly positive defined matrix that satisfies  
$0 < \underline{\lambda} \leq \lambda(x) \leq \overline{\lambda} \leq +\infty$,  
$\underline{\lambda}, \overline{\lambda} \in \mathbb{R}$, 
with uniformly bounded eigenvalues $\lambda(x)$.
Then, for the product $({\bf u}, \boldsymbol{v})_{{\mathbb{\bf A}}} := ({\mathbb{\bf A}} {\bf u}, \boldsymbol{v})$, 
we have 
%
$$\underline{\lambda} \, \| \boldsymbol{v} \| \leq \|  \boldsymbol{v} \|_{{\mathbb{\bf A}}} \leq  \overline{\lambda} \, \| \boldsymbol{v} \| 
\quad \mbox{and} \quad 
({\bf u}, \boldsymbol{v}) \leq \| \cdot \|_{{\mathbb{\bf A}}} \, \| \cdot \|_{{\mathbb{\bf A}^{-1}}}, \quad 
\forall {\bf u}, \boldsymbol{v} \in [L^2(\Omega)]^d.$$
We use the standard notation for the Sobolev space of vector-valued functions having square-summable derivatives
\begin{equation*}
H^1(\Omega) := \big\{\, v \in L^2(\Omega)\,\mid\, \nabla\,v \in [L^2(\Omega)]^d \,\big\},
\end{equation*}
equipped with the norm $\|\,v\,\|_{H^1(\Omega)} := \big(\|\, v\, \|^2_{\Omega} + |\, v\, |_{\Omega}^2\big)^{\rfrac{1}{2}}.$
Also, we use the  semi-norm $|\, v\, |_{\Omega} := |\, v\, |_{H^1(\Omega)} := \|\,\nabla\,v\,\|_{\Omega}$, 
and for the the vector-valued functions with square-summable divergence introduce the Hilbert space
%
%
%
\begin{equation*}
H(\Omega, \dvrg) 
:= \big\{\, \vectorv \in [L^2(\Omega)]^d\,\mid\, \dvrg \vectorv \in L^2(\Omega) \,\big\},
\end{equation*}
endowed with the norm
$\|\, \vectorv\, \|^2_{H(\Omega, \dvrg)} 
:= \|\, \vectorv \, \|^2_{[L^2(\Omega)]^d} + \| \, \dvrg \vectorv \, \|^2_{L^2(\Omega)}.$

Let $\Sigma$ be a part of the boundary  such that $\meas_{d - 1} \Sigma > 0$  (in particular, may coincide with 
$\partial \Omega$). For functions in $H^1_{0,\Sigma}(\Omega) := \big\{ \, v \in H^1(\Omega) \,\mid\, v\!\mid_{\Sigma} = 0 \, \big\}$, 
the Friedrichs'-type inequality reads:
%
\begin{equation}
\|v\|_{\Omega} \leq \CF |v|_{\Omega}, \quad 
\forall v \in H^1_{0, \Sigma}(\Omega).
\label{eq:Poincare}
\end{equation}
%
%
The corresponding trace operator $\gamma: H^1(\Omega) \rightarrow H^{\tfrac{1}{2}} (\Omega)$ is bounded and
satisfies the estimate
\begin{equation}
v |_{\Sigma} := \gamma v, \quad 
\|\, v \,\|_{\Sigma} 
\leq C^{\rm tr}_{\Sigma\Omega} \,\| \, v \, \|_{H^1(\Omega)}, \quad 
\forall v \in H^1(\Omega),
\label{eq:trace}
\end{equation}
where $\|\,v\,\|_{\Sigma}$ is the norm of $L^2(\Sigma)$.

$C_{\rm K}$ denotes the constant in the Korn inequality
\begin{equation}
\| \vectorw \|_{[H^1(\Omega)]^d} 
\leq C_{\rm K}
\|\strain(\vectorw)\|_{[L^2(\Omega)]^{{d} \times {d}}}, \quad 
\forall \vectorw \in [H^1(\Omega)]^d,
\label{eq:korn}
\end{equation}
%
Also, we use the inequality
\begin{equation}
\| \, \dvrg \, \vectorw \, \| 
= \| \, {\rm tr} \, \strain(\vectorw) \, \| 
= \| \, \tensorI : \strain(\vectorw) \, \| 
\leq \sqrt{d}\, \| \, \strain(\vectorw) \, \|, \quad 
\forall \vectorw \in [H^1(\Omega)]^d, 
\label{eq:inequality-divu-strain}
\end{equation}
where ${\mathds{I}} \in {\mathds{M}}^{d \times d}$ is the unit tensor of ${{\mathds{M}}}^{{d} \times {d}}$, 
and $\strain(\vectorw) \in {\mathds{M}}^{{d} \times {d}}$ denotes the symmetric part of $\nabla \vectorw$.

%
Next, let $Q := \Omega \times (0, T)$ denote a space-time 
cylinder (with given time-interval $(0, T)$, 
$0 < T < +\infty$), and let $\Sigma = \partial \Omega \times (0, T)$ be a lateral surface of $Q$, whereas 
$\Sigma_0 := \partial \Omega \times \{0\}$ and $\Sigma_T := \partial \Omega \times \{T\}$ define the bottom 
and the top parts of the mantel (so that $\partial Q = \Sigma \cup \Sigma_0 \cup \Sigma_T$). 
Consider functions defined in $(0, T)$ with values in a functional space $X$
(cf. ~\cite{Ladyzhenskaya1985, Ladyzhenskayaetall1967, Zeidler1990A}). Let $\|\cdot\|_X$ denote the norm 
in $X$, then for $r = 2$, we define the Bochner space 
\begin{equation*}
L^2(0, T; X) := 
\Big\{\, f \mbox{ measurable in } [ \, 0, T \,]\;
      \Big| \; \int_0^T \|f(t)\|_X^2 \dt < \infty \, \Big\}, 
\end{equation*}
and respective norm
$\|f\|_{L^2(0, T; X)} := \Big( \int_0^T \|f(t)\|_X^2 \dt \Big)^{\rfrac{1}{2}}.$
%
It is a Hilbert space space if $X$ is a Hilbert space. 
Throughout the paper, we also use the spaces  
\begin{equation}
H^1(0,T; X) := \big\{\, f \in L^2(0,T; X)\,\;|\;  {\partial_t} f \in L^2(0,T; X) \,\big\}
\label{eq:bochner-H1-a-b-X}
\end{equation}
%
%
equipped with norm
$\| u \|_{H^1(0,T; X)} := \Big( \int_0^T \big(\| \partial_t f(t)\|_X^2 + \|f(t)\|_X^2\big) \dt \Big)^{\rfrac{1}{2}}.$
%
%

We assume that $\Tau_h$ is a regular mesh satisfying angle condition defined on $\Omega$. 
Then, the corresponding discretization spaces
with the Lagrangian finite elements of order $0$ or $1$ are defined as
\begin{alignat}{2}
\mathds{P}_0 & := \{ v_h \in L^2 (\Omega) \, | \, \forall T \in \Tau_h, \, v_h  |_T \in \mathds{P}_0 \}, \quad 
\mathds{P}_1 & := \{ v_h \in H^1 (\Omega) \, | \, \forall T \in \Tau_h, \,  v_h |_T \in \mathds{P}_1 \}, 
\label{eq:p0-p1}
\end{alignat}
where $\mathds{P}_k$ denotes the space of polynomials of the order $k \in \mathds{N} \cup 0$.
The Raviart-thomas elements of the lowest and first order are denoted by
\begin{align*}
\mathds{RT}_0 := \{& \boldsymbol{y}_h \in H(\dvrg, \Omega): 
\, \forall \, T \in \Tau_h, \,  \boldsymbol{y}_h |_T = \boldsymbol{a} + b \, \boldsymbol{x}, 
\boldsymbol{a} \in \mathds{R}^d, b \in \mathds{R}\, \}, \\
\mathds{RT}_1 := \{&\boldsymbol{y}_h \in H(\dvrg, \Omega): 
\, \forall \, T \in \Tau_h,   \boldsymbol{y}_h(\boldsymbol{x}) |_T = \boldsymbol{q}(\boldsymbol{x}) 
+ \boldsymbol{x} \, r(\boldsymbol{x}) , \boldsymbol{q} \in [\mathds{P}_1]^d, r \in \mathds{P}_1 \},
\end{align*}
respectively. Finally, the table below presents notation used on the paper for the physical quantities.

\begin{table}[htbp]
\footnotesize
\begin{center}
\begin{tabular}{r c p{10cm} }
\toprule
$\stress_{\rm por}$ & $$ & poroelastic Cauchy stress (total stress) tensor \\
$\vectoru$ & $$ & displacement of the solid\\
$p$ & $ $ & fluid pressure \\
$\stress$ & $ $ & linear elastic (effective) stress tensor \\  
$\strain (\vectoru)$ & $ $ & strain tensor \\  
$\lambda, \mu$ & $ $ & Lam\'e paramteres \\  
$\boldsymbol{f}$ & $ $ & volumetric body force \\ 
$\vectorw$  & $ $ & Darcy velocity \\  
$\mu_f$ & $$ & fluid viscosity \\
$\tensorK$ & $$ & permeability tensor\\  
$g$ & $$ & gravitation constant\\ 
$\rho_f$ & $$ & fluid phase density \\
$\alpha$ & $$ & Biot-Willis coefficient \\
$\beta = \tfrac{1}{M} + c_f \varphi_0$ & $$ & storage coefficient \\
$M$ & $$ & Biot constant \\
$c_f$ & $$ &  fluid compressibility \\
$\varphi_0$ & $$ & initial porosity \\ 
\bottomrule
\end{tabular}
\end{center}
\vskip -10pt
\caption{Table of notation}
\label{tab:table-of-notation}
\end{table}

\paragraph{Contraction theorem} 
An essential part of our analysis is based on the contraction theorem, which was proven in \cite{MikelicWheeler2013} 
for the poroelastic problem in question. However, Lemma \ref{lem:pressure-full-norm-majorant} uses the contraction 
of the sequence $\{ \delta(\eta - \eta_h) \}^i$, $\forall \eta \in W_h$, which is proven blow.
\vskip 8pt
\begin{theorem}
\label{th:theorem-contraction-appendix}
With  $\gamma = \tfrac{\alpha}{\sqrt{\lambda}}$ and $L = \tfrac{\alpha^2}{2 \,\lambda}$, the sequence 
$\{ \delta(\eta - \eta_h) \}^i \in W_h$, where $\eta^i \in W$ is generated by the fixed-stress split iterative scheme defined in 
\eqref{eq:flow-fully-discrete}--\eqref{eq:mechanics-fully-discrete} and $\eta^i_h \in W_h$ is discretization of the latter 
sequence, is a contraction given by
\begin{equation}
\| \strain (\delta (\vectoru -\vectoru_h)^{i}) \|^2_{2\mu}
+ q \, \| \nabla \delta (p - p_h)^{i} \|^2_{\tensorK_\tau} 
+ \| \delta (\eta - \eta_h)^{i} \|^2
\leq q^2 \| \delta (\eta - \eta_h)^{i-1} \|^2, \quad q = \tfrac{L}{\beta + L}.
\label{eq:contration-general-II}
\end{equation}
\end{theorem}
\ProofBegin
Consider the difference between $(i\minus1)$th and $i$th iterations in \eqref{eq:flow-fully-discrete} and 
\eqref{eq:mechanics-fully-discrete}. Assuming that $\delta p^{i}  = p^{i} - p^{i-1}$ and 
$\delta \vectoru^{i} = \vectoru^{i} - \vectoru^{i-1}$, as well as particular chosen 
$w = w_h \in W_{0h} \subset W_0$ and 
$\vectorv = \vectorv_h \in \boldsymbol{V_{0h}} \subset \boldsymbol{V_{0}}$ with conforming Galerkin 
disctretization 
spaces $W_{0h}$ and $\boldsymbol{V_{0h}}$ of $W_0$ and $\boldsymbol{V_{0}}$, respectively. 
Thus, we obtain
\begin{alignat}{2}
(\tensorK_{\tau} \nabla \delta {p^{i}}, \nabla w_h) +  (\beta + L) (\delta {p^{i}}, w_h)
  & = (- \gamma \,\delta  {\eta}^{i-1} , w_h) , \qquad \qquad \;\, \forall \, w \in W_{0h}, 
\label{eq:flow-fully-discrete-diff} \\
\big(2 \, \mu \, \strain (\delta {\vectoru^{i}}), \strain (\vectorv_h) \big) + (\lambda \, \dvrg \delta {\vectoru^{i}}, \dvrg \vectorv_h)
  & = (- \alpha  \nabla \delta p^{i}, \vectorv_h), \qquad \qquad \,\,\forall \, \vectorv_h \in \boldsymbol{V_{0h}},
\label{eq:mechanics-fully-discrete-diff} 
\end{alignat}
For Galerkin approximations $(\vectoru, p)^{i}_h \in \boldsymbol{V_{0h}} \times W_{0h}$, the system above can be rewritten as 
\begin{alignat}{2}
(\tensorK_{\tau} \nabla \delta {p^{i}_h}, \nabla w_h) + 
(\beta + L) (\delta p^{i}_h, w_h)
	& = (- \gamma \,\delta  {\eta}^{i-1}_h , w_h) , \qquad \qquad \;\, \forall \, w_h \in W_{0h}, 
	\label{eq:flow-fully-discrete-diff-h} \\
\big(2 \, \mu \, \strain (\delta {\vectoru^{i}_h}), \strain (\vectorv_h) \big) 
	+ (\lambda \, \dvrg \delta {\vectoru^{i}_h}, \dvrg \vectorv_h)
  & = (- \alpha  \nabla \delta p^{i}_h, \vectorv_h), \qquad \qquad \,\,\forall \, \vectorv \in \boldsymbol{V_{0h}}.
\label{eq:mechanics-fully-discrete-diff-h} 
\end{alignat}
Difference of \eqref{eq:flow-fully-discrete-diff}--\eqref{eq:mechanics-fully-discrete-diff} and 
\eqref{eq:flow-fully-discrete-diff-h}--\eqref{eq:mechanics-fully-discrete-diff-h}, 
substitution of $w_h = \delta (p - p_h)^{i}$ in flow part and 
$\vectorv_h= \delta (\vectoru - \vectoru_h)^{i}$ in the    corresponding mechanics part yield
\begin{alignat}{2}
\| \nabla \delta (p - p_h)^{i} \|^2_{\tensorK_\tau}  + (\beta + L) \, \| \delta (p - p_h)^{i} \|^2
	& = - \gamma \,(\delta  (\eta - \eta_h)^{i-1}, \delta (p - p_h)^{i}) , \;\, \forall \, w_h \in W_{0h}, 
	\label{eq:flow-fully-discrete-diff-h-diff} \\
\| \strain (\delta (\vectoru - \vectoru_h)^{i}) \|^2_{2\mu} + \| \dvrg \delta (\vectoru - \vectoru_h)^{i}\|^2_{\lambda} 
& = - \alpha \, (\delta (p - p_h)^{i}, \dvrg \delta (\vectoru - \vectoru_h)^{i}), \,\,\forall \, \vectorv \in \boldsymbol{V_{0h}}.
\label{eq:mechanics-fully-discrete-diff-h-diff} 
\end{alignat}
Application of the Young inequality in \eqref{eq:flow-fully-discrete-diff-h-diff} provide the relation
\begin{alignat}{2}
	(\beta + L) \, \| \delta (p - p_h)^{i} \|^2
	+ \| \nabla \delta (p - p_h)^{i} \|^2_{\tensorK_\tau} 
	\leq 
	\tfrac{\epsilon}{2} \, \| \delta (p - p_h)^{i} \|^2 + 
	\tfrac{\gamma^2}{2 \epsilon} \, \| \delta (\eta - \eta_h)^{i-1} \|^2, \quad \epsilon > 0.
	\label{eq:flow-iterative-3}
\end{alignat}
Regrouping similar terms in \eqref{eq:flow-iterative-3} implies
\begin{alignat*}{2}
	(\beta + L - \tfrac{\epsilon}{2}) \, \| \delta (p - p_h)^{i} \|^2 
	+ \| \nabla \delta (p - p_h)^{i} \|^2_{\tensorK_\tau}
	\leq 
	\tfrac{\gamma^2}{2 \epsilon} \, \| \delta ({\eta - \eta_h})^{i-1} \|^2. 
\end{alignat*}
%
Substitution of the optimal $\epsilon = \beta + L$, obtained from the minimization problem
$\min\limits_{\eps > 0} \big(2 \, \epsilon \,(\beta + L - \tfrac{\epsilon}{2})\big)^{-1},$
yields 
\begin{equation}
(\beta + L) \, \| \delta (p - p_h)^{i}  \|^2 
+ 2\, \| \nabla \delta (p - p_h)^{i}  \|^2_{\tensorK_\tau} 
\leq \tfrac{\gamma^2}{\beta + L} \, \| \delta (\eta - \eta_h)^{i-1} \|^2. 
\label{eq:flow-iterative-5}
\end{equation}
%
%
By summing \eqref{eq:flow-iterative-5}, multiplied by free parameter $c_0 > 0$, 
and \eqref{eq:mechanics-fully-discrete-diff-h-diff}, we arrive at the following inequality
%
\begin{alignat}{2}
  \Big\{ c_0 \, (\beta + L) \, \| \delta (p - p_h)^{i} \|^2
	& + \| \dvrg \delta (\vectoru - \vectoru_h)^{i} \|^2_\lambda 
	    - \alpha \, (\delta (p - p_h)^{i} , \dvrg \delta (\vectoru - \vectoru_h)^{i}) \Big\} \nonumber\\
	& + \| \strain (\delta (\vectoru - \vectoru_h)^{i}) \|^2_{2\mu}
	+ 2 c_0 \, \| \nabla \delta (p - p_h)^{i} \|^2_{\tensorK_\tau} 
	\leq c_0 \, \tfrac{\gamma^2}{\beta + L} \| \delta (\eta - \eta_h)^{i-1} \|^2. 
	\label{eq:sum-flow-and-mechanics-iterative}
\end{alignat}
%
Let us determine the values of parameters $c_0$, $\gamma$, and $L$ such that the terms in 
the left-hand side of \eqref{eq:sum-flow-and-mechanics-iterative} are positive and contraction in 
$\| \delta (\eta - \eta_h)^{i-1}\|^2$ is achieved. 
It follows from 
$ \delta (\eta - \eta_h)^{i} 
= \tfrac{\alpha}{\gamma} \, \dvrg \delta (\vectoru - \vectoru_h)^{i} - \tfrac{L}{\gamma} \, \delta (p - p_h)^{i}$, 
that 
\begin{equation}
\| \delta (\eta - \eta_h)^{i}\|^2 
= \tfrac{\alpha^2}{\gamma^2} \, \| \dvrg \delta (\vectoru - \vectoru_h)^{i} \|^2 
+ \tfrac{L^2}{\gamma^2} \, \| \delta (p - p_h)^{i} \|^2
- \tfrac{2 \alpha \, L}{\gamma^2} \, (\dvrg \delta (\vectoru - \vectoru_h)^{i}, \delta (p - p_h)^{i}).
\label{eq:relation-sigma-div-u-p-squared}
\end{equation}
%
Comparing \eqref{eq:relation-sigma-div-u-p-squared} and \eqref{eq:sum-flow-and-mechanics-iterative}, 
we arrive at the following condition on the free parameters:
$$
\begin{cases}
\; \tfrac{\alpha^2}{\gamma^2} \leq \lambda, \label{eq:coeff-dvrg-u}\\
\; \tfrac{L^2}{\gamma^2} \leq c_0 \, (\beta + L), \label{eq:coeff-p}\\
\; \tfrac{2 \alpha \, L}{\gamma^2} = \alpha, \label{eq:coeff-cross-term}
\end{cases}
%
\quad \mbox{which yields} \qquad
%
\begin{cases}
\; L \geq \tfrac{\alpha^2}{2 \,\lambda},\\
\; c_0 \geq \tfrac{L}{2\, (\beta + L)}, \\
\; \gamma^2 = 2 L. 
\end{cases}
$$
%
Then, the contraction rate $q = c_0 \, \tfrac{\gamma^2}{\beta + L}$ is monotone 
w.r.t. to $L$ and attains its minimum at
$$L = \tfrac{\alpha^2}{2 \,\lambda} 
\quad \mbox{and} 
\quad c_0 = \tfrac{L}{2\, (\beta + L)}.$$
By using condition $\gamma^2 = 2 L$, we obtain
\begin{equation}
\| \strain (\delta (\vectoru -\vectoru_h)^{i}) \|^2_{2\mu}
+ q \, \| \nabla \delta (p - p_h)^{i} \|^2_{\tensorK_\tau} 
+ \| \delta (\eta - \eta_h)^{i} \|^2
\leq q^2 \| \delta (\eta - \eta_h)^{i-1} \|^2
\end{equation}
with
\begin{equation*}
q = \tfrac{L}{\beta + L} \quad \mbox{and} \quad L = \tfrac{\alpha^2}{2 \,\lambda}.
\end{equation*}
\ProofEnd

\end{document}

%% file: customized_packages_and_macros.tex
\usepackage{amsmath}
\usepackage[position=top]{subfig} 
\usepackage{amsfonts} 							
\usepackage{dsfont}
\usepackage{algorithm}              
\usepackage{algorithmic}
\usepackage{booktabs} 							
\usepackage{multirow}								
\usepackage{setspace}
\usepackage{mathrsfs}
\usepackage{amssymb}
\usepackage{tikz}
\usepackage{xcolor}
\usepackage{comment}

\usepackage{framed}
\usepackage{afterpage}
\usepackage{capt-of}
\numberwithin{equation}{section}
\usepackage{enumerate}

\definecolor{formalshade}{rgb}{0.95,0.95,1}

\usepackage{geometry}
\usepackage{adjustbox}
\usepackage{sidecap}

\definecolor{gray(x11gray)}{rgb}{0.75, 0.75, 0.75}
\newcommand{\mnote}[1]{\!^\textrm{\scriptsize\color{gray(x11gray)}#1}}

\newcommand*\rfrac[2]{{}^{#1}\!/_{#2}}
\newcommand{\minus}{\scalebox{0.8}[1.0]{-}}

\newcommand {\eps}   {\varepsilon}
\newcommand {\Tau}   {\mathcal{T}}

\newcommand {\strain}   {\boldsymbol {\varepsilon}}
\newcommand {\stress}   {\boldsymbol {\sigma}}
\newcommand {\vectorv}  {\boldsymbol {v}}
\newcommand {\vectoru}  {\boldsymbol {u}}
\newcommand {\vectorw}  {\boldsymbol {w}}
\newcommand {\vectorn}  {\boldsymbol {n}}
\newcommand {\vectorz}  {\boldsymbol {z}}
\newcommand {\vectortau}  {\boldsymbol {\tau}}

\newcommand {\vectort}  {\boldsymbol {t}}

\newcommand {\tensorI}  {{\mathds{I}}}
\newcommand {\tensorL}  {\mathds L}
\newcommand {\tensorK}  {\mathds {K}}

\usepackage{array}

\newcommand {\CF}    {{C^{\rm F}_{\Gamma}}}

\newcommand {\Rtwo}  {{\mathds{R}}^2}

\newcommand {\Ieff}  {{I_{\rm eff}}}

\def \NormA#1  {{\mid\!\mid\!\mid #1 \mid\!\mid\!\mid}^2_A }   
\def \NormsqC#1  {{\| #1 \|}^2_{\mathds{C}} }   
\def \NormC#1  {{\mid\!\mid\!\mid #1 \mid\!\mid\!\mid}_{\mathds{C}} }   

\def \NormAinverse#1  { {\mid\!\mid\!\mid #1 \mid\!\mid\!\mid}^2_{A^{-1}} }   
\def \Normenergy#1  {\mid\!\mid\!\mid #1 \mid\!\mid\!\mid}   

\def \NormQT#1 {{\mid\!\mid\!\mid #1 \mid\!\mid\!\mid}^2_{Q_T}}   
\def \NormQO#1 {{\mid\!\mid\!\mid #1 \mid\!\mid\!\mid}^2_{Q^0}}   
\def \NormQk#1 {{\mid\!\mid\!\mid #1 \mid\!\mid\!\mid}^2_{Q^k}}   
\def \Normf#1  {\Big \lceil #1 \Big\rceil_\Omega}
\def \dvrg     {\mathrm{div}}	

\def \dt       {\mathrm{\:d}t}
\def \dx       {\mathrm{\:d}x}

\def \Normt#1  {\mid\!\mid\!\mid #1 \mid\!\mid\!\mid}   

\def\L#1{L^{#1}}
\def\H#1{H^{#1}}

\def\H#1{H^{#1}}

\newcommand {\meas} {\mathrm{meas}}

\def \Mean#1#2 {{ \Big \{ #1 \Big\} }_{#2}}
\def \smallMean#1#2 {{ \big \{  #1 \big\} }_{#2}}

\def\ProofBegin{\noindent{\bf Proof:} \:}
\def\ProofEnd{{\hfill $\square$} \vskip 2pt
\vskip 10pt}

\newtheorem{theorem}{Theorem}

\newtheorem{lemma}{Lemma}
\newtheorem{remark}{Remark}
\newtheorem{corollary}{Corollary}
\newtheorem{example}{Example}

\usepackage{color}
\definecolor{deepblue}{rgb}{0,0,0.5}
\definecolor{deepred}{rgb}{0.6,0,0}
\definecolor{deepgreen}{rgb}{0,0.5,0}
\usepackage{listings}

\newcommand\pythonstyle{
\lstset{
		language=Python,
		basicstyle=\footnotesize\ttfamily,
		otherkeywords={self},             
		keywordstyle=\bfseries\color{deepblue}\bfseries,
		commentstyle=\itshape\color{purple!40!black},
		stringstyle=\color{deepgreen},
		frame=single,                         
		showstringspaces=false,            %
		belowcaptionskip=.75\baselineskip
}}
\lstnewenvironment{python}[1][]
{
\pythonstyle
\lstset{#1}
}
{}

\newcommand\pythoninline[1]{{\pythonstyle\lstinline!#1!}}

\usepackage{colortbl}
\definecolor{dollarbill}{rgb}{0.52, 0.73, 0.4}
\definecolor{grannysmithapple}{rgb}{0.66, 0.89, 0.63}
\definecolor{junbud}{rgb}{0.7, 0.93, 0.36}
\definecolor{olivine}{rgb}{0.6, 0.73, 0.45}
\definecolor{pistachio}{rgb}{0.58, 0.77, 0.45}
\definecolor{teagreen}{rgb}{0.82, 0.94, 0.75}
\definecolor{yellow-green}{rgb}{0.6, 0.8, 0.2}
\definecolor{babypink}{rgb}{0.96, 0.76, 0.76}
\definecolor{beaublue}{rgb}{0.74, 0.83, 0.9}
\definecolor{corn}{rgb}{0.98, 0.93, 0.36}
\definecolor{arylideyellow}{rgb}{0.91, 0.84, 0.42}
\definecolor{cadmiumgreen}{rgb}{0.0, 0.42, 0.24}
\definecolor{darkspringgreen}{rgb}{0.09, 0.45, 0.27}
\definecolor{dartmouthgreen}{rgb}{0.05, 0.5, 0.06}
\definecolor{forestgreen(web)}{rgb}{0.13, 0.55, 0.13}
\definecolor{blue(pigment)}{rgb}{0.2, 0.2, 0.6}
\definecolor{junebud}{rgb}{0.74, 0.85, 0.34}

\definecolor{DeepSkyBlue}{rgb}{0,0.75,1}
\definecolor{DeepSkyBlue1}{rgb}{0,0.75,1}
\definecolor{DeepSkyBlue2}{rgb}{0,0.7,0.93}
\definecolor{DeepSkyBlue3}{rgb}{0,0.6,0.8}
\definecolor{DeepSkyBlue4}{rgb}{0,0.41,0.54}
\definecolor{DarkSeaGeen4}{rgb}{0.,0.45,0.}
\definecolor{darkcandyapplered}{rgb}{0.64, 0.0, 0.0}
\definecolor{byzantium}{rgb}{0.58, 0.0, 0.83}
\definecolor{lightblue}{rgb}{0.68, 0.85, 0.9}
\definecolor{darkblue}{rgb}{0.,0.,0.45}
\definecolor{verydarkblue}{rgb}{25,25,112}
\definecolor{darkspring}{rgb}{0.09, 0.45, 0.27}
\definecolor{darkcoral}{rgb}{0.8, 0.36, 0.27}
\definecolor{crimsonglory}{rgb}{0.75, 0.0, 0.2}
\definecolor{darkblue}{rgb}{0.0, 0.0, 0.55}
\definecolor{darkelectricblue}{rgb}{0.33, 0.41, 0.47}

\definecolor{dollarbill}{rgb}{0.52, 0.73, 0.4}
\definecolor{grannysmithapple}{rgb}{0.66, 0.89, 0.63}
\definecolor{junbud}{rgb}{0.7, 0.93, 0.36}
\definecolor{olivine}{rgb}{0.6, 0.73, 0.45}
\definecolor{pistachio}{rgb}{0.58, 0.77, 0.45}
\definecolor{teagreen}{rgb}{0.82, 0.94, 0.75}
\definecolor{yellow-green}{rgb}{0.6, 0.8, 0.2}
\definecolor{babypink}{rgb}{0.96, 0.76, 0.76}
\definecolor{beaublue}{rgb}{0.74, 0.83, 0.9}
\definecolor{corn}{rgb}{0.98, 0.93, 0.36}
\definecolor{arylideyellow}{rgb}{0.91, 0.84, 0.42}
\definecolor{cadmiumgreen}{rgb}{0.0, 0.42, 0.24}
\definecolor{applegreen}{rgb}{0.29, 0.33, 0.13}
\definecolor{armygreen}{rgb}{0.29, 0.33, 0.13}
\definecolor{dartmouthgreen}{rgb}{0.05, 0.5, 0.06}
\definecolor{forestgreen(web)}{rgb}{0.13, 0.55, 0.13}
\definecolor{blue(pigment)}{rgb}{0.2, 0.2, 0.6}
\definecolor{junebud}{rgb}{0.74, 0.85, 0.34}

\definecolor{byzantium}{rgb}{0.44, 0.16, 0.39}
\definecolor{darkbyzantium}{rgb}{0.36, 0.22, 0.33}
\definecolor{DeepSkyBlue}{rgb}{0,0.75,1}
\definecolor{DeepSkyBlue1}{rgb}{0,0.75,1}
\definecolor{DeepSkyBlue2}{rgb}{0,0.7,0.93}
\definecolor{DeepSkyBlue3}{rgb}{0,0.6,0.8}
\definecolor{DeepSkyBlue4}{rgb}{0,0.41,0.54}
\definecolor{DarkSeaGeen4}{rgb}{0.,0.45,0.}
\definecolor{darkcandyapplered}{rgb}{0.64, 0.0, 0.0}
\definecolor{lightblue}{rgb}{0.,0.,0.9}
\definecolor{darkblue}{rgb}{0.,0.,0.45}
\definecolor{verydarkblue}{rgb}{25,25,112}
\definecolor{darkblue}{rgb}{0.0, 0.0, 0.55}
\definecolor{darkelectricblue}{rgb}{0.33, 0.41, 0.47}
\definecolor{gainsboro}{rgb}{0.86, 0.86, 0.86}
\definecolor{lightgray}{rgb}{0.83, 0.83, 0.83}
\definecolor{silver}{rgb}{0.75, 0.75, 0.75}
\definecolor{ashgrey}{rgb}{0.7, 0.75, 0.71}
\definecolor{battleshipgrey}{rgb}{0.52, 0.52, 0.51}
\definecolor{prussianblue}{rgb}{0.0, 0.19, 0.33}
\definecolor{uclablue}{rgb}{0.33, 0.41, 0.58}

\definecolor{midnightblue}{rgb}{0.1, 0.1, 0.44}
\definecolor{upforestgreen}{rgb}{0.0, 0.27, 0.13}
\definecolor{midnightgreen(eaglegreen)}{rgb}{0.0, 0.29, 0.33}

\renewcommand{\(}{\begin{columns}}
\renewcommand{\)}{\end{columns}}
\newcommand{\<}[1]{\begin{column}{#1}}
\renewcommand{\>}{\end{column}}

\tikzstyle{every picture}+=[remember picture]

\usetikzlibrary{%
  arrows,%
  shapes.misc,
  shapes.arrows,%
  chains,%
  matrix,%
  positioning,
  scopes,%
  decorations.pathmorphing,
  shadows,%
  backgrounds,
  fit,positioning,shapes.symbols,chains
}
\usetikzlibrary{shapes.geometric,calc}
\usetikzlibrary{positioning, fit, scopes}
\tikzset{label/.style={align=left,text width=5.5cm,inner sep=1mm,outer sep=0mm,font=\scriptsize}}

%% file: biot-cam-paper-kknr.bbl
\def\cprime{$'$} \def\cprime{$'$} \def\cprime{$'$}